\documentclass[fleqn,twoside]{article}
\usepackage{espcrc1}

\newcommand{\AmS}{{\protect\the\textfont2
  A\kern-.1667em\lower.5ex\hbox{M}\kern-.125emS}}
\title{TOWARDS DESIGN OF SYSTEM HIERARCHY (RESEARCH SURVEY)}

\author{Mark Sh. Levin
\thanks{
 Mark Sh. Levin:~
 http://www.mslevin.iitp.ru;~
 email: mslevin@acm.org
  }
  }

\begin{document}

\maketitle

\begin{abstract}
  The paper addresses design/building frameworks for
 some kinds of tree-like and hierarchical structures of systems.
 The following approaches are examined:
  (1) expert-based procedure,
 (2)  hierarchical clustering;
 (3) spanning problems
  (e.g., minimum spanning tree,
  minimum Steiner tree,
  maximum leaf spanning tree problem;
 (4) design of organizational ``optimal'' hierarchies;
 (5) design of multi-layer
 (e.g., three-layer)
  k-connected network;
 (6) modification of hierarchies:
 (i) modification of tree via condensing of neighbor nodes,
 (ii) hotlink assignment,
 (iii) transformation of tree into Steiner tree,
 (iv) restructuring as modification of an initial structural
 solution into a solution that is the most close to a goal solution while taking into account
 a cost of the modification.
  Combinatorial optimization problems are considered
 as basic ones (e.g.,
 classification,
 knapsack problem,
 multiple choice problem,
 assignment problem).
 Some numerical examples illustrate the suggested problems and
 solving frameworks.

~~~~~~~~~~~

 {\it Keywords:}~
 hierarchy,
  multi-layer structure,
  tree, networks, ontology,
  information systems,
 solving strategy,
 heuristics,
 combinatorial optimization,
 knapsack,
 multiple choice problem,
 spanning tree,
 Steiner tree,
 hotlink assignment,
 multicriteria ranking

\vspace{1pc}

\end{abstract}

%%%%%%%%%%%%%%%%%%%%%%%%%%%%%%%%%%%%%%%%%%%%%%%%%%%%
\tableofcontents

\newcounter{cms}
\setlength{\unitlength}{1mm}

\section{INTRODUCTION}

 Hierarchies play a crucial role as a very useful model for complex systems
  in all domains:
 (1) hierarchies leads to a decomposition (i.e., partitioning)
 of the system representation
 and corresponding system problems
 (e.g., design, improvement, maintenance);
 (2) a hierarchy is often an excellent basis
 for fast algorithms and/or solving procedures design
 (i.e., an analogue of linear or convex functions
 in continuous mathematics);
 (3) hierarchical approach is the best one for
  structures of information/knowledge (e.g., ontology);
 (4) hierarchical approach is the best one for
  problem solving strategies (e.g., decision trees);
 (5) hierarchies are
 are very understandable for humans and
 can be used as a basis of easy learning/teaching interactive procedures.
 Evidently, many years various methods have been used
 to design the hierarchical structures
 (e.g., \cite{abts12},  \cite{ait08},
 \cite{bley04},
% \cite{brun03},
 \cite{corcho03}, \cite{eng88}, \cite{gar79},
 \cite{gehrke09},
 \cite{good88}, \cite{good90},
  \cite{gubko06}, \cite{haimes81},  \cite{haimes90}, \cite{holler04},
  \cite{johnson07},
 \cite{knuth68},
 \cite{kuz05}, \cite{mag95},
 \cite{mar01}, \cite{mirkin96}, \cite{molnar12}, \cite{obr10}, \cite{voronin02},
 \cite{waller76}, \cite{yel76}).

 In the article,
 the following two basic problems are examined:

 {\bf I.} Design approaches to design the hierarchies
 (including expert-based procedure, hierarchical clustering and spanning trees).

  {\bf II.} Schemes for transformation of hierarchies
  (e.g., hotlink assignment, transformation of tree into Steiner tree,
  restructuring).

 Mainly, the considered problems are briefly described via a framework:
 (i) engineering description,
 (ii) problem formulation(s),
 (iii) some problem versions,
 (iv) basic solving schemes,
 (v) new prospective problem versions
 (e.g., multicriteria problems,
 problems under uncertainty).
 Numerical examples illustrate the described approaches.

\section{BASIC HIERARCHIES AND SOME DESIGN APPROACHES}

%%%%%%%%%%%%%%%%%%%%%%%%%%%%%%%%%%%%%%%%%%%%%%%%%%%%%%%%%%%%%%%%%%%%

\subsection{Basic Types of Hierarchies}

 Generally, it is reasonable to point out some
 basic types  of hierarchies
 (e.g., \cite{gar79},
 \cite{haimes90}, \cite{knuth68}):

 (1) various kinds of trees (e.g., Fig. 1, Fig. 2, Fig. 3)
 (e.g.,  \cite{furnas94},\cite{gar79},\cite{knuth68});

 (2) organic hierarchy
 (i.e., with organic interconnection among children-vertices,
 Fig. 4)
  \cite{conant74};

 (3) ``basic'' hierarchy as
  a tree with additional edges (Fig. 5)
   (e.g., \cite{lev89});

 (4) ``morphological hierarchy''
  (e.g., \cite{lev98},\cite{lev06},\cite{lev11agg})
  (Fig. 6);

 (5) multi-layer structures
 (e.g.,
 multi-layer networks, hierarchical networks)
 (Fig. 7)
  (e.g., \cite{abts12},\cite{bel97},\cite{bley04}, \cite{holler04},\cite{lev11ADES},\cite{obr10}).

 Here it is reasonable to point out some important research
 directions in modeling of various multi-layer graphs/networks, for example:
  (a) hypergraphs
 (e.g.,
   \cite{berge73},\cite{berge89})
 and
  hypernetworks
  (e.g.,
  \cite{ha10},\cite{johnson07});
  (b) multi-layer social networks
 (e.g., \cite{kazienko10},\cite{magnani11});
  (c) multi-stratum networks
 (e.g., \cite{magnani12});
 (d) multi-layer computer systems
 \cite{tanen06};
 (e) multi-layer communications \cite{tanen02};
 (f) multi-layer (hierarchical) information-communication networks
  (e.g.,
   \cite{abts12},\cite{bel97},\cite{bley04},\cite{kuz05},\cite{lev11ADES},\cite{mar01},\cite{mur99},\cite{obr10},\cite{shio04}).

\begin{center}
\begin{picture}(70,35)
\put(21,00){\makebox(0,0)[bl] {Fig. 1. Tree}}

%===================================================
\put(37,30){\makebox(0,0)[bl] {Root}}
%%%%%%%%%%%%%%%%%%%%%%%%%%%%%%%%
\put(35,30){\circle*{2.4}}

\put(15,25){\line(4,1){20}} \put(55,25){\line(-4,1){20}}

%%%%%%%%%%%%%%%%%%%%%%%%%%%%%%%%
\put(15,25){\circle*{2.0}}

\put(05,20){\line(2,1){10}} \put(22.5,20){\line(-3,2){7.5}}

%--Branch 1
\put(05,20){\circle*{1.4}}

\put(00,15){\line(1,1){05}} \put(05,15){\line(0,1){05}}
\put(10,15){\line(-1,1){05}}

\put(00,15){\circle*{1}} \put(05,15){\circle*{1}}
 \put(10,15){\circle*{1}}

%--Branch 2
\put(22.5,20){\circle*{1.4}}

\put(20,15){\line(1,2){2.5}} \put(25,15){\line(-1,2){2.5}}

\put(20,15){\circle*{1}} \put(25,15){\circle*{1}}

%%%%%%%%%%%%%%%%%%%%%%%%%%%%%%%%
\put(55,25){\circle*{2.0}}

\put(40,20){\line(3,1){15}} \put(62.5,20){\line(-3,2){7.5}}

%--Branch 3
\put(40,20){\circle*{1.4}}

\put(35,15){\line(1,1){05}} \put(40,15){\line(0,1){05}}
\put(45,15){\line(-1,1){05}}

\put(35,15){\circle*{1}} \put(40,15){\circle*{1}}
 \put(45,15){\circle*{1}}

%--Branch 4
\put(62.5,20){\circle*{1.4}}

\put(55,15){\line(3,2){7.5}} \put(60,15){\line(1,2){2.5}}
\put(65,15){\line(-1,2){2.55}} \put(70,15){\line(-3,2){7.5}}

\put(55,15){\circle*{1}} \put(60,15){\circle*{1}}
 \put(65,15){\circle*{1}}  \put(70,15){\circle*{1}}

%--

\put(20,09.5){\vector(0,1){4}} \put(25,09.5){\vector(0,1){4}}

%\put(30,10){\vector(0,1){4}}

\put(35,09.6){\vector(0,1){4}} \put(40,09.5){\vector(1,1){4}}

\put(21,06){\makebox(0,0)[bl] {Leaf vertices}}

\end{picture}
\end{center}

\begin{center}
\begin{picture}(76,28)
\put(23,00){\makebox(0,0)[bl] {Fig. 2. Forest}}

%%%%%%%%%%%%%%%%%%%%%%%%%%%%%%%% Component 1
\put(15,25){\circle*{2.0}}

\put(05,20){\line(2,1){10}} \put(17.5,20){\line(-1,2){2.5}}

%--Branch 1
\put(05,20){\circle*{1.4}}

\put(00,15){\line(1,1){05}} \put(05,15){\line(0,1){05}}
\put(10,15){\line(-1,1){05}}

\put(00,15){\circle*{1}} \put(05,15){\circle*{1}}
 \put(10,15){\circle*{1}}

%--Branch 2
\put(17.5,20){\circle*{1.4}}

\put(15,15){\line(1,2){2.5}} \put(20,15){\line(-1,2){2.5}}

\put(15,15){\circle*{1}} \put(20,15){\circle*{1}}

%%%%%%%%%%%%%%%%%%%%%%%%%%%%%%%%Component 2
\put(40,25){\circle*{2.0}}

\put(30,20){\line(2,1){10}} \put(47.5,20){\line(-3,2){7.5}}

%--Branch 3
\put(30,20){\circle*{1.4}}

\put(25,15){\line(1,1){05}} \put(30,15){\line(0,1){05}}
\put(35,15){\line(-1,1){05}}

\put(25,15){\circle*{1}} \put(30,15){\circle*{1}}
 \put(35,15){\circle*{1}}

%--Branch 4
\put(47.5,20){\circle*{1.4}}

\put(40,15){\line(3,2){7.5}} \put(45,15){\line(1,2){2.5}}
\put(50,15){\line(-1,2){2.55}} \put(55,15){\line(-3,2){7.5}}

\put(40,15){\circle*{1}} \put(45,15){\circle*{1}}
 \put(50,15){\circle*{1}}  \put(55,15){\circle*{1}}

%%%%%%%%%%%%%%%%%%%%%%%%%%%%%%%% Component 3
\put(67.5,25){\circle*{2.0}}

\put(62.5,20){\line(1,1){5}} \put(72.5,20){\line(-1,1){5}}

%--Branch 5
\put(62.5,20){\circle*{1.4}}

\put(60,15){\line(1,2){2.5}} \put(65,15){\line(-1,2){2.5}}

\put(60,15){\circle*{1}} \put(65,15){\circle*{1}}

%--Branch 6
\put(72.5,20){\circle*{1.4}}

\put(70,15){\line(1,2){2.5}} \put(75,15){\line(-1,2){2.5}}

\put(70,15){\circle*{1}} \put(75,15){\circle*{1}}

%%%%%%%%%%%%%%%%%%%%%%%%%%%%%%%%%%%%%%%%%%%%%%%%%%%%%%%%%%%%
%--

\put(12,09.5){\vector(0,1){4}} \put(39,09.5){\vector(0,1){4}}
\put(68,09.6){\vector(0,1){4}}

\put(00,06){\makebox(0,0)[bl] {Component \(1\)}}
\put(29,06){\makebox(0,0)[bl] {Component \(2\)}}
\put(58,06){\makebox(0,0)[bl] {Component \(3\)}}

\end{picture}
\end{center}

\begin{center}
\begin{picture}(78.5,23)
\put(24,00){\makebox(0,0)[bl] {Fig. 3. Polytree}}

\put(25,21){\circle*{2.4}}

\put(25,21){\vector(-2,-1){09}} \put(25,21){\vector(3,-1){14}}

%%%%%%%%%%%%%%%%%%%%%%%%%%%%%%%% Component 1
\put(15,16){\circle*{2.0}}

%\put(05,11){\line(2,1){10}} \put(17.5,11){\line(-1,2){2.5}}

\put(15,16){\vector(-2,-1){09.4}} \put(15,16){\vector(1,-2){2.3}}

%--Branch 1
\put(05,11){\circle*{1.4}}

\put(05,11){\vector(-1,-1){04.5}} \put(05,11){\vector(0,-1){04.5}}
\put(05,11){\vector(1,-1){04.5}}

\put(00,6){\circle*{1}} \put(05,6){\circle*{1}}
 \put(10,6){\circle*{1}}

%--Branch 2
\put(17.5,11){\circle*{1.4}}

\put(17.5,11){\vector(-1,-2){2.3}}
\put(17.5,11){\vector(1,-2){2.3}}

\put(15,6){\circle*{1}} \put(20,6){\circle*{1}}

%%%%%%%%%%%%%%%%%%%%%%%%%%%%%%%%
\put(60,21){\circle*{2.4}}

\put(60,21){\vector(-4,-1){19}} \put(60,21){\vector(3,-2){7}}

%%%%%%%%%%%%%%%%%%%%%%%%%%%%%%%%Component 2
\put(40,16){\circle*{2.0}}

\put(40,16){\vector(-2,-1){09}} \put(40,16){\vector(3,-2){7}}

%--Branch 3
\put(30,11){\circle*{1.4}}

\put(30,11){\vector(-1,-1){04.5}} \put(30,11){\vector(0,-1){04.5}}
\put(30,11){\vector(1,-1){04.5}}

\put(25,6){\circle*{1}} \put(30,6){\circle*{1}}
 \put(35,6){\circle*{1}}

%--Branch 4
\put(47.5,11){\circle*{1.4}}

\put(47.5,11){\vector(-3,-2){7}}
\put(47.5,11){\vector(-1,-2){2.3}}
\put(47.5,11){\vector(1,-2){2.3}} \put(47.5,11){\vector(3,-2){7}}

\put(40,6){\circle*{1}} \put(45,6){\circle*{1}}
 \put(50,6){\circle*{1}}  \put(55,6){\circle*{1}}

%%%%%%%%%%%%%%%%%%%%%%%%%%%%%%%% Component 3'
\put(77.5,16){\circle*{2.0}}

\put(77.5,16){\vector(-1,-1){4.5}}

%%%%%%%%%%%%%%%%%%%%%%%%%%%%%%%% Component 3
\put(67.5,16){\circle*{2.0}}

\put(67.5,16){\vector(-1,-1){4.5}}
\put(67.5,16){\vector(1,-1){4.5}}

%--Branch 5
\put(62.5,11){\circle*{1.4}}

\put(62.5,11){\vector(-1,-2){2.3}}
\put(62.5,11){\vector(1,-2){2.3}}

\put(60,6){\circle*{1}} \put(65,6){\circle*{1}}

%--Branch 6
\put(72.5,11){\circle*{1.4}}

\put(72.5,11){\vector(-1,-2){2.3}}
\put(72.5,11){\vector(1,-2){2.3}}

\put(70,6){\circle*{1}} \put(75,6){\circle*{1}}

\end{picture}
\end{center}

\begin{center}
\begin{picture}(55,38)
\put(10,00){\makebox(0,0)[bl] {Fig. 4. Organic hierarchy}}

%===================================================

%%%%%%%%%%%%%%%%%%%%%%%%%%%%%%%%
\put(25,35){\circle*{2.0}}

\put(10,25){\line(3,2){15}} \put(40,25){\line(-3,2){15}}

%--------------

\put(10,25){\line(1,0){3}} \put(15,25){\line(1,0){3}}
\put(20,25){\line(1,0){3}} \put(25,25){\line(1,0){3}}
\put(30,25){\line(1,0){3}} \put(35,25){\line(1,0){3}}
%\put(40,25){\line(1,0){3}}

%--Branch 1
\put(10,25){\circle*{1.4}}

\put(00,15){\line(1,1){10}} \put(10,15){\line(0,1){10}}
\put(20,15){\line(-1,1){10}}

\put(00,15){\circle*{1}} \put(10,15){\circle*{1}}
 \put(20,15){\circle*{1}}

\put(00,15){\line(1,0){3}} \put(05,15){\line(1,0){3}}
\put(10,15){\line(1,0){3}} \put(15,15){\line(1,0){3}}

%--Branch 2
\put(40,25){\circle*{1.4}}

\put(30,15){\line(1,1){10}} \put(40,15){\line(0,1){10}}
\put(50,15){\line(-1,1){10}} \put(60,15){\line(-2,1){20}}

\put(30,15){\circle*{1}} \put(40,15){\circle*{1}}
 \put(50,15){\circle*{1}}  \put(60,15){\circle*{1}}

\put(30,15){\line(1,0){3}} \put(35,15){\line(1,0){3}}
\put(40,15){\line(1,0){3}} \put(45,15){\line(1,0){3}}
\put(50,15){\line(1,0){3}} \put(55,15){\line(1,0){3}}
%--

\put(15,09.5){\vector(0,1){4}} \put(25,09.5){\vector(0,1){14}}

\put(45,09.5){\vector(0,1){4}}

\put(11,06){\makebox(0,0)[bl] {Organic interconnection}}

\end{picture}
\end{center}

\begin{center}
\begin{picture}(70,35)

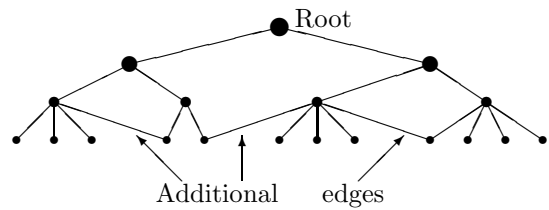
\put(00,00){\makebox(0,0)[bl] {Fig. 5. Hierarchy (tree with
additional edges)}}
%===================================================
\put(37,30){\makebox(0,0)[bl] {Root}}
%%%%%%%%%%%%%%%%%%%%%%%%%%%%%%%%
\put(35,30){\circle*{2.4}}

\put(15,25){\line(4,1){20}} \put(55,25){\line(-4,1){20}}

%%%%%%%%%%%%%%%%%%%%%%%%%%%%%%%%
\put(15,25){\circle*{2.0}}

\put(05,20){\line(2,1){10}} \put(22.5,20){\line(-3,2){7.5}}

%--Branch 1
\put(05,20){\circle*{1.4}}

\put(00,15){\line(1,1){05}} \put(05,15){\line(0,1){05}}
\put(10,15){\line(-1,1){05}}

\put(00,15){\circle*{1}} \put(05,15){\circle*{1}}
 \put(10,15){\circle*{1}}

%--Branch 2
\put(22.5,20){\circle*{1.4}}

\put(20,15){\line(1,2){2.5}} \put(25,15){\line(-1,2){2.5}}

\put(20,15){\circle*{1}} \put(25,15){\circle*{1}}

%--Additional edge

\put(20,15){\line(-3,1){15}}

%--Additional edge

\put(25,15){\line(3,1){15}}

%%%%%%%%%%%%%%%%%%%%%%%%%%%%%%%%
\put(55,25){\circle*{2.0}}

\put(40,20){\line(3,1){15}} \put(62.5,20){\line(-3,2){7.5}}

%--Branch 3
\put(40,20){\circle*{1.4}}

\put(35,15){\line(1,1){05}} \put(40,15){\line(0,1){05}}
\put(45,15){\line(-1,1){05}}

\put(35,15){\circle*{1}} \put(40,15){\circle*{1}}
 \put(45,15){\circle*{1}}

%--Branch 4
\put(62.5,20){\circle*{1.4}}

\put(55,15){\line(3,2){7.5}} \put(60,15){\line(1,2){2.5}}
\put(65,15){\line(-1,2){2.55}} \put(70,15){\line(-3,2){7.5}}

\put(55,15){\circle*{1}} \put(60,15){\circle*{1}}
 \put(65,15){\circle*{1}}  \put(70,15){\circle*{1}}

%--Additional edge

\put(55,15){\line(-3,1){15}}

%--

\put(22,09.5){\vector(-1,1){6}} \put(30,09.5){\vector(0,1){6}}

%\put(30,10){\vector(0,1){4}}

%\put(35,09.6){\vector(0,1){4}}

\put(45.5,09.5){\vector(1,1){6}}

\put(18.5,06){\makebox(0,0)[bl] {Additional ~~~ edges}}

\end{picture}
\end{center}

\begin{center}
\begin{picture}(74,53)
\put(00,00){\makebox(0,0)[bl]{Fig. 6.
 ``Morphological'' system hierarchy
 (\cite{lev11agg},\cite{lev11ADES})}}

%-----

\put(37,51){\circle*{2.5}}

%-----

\put(23.6,43){\makebox(0,0)[bl]{System hierarchy}}
\put(21.5,40){\makebox(0,0)[bl]{(tree-like structure)}}

%----

\put(01,34){\makebox(0,0)[bl]{Leaf vertices}}

\put(04,33){\line(0,-1){8}}

\put(05,33){\line(2,-3){5.4}}

\put(6,33){\line(3,-1){23}}

%-------------------------

\put(24,36){\makebox(0,0)[bl]{Alternatives for}}
\put(27.5,33.5){\makebox(0,0)[bl]{leaf vertices}}

\put(30,33){\line(-1,-1){14}} \put(43,33){\line(1,-1){15}}

%-------------------------

\put(52,36){\makebox(0,0)[bl]{Compatibility}}
\put(52,33){\makebox(0,0)[bl]{among}}
\put(52,30.5){\makebox(0,0)[bl]{alternatives}}

\put(60.4,30){\line(-1,-1){20.5}}

%------------------------------------

\put(00,23){\line(1,0){74}}

\put(00,39){\line(3,1){37}} \put(74,39){\line(-3,1){37}}

\put(00,23){\line(0,1){16}} \put(74,23){\line(0,1){16}}

%======================
%-- A1
\put(02,24.5){\makebox(0,0)[bl]{\(1\)}}

\put(03,23){\circle*{1.8}} \put(03,18){\oval(5,8)}

%-- A2
\put(12,24.5){\makebox(0,0)[bl]{\(2\)}}

\put(13,23){\circle*{1.8}} \put(13,18){\oval(5,8)}

%-- A(k-1)
\put(29.5,24){\makebox(0,0)[bl]{\(\tau-1\)}}

\put(32,23){\circle*{1.8}} \put(32,18){\oval(5,8)}

%-- Ak
\put(41,24.5){\makebox(0,0)[bl]{\(\tau\)}}

\put(42,23){\circle*{1.8}} \put(42,18){\oval(5,8)}

%-- A(m-1)
\put(57,24.5){\makebox(0,0)[bl]{\(m-1\)}}

\put(61,23){\circle*{1.8}} \put(61,18){\oval(5,8)}

%-- Am
\put(70,24.5){\makebox(0,0)[bl]{\(m\)}}

\put(71,23){\circle*{1.8}} \put(71,18){\oval(5,8)}

%=================
\put(19.5,17.5){\makebox(0,0)[bl]{. . .}}
\put(48.5,17.5){\makebox(0,0)[bl]{. . .}}

\put(19.5,9){\makebox(0,0)[bl]{. . .}}
\put(48.5,9){\makebox(0,0)[bl]{. . .}}
%================
%--- 1-2
\put(06,06){\line(1,0){04}} \put(06,13){\line(1,0){04}}
\put(06,06){\line(0,1){7}} \put(10,06){\line(0,1){7}}

\put(07,06){\line(0,1){7}} \put(08,06){\line(0,1){7}}
\put(09,06){\line(0,1){7}}

%--- (k-1)-k
\put(35,06){\line(1,0){04}} \put(35,13){\line(1,0){04}}
\put(35,06){\line(0,1){7}} \put(39,06){\line(0,1){7}}

\put(36,06){\line(0,1){7}} \put(37,06){\line(0,1){7}}
\put(38,06){\line(0,1){7}}

%--- (m-1)-m
\put(64,06){\line(1,0){04}} \put(64,13){\line(1,0){04}}
\put(64,06){\line(0,1){7}} \put(68,06){\line(0,1){7}}

\put(65,06){\line(0,1){7}} \put(66,06){\line(0,1){7}}
\put(67,06){\line(0,1){7}}

%----------------

\end{picture}
\end{center}

\begin{center}
\begin{picture}(55,35)

\put(07,00){\makebox(0,0)[bl]{Fig. 7. Multilayer structure}}

%---- Top- layer

\put(39,30){\makebox(0,0)[bl]{Top}}
\put(39,27){\makebox(0,0)[bl]{layer}}

\put(20,30){\oval(32,6)}

\put(09,30){\circle*{1.8}} \put(13,30){\circle*{1.8}}
\put(17.5,29.5){\makebox(0,0)[bl]{{\bf .~.~.}}}

\put(27,30){\circle*{1.8}}  \put(31,30){\circle*{1.8}}

%--

\put(09,30){\vector(0,-1){06}} \put(09,30){\vector(-1,-2){03}}

\put(13,30){\vector(0,-1){06}}

\put(27,30){\vector(0,-1){06}} \put(27,30){\vector(-1,-2){03}}

\put(31,30){\vector(0,-1){06}} \put(31,30){\vector(-1,-2){03}}

\put(31,30){\vector(1,-2){03}}

%---- Intermediate

\put(40,20){\makebox(0,0)[bl]{Intermediate}}
\put(40,17){\makebox(0,0)[bl]{layer}}

\put(20,20){\oval(36,6)}

\put(07,20){\circle*{1.5}} \put(12,20){\circle*{1.5}}
\put(17.5,19.5){\makebox(0,0)[bl]{{\bf .~.~.}}}

\put(28,20){\circle*{1.5}}  \put(33,20){\circle*{1.5}}

%--

\put(07,20){\vector(0,-1){06}} \put(07,20){\vector(-1,-2){03}}
\put(07,20){\vector(1,-2){03}}

\put(12,20){\vector(0,-1){06}} \put(12,20){\vector(1,-2){03}}

\put(28,20){\vector(0,-1){06}} \put(28,20){\vector(-1,-2){03}}

\put(33,20){\vector(0,-1){06}} \put(33,20){\vector(-1,-2){03}}
\put(33,20){\vector(1,-2){03}}

%---- Bottom layer

\put(41,10){\makebox(0,0)[bl]{Bottom}}
\put(41,07){\makebox(0,0)[bl]{layer}}

\put(20,10){\oval(40,6)}

\put(05,10){\circle*{1}} \put(10,10){\circle*{1}}
\put(17.5,09.5){\makebox(0,0)[bl]{{\bf .~.~.}}}

\put(30,10){\circle*{1}}  \put(35,10){\circle*{1}}

\end{picture}
\end{center}

%%%%%%%%%%%%%%%%%%%%%%%%%%%%%%%%%%%%%%%%%%%%%%%%%%%%

 In applied domains,
 many special types of tree-like structures or hierarchies
 are widely used, for example:
  (1) hierarchical schemes for data, for information systems
 (e.g., \cite{bota92}, \cite{knuth68}, \cite{levach07}, \cite{lev89},
 \cite{ras85});
 (2) thesauri and concept spaces
 (e.g., \cite{chen93},\cite{chen96});
 (3) organizational hierarchies
 (e.g.,
  \cite{baligh06}, \cite{duncan79},
  \cite{gubko06}, \cite{mishin07}, \cite{voronin02});
 (4) multi-level complex systems
 (e.g., \cite{mes70});
 (5) phylogenetic trees
  (e.g., \cite{mirkin96},\cite{rob81},\cite{sai89})
  and evolutionary trees
  (e.g., \cite{amir97},\cite{mirkin96});
 (6) ontologies
 (e.g.,
 \cite{corcho03},
 \cite{hol02},\cite{noy97},\cite{usch96});
 (7) statecharts
 (e.g., \cite{bog01}, \cite{harel87});
 (8) decision trees
 (e.g.,
  \cite{ait08},
   \cite{gehrke09}, \cite{gehrke00},
  \cite{good88}, \cite{good90}, \cite{qui86}, \cite{qui90});
 (9) hierarchy of criteria in decision making
 (Analytic Hierarchy Process)
  \cite{saaty88};
  and
 (10) hierarchical access networks
   (e.g., \cite{godor05}).

 Fig. 8 depicts a basic design process to obtain
 a hierarchical structure.

\begin{center}
\begin{picture}(74,61)

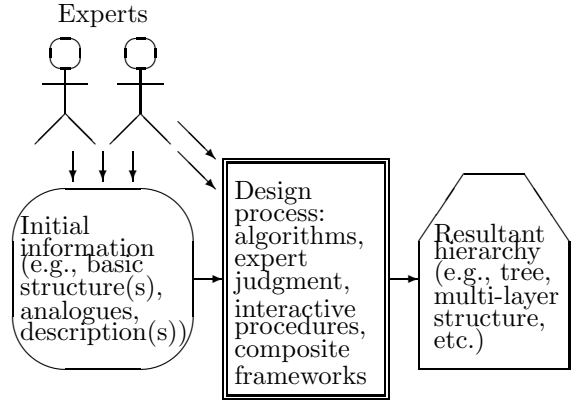
\put(15,00){\makebox(0,0)[bl]{Fig. 8. Building of hierarchy}}

%---- 2 experts

\put(06,55.5){\makebox(0,0)[bl]{Experts}}

%--
\put(07,44){\line(0,1){6}}

\put(03,40){\line(1,1){4}} \put(11,40){\line(-1,1){4}}

\put(04,48){\line(1,0){6}}

\put(07,52){\oval(4,4)}

%--
\put(17,44){\line(0,1){6}}

\put(13,40){\line(1,1){4}} \put(21,40){\line(-1,1){4}}

\put(14,48){\line(1,0){6}}

\put(17,52){\oval(4,4)}

%-

\put(08,39){\vector(0,-1){04}} \put(12,39){\vector(0,-1){04}}
\put(16,39){\vector(0,-1){04}}

%-

%\put(23,48){\vector(1,-1){04}}

\put(22,43){\vector(1,-1){05}} \put(22,39){\vector(1,-1){05}}

%===============================================

\put(12,22){\oval(24,24)}

%\put(40,10){\oval(33,7)}

\put(01,28){\makebox(0,0)[bl]{Initial}}
\put(01,25){\makebox(0,0)[bl]{information}}
\put(01,22){\makebox(0,0)[bl]{(e.g., basic}}
\put(01,19){\makebox(0,0)[bl]{structure(s),}}
\put(01,16){\makebox(0,0)[bl]{analogues,}}
\put(01,13){\makebox(0,0)[bl]{description(s))}}

\put(24,22){\vector(1,0){04}}

%--------------------------- Aggregation process

\put(28,06){\line(1,0){22}} \put(28,38){\line(1,0){22}}
\put(28,06){\line(0,1){32}} \put(50,06){\line(0,1){32}}

\put(28.5,6.5){\line(1,0){21}} \put(28.5,37.5){\line(1,0){21}}
\put(28.5,6.5){\line(0,1){31}} \put(49.5,6.5){\line(0,1){31}}

\put(29.5,32){\makebox(0,0)[bl]{Design}}
\put(29.5,29){\makebox(0,0)[bl]{process:}}
\put(29.5,26){\makebox(0,0)[bl]{algorithms,}}
\put(29.5,23){\makebox(0,0)[bl]{expert}}
\put(29.5,20){\makebox(0,0)[bl]{judgment,}}
\put(29.5,17){\makebox(0,0)[bl]{interactive}}
\put(29.5,14){\makebox(0,0)[bl]{procedures,}}
\put(29.5,11){\makebox(0,0)[bl]{composite}}
\put(29.5,8){\makebox(0,0)[bl]{frameworks}}

\put(50,22){\vector(1,0){04}}

%===============================================

%\put(64,22){\oval(20,24)} \put(64,22){\oval(19,23)}

\put(56,27){\makebox(0,0)[bl]{Resultant}}
\put(56,24){\makebox(0,0)[bl]{hierarchy}}
\put(56,21){\makebox(0,0)[bl]{(e.g., tree,}}
\put(56,18){\makebox(0,0)[bl]{multi-layer}}
\put(56,15){\makebox(0,0)[bl]{structure,}}
\put(56,12){\makebox(0,0)[bl]{etc.)}}

\put(54,10){\line(1,0){20}}

\put(54,10){\line(0,1){17}} \put(74,10){\line(0,1){17}}

\put(54,27){\line(2,3){6}} \put(74,27){\line(-2,3){6}}

\put(60,36){\line(1,0){08}}

\end{picture}
\end{center}

 It is reasonable to list our main types of
 hierarchy design problems as the following:

~~

  {\bf Problem 1.} Expert-based design:

 {\it Input}: description of a system.

 {\it Output}: hierarchy as a result of system partitioning
 (tree, hierarchy, layered structure).

 {\bf Problem 2.} Basic design:

 {\it Input}: set of elements, element attributes or element
 interconnection.

 {\it Output}: set of clusters or
  spanning/covering structure over the initial element set
 (tree, hierarchy, layered structure).

 {\bf Problem 3.} Spanning/covering:

 {\it Input}: initial network/graph (i.e., set of elements and set
 of edges).

 {\it Output}:  spanning/covering hierarchical
 structure (tree, hierarchy, layered structure).

 {\bf Problem 4.} Redesign (modification, transformation,
 improvement):

 {\it Input}: initial hierarchical structure
 (tree, hierarchy, layered structure).

 {\it Output}: new hierarchical structure
 with some required features
 (tree, hierarchy, layered structure).

 {\bf Problem 5.} Restructuring (special case of modification):

 {\it Input}:
 (i) initial hierarchical structure
 (tree, hierarchy, layered structure),
 (ii) goal hierarchical structure
 (tree, hierarchy, layered structure).

 {\it Output}: new hierarchical structure
 (tree, hierarchy, layered structure)
 while taking into account the following:
 (a) ``cheap'' transformation of the initial structure,
 (b) ``small'' proximity between the new structure and the goal structure.

~~

 Note,
 an important problem corresponds to  aggregating some several initial hierarchies into a resultant
 aggregated structure
 (e.g., aggregation/integration/merging of information systems,
 data base schemas,
 catalogs, ontologies,
 knowledge bases,
 organizational structures,
 experts preferences)
 has been intensively examined
 (e.g.,  \cite{agrawal01}, \cite{batini86}, \cite{chen93},\cite{cook96},
 \cite{dan01}, \cite{lin96},
   \cite{noy00}, \cite{noy05},
  \cite{pinto01},
   \cite{wache01},
   \cite{yager88}).
 An author recent survey
 on aggregation of some structures
 is presented in \cite{lev11agg}.

%%%%%%%%%%%%%%%%%%%%%%%%%%%%%%%%%%%%%%%%%%%%%%
\subsection{Expert-based 'Top-Down' Procedure}

 Mainly, expert-based procedures for building a system hierarchy
 are based on domain experts
 and some typical ``technological' frames
 (e.g., product life cycle as the following:
 design, manufacturing, testing, maintenance, utilization,
 recycling).
 This procedure consists of the following phases
 (it is the 'divisive' strategy of hierarchical clustering):
 {\it 1.} dividing a system into its subsystems;
 {\it 2.} dividing each subsystem into its parts;
 {\it 3.} dividing each subsystem part into its components;
 etc.
 It is reasonable to point out the basic algorithmic rules for
 dividing the system (subsystem, subsystem parts, etc.):
 (a) dividing (partitioning) by physical parts,
 (b) dividing by system functions,
 (c) dividing by time stages of data processing.
 Three applied illustrative examples are presented to illustrate the
 procedure above:

 (1) for concrete macrotechnology (Fig. 9)
 (\cite{lev06}, \cite{levnis01}).;

 (2) for a two-floor building (Fig. 10)
 (\cite{lev06}, \cite{levdan05});
 and

 (3) for medical treatment (children asthma) (Fig. 11)
 (\cite{lev06}, \cite{levsok04}).

\begin{center}
\begin{picture}(95,55)

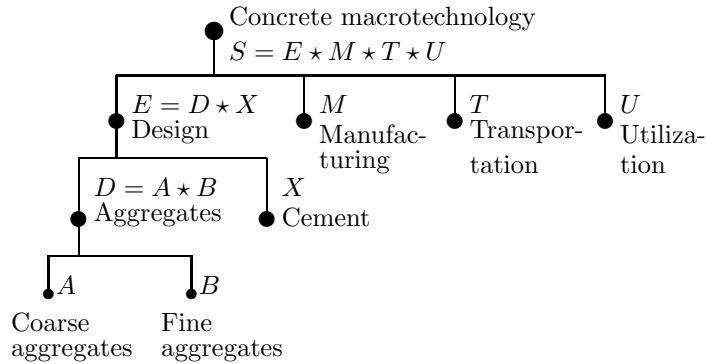
\put(01,00){\makebox(0,0)[bl] {Fig. 9. Hierarchy of concrete
 macrotechnology (\cite{lev06},\cite{levnis01})}}
%===================================================

\put(30,50){\makebox(0,0)[bl]{Concrete macrotechnology}}
\put(30,46){\makebox(0,0)[bl]{\(S = E \star M \star T \star U\)}}

\put(28,50){\circle*{2.5}}

\put(28,44){\line(0,1){6}}

\put(15,44){\line(1,0){65}}

%- - - - - - - -

\put(80,38){\line(0,1){06}}

\put(82,39){\makebox(0,0)[bl]{\(U\)}}
\put(82,35){\makebox(0,0)[bl]{Utiliza-}}
\put(82,31){\makebox(0,0)[bl]{tion}}

\put(80,38){\circle*{2}}

%- - - - - - - -

\put(60,38){\line(0,1){06}}

\put(62,39){\makebox(0,0)[bl]{\(T\)}}
\put(62,35){\makebox(0,0)[bl]{Transpor-}}
\put(62,31){\makebox(0,0)[bl]{tation}}

\put(60,38){\circle*{2}}

%- - - - - - - -

\put(40,38){\line(0,1){06}}

\put(42,39){\makebox(0,0)[bl]{\(M\)}}
\put(42,35){\makebox(0,0)[bl]{Manufac-}}
\put(42,31){\makebox(0,0)[bl]{turing}}

 \put(40,38){\circle*{2}}

%- - - - - - - -

\put(15,38){\line(0,1){06}}

\put(17,39){\makebox(0,0)[bl]{\(E=D\star X \)}}
\put(17,35){\makebox(0,0)[bl]{Design}}

\put(15,38){\circle*{2}}

\put(15,33){\line(0,1){10}}
%---

\put(10,33){\line(1,0){25}}

%- - - - - - - -

\put(35,25){\line(0,1){08}}

\put(37,28){\makebox(0,0)[bl]{\(X \)}}
\put(37,24){\makebox(0,0)[bl]{Cement}}

\put(35,25){\circle*{2}}

%- - - - - - - -

\put(10,25){\line(0,1){08}}

\put(12,28){\makebox(0,0)[bl]{\(D=A \star B \)}}
\put(12,24){\makebox(0,0)[bl]{Aggregates}}

\put(10,25){\circle*{2}}

%--

\put(10,20){\line(0,1){05}}

%--

\put(06,20){\line(1,0){19}}

%--

\put(06,15){\line(0,1){05}} \put(06,15){\circle*{1.5}}
\put(07,15){\makebox(0,0)[bl]{\(A\)}}
\put(01,10){\makebox(0,0)[bl]{Coarse}}
\put(01,06){\makebox(0,0)[bl]{aggregates}}

%--

\put(25,15){\line(0,1){05}} \put(25,15){\circle*{1.5}}
\put(26,15){\makebox(0,0)[bl]{\(B\)}}
\put(21,10){\makebox(0,0)[bl]{Fine}}
\put(21,06){\makebox(0,0)[bl]{aggregates}}

\end{picture}
\end{center}

\begin{center}
%\begin{picture}(83,36)
\begin{picture}(82,52)

\put(07,00){\makebox(0,0)[bl] {Fig. 10. Hierarchy of
 building (\cite{lev06},\cite{levdan05})}}
%===================================================

\put(27,47){\makebox(0,0)[bl]{Building \(S = A \star B \star C
\)}}

\put(25,49){\circle*{2.5}}

\put(25,45){\line(0,1){4}}

\put(15,45){\line(1,0){55}}

%- - - - - - - -

\put(70,38){\line(0,1){07}}

\put(72,41){\makebox(0,0)[bl]{\(C\)}}
\put(72,37){\makebox(0,0)[bl]{Floors}}

\put(70,38){\circle*{2}}

%- - - - - - - -

\put(40,38){\line(0,1){07}}

\put(42,41){\makebox(0,0)[bl]{\(B=D\star F\)}}
\put(42,37){\makebox(0,0)[bl]{Basic structure}}

 \put(40,38){\circle*{2}}

 \put(40,38){\line(0,-1){05}}

%- - - - - - - -

\put(15,38){\line(0,1){07}}

\put(17,41){\makebox(0,0)[bl]{\(A\)}}
\put(17,37){\makebox(0,0)[bl]{Foundation}}

\put(15,38){\circle*{2}}

%\put(15,33){\line(0,1){10}}
%---

\put(10,33){\line(1,0){45}}

%- - - - - - - -

\put(55,25){\line(0,1){08}}

\put(57,28){\makebox(0,0)[bl]{\(F=I\star J \)}}
\put(57,24){\makebox(0,0)[bl]{Nonbearing}}
\put(57,21){\makebox(0,0)[bl]{structures}}

\put(55,25){\circle*{2}}

\put(55,25){\line(0,-1){05}}

%- - - - - - - -

\put(10,25){\line(0,1){08}}

\put(12,28){\makebox(0,0)[bl]{\(D=E \star G \star H \)}}
\put(12,24){\makebox(0,0)[bl]{Bearing}}
\put(12,21){\makebox(0,0)[bl]{structures}}

\put(10,25){\circle*{2}}

%--

\put(10,20){\line(0,1){05}}

%--

\put(06,20){\line(1,0){25}}

%--

\put(06,15){\line(0,1){05}} \put(06,15){\circle*{1.5}}
\put(07,15){\makebox(0,0)[bl]{\(E\)}}
\put(01,10){\makebox(0,0)[bl]{Frame}}

%--

\put(17,15){\line(0,1){05}} \put(17,15){\circle*{1.5}}
\put(18,15){\makebox(0,0)[bl]{\(G\)}}
\put(13,09.5){\makebox(0,0)[bl]{Rigidity}}
\put(13,06){\makebox(0,0)[bl]{core}}

%--

\put(31,15){\line(0,1){05}} \put(31,15){\circle*{1.5}}
\put(32,15){\makebox(0,0)[bl]{\(H\)}}
\put(27,10){\makebox(0,0)[bl]{Stair-}}
\put(27,06.5){\makebox(0,0)[bl]{case}}

%%%%%%%%%%%%%%%%%%%%%

%--

\put(46,20){\line(1,0){14}}

%--

\put(46,15){\line(0,1){05}} \put(46,15){\circle*{1.5}}
\put(47,15){\makebox(0,0)[bl]{\(I\)}}
\put(41,10){\makebox(0,0)[bl]{Filler}}
\put(41,06){\makebox(0,0)[bl]{walls}}

%--

\put(60,15){\line(0,1){05}} \put(60,15){\circle*{1.5}}
\put(61,15){\makebox(0,0)[bl]{\(J\)}}
\put(53,09.5){\makebox(0,0)[bl]{Partitioning}}
\put(53,06){\makebox(0,0)[bl]{walls}}

\end{picture}
\end{center}

\begin{center}
\begin{picture}(111,42)

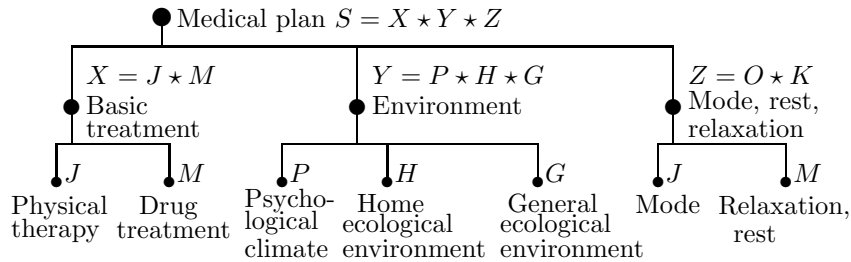
\put(05,00){\makebox(0,0)[bl] {Fig. 11. Hierarchy of medical
 plan (children asthma) (\cite{lev06},\cite{levsok04})}}
%===================================================

\put(22,36){\makebox(0,0)[bl]{Medical plan \(S = X \star Y \star Z
\)}}

\put(20,38){\circle*{2.5}} \put(20,34){\line(0,1){04}}

\put(08,34){\line(1,0){80}}

%- - - - - - - -

\put(08,26){\line(0,1){08}}

\put(10,29){\makebox(0,0)[bl]{\(X=J \star M \)}}
\put(10,25){\makebox(0,0)[bl]{Basic}}
\put(10,22){\makebox(0,0)[bl]{treatment}}

\put(08,26){\circle*{2}}

%--

\put(08,21){\line(0,1){05}}

%--

\put(06,21){\line(1,0){15}}

%--

\put(06,16){\line(0,1){05}} \put(06,16){\circle*{1.5}}
\put(07,16){\makebox(0,0)[bl]{\(J\)}}
\put(00,11){\makebox(0,0)[bl]{Physical}}
\put(00,08){\makebox(0,0)[bl]{therapy}}

%--

\put(21,16){\line(0,1){05}} \put(21,16){\circle*{1.5}}
\put(22,16){\makebox(0,0)[bl]{\(M\)}}
\put(17,11){\makebox(0,0)[bl]{Drug}}
\put(14,08.5){\makebox(0,0)[bl]{treatment}}

%- - - - - - - -

\put(46,26){\line(0,1){08}}

\put(48,29){\makebox(0,0)[bl]{\(Y=P \star H\star G \)}}
\put(48,25){\makebox(0,0)[bl]{Environment}}

\put(46,26){\circle*{2}}

%--

\put(46,21){\line(0,1){05}}

%--

\put(36,21){\line(1,0){34}}

%--

\put(36,16){\line(0,1){05}} \put(36,16){\circle*{1.5}}
\put(37,16){\makebox(0,0)[bl]{\(P\)}}
\put(31,12){\makebox(0,0)[bl]{Psycho-}}
\put(31,09){\makebox(0,0)[bl]{logical}}
\put(31,06){\makebox(0,0)[bl]{climate}}

%--

\put(50,16){\line(0,1){05}} \put(50,16){\circle*{1.5}}
\put(51,16){\makebox(0,0)[bl]{\(H\)}}
\put(46,12){\makebox(0,0)[bl]{Home }}
\put(44,08.6){\makebox(0,0)[bl]{ecological}}
\put(44,06){\makebox(0,0)[bl]{environment}}

%--

\put(70,16){\line(0,1){05}} \put(70,16){\circle*{1.5}}
\put(71,16){\makebox(0,0)[bl]{\(G\)}}
\put(66,12){\makebox(0,0)[bl]{General}}
\put(65,08.6){\makebox(0,0)[bl]{ecological}}
\put(65,06){\makebox(0,0)[bl]{environment}}

%%%%%%%%%%%%%%%%%%%%%%%%%%%%% Z=O-K

%- - - - - - - -

\put(88,26){\line(0,1){08}}

\put(90,29){\makebox(0,0)[bl]{\(Z=O \star K \)}}
\put(90,25){\makebox(0,0)[bl]{Mode, rest,}}
\put(90,22){\makebox(0,0)[bl]{relaxation}}

\put(88,26){\circle*{2}}

%--

\put(88,21){\line(0,1){05}}

%--

\put(86,21){\line(1,0){17}}

%--

\put(86,16){\line(0,1){05}} \put(86,16){\circle*{1.5}}
\put(87,16){\makebox(0,0)[bl]{\(J\)}}
\put(83,11.8){\makebox(0,0)[bl]{Mode}}

%--

\put(103,16){\line(0,1){05}} \put(103,16){\circle*{1.5}}
\put(104,16){\makebox(0,0)[bl]{\(M\)}}
\put(94,11){\makebox(0,0)[bl]{Relaxation,}}
\put(96,08){\makebox(0,0)[bl]{rest}}

\end{picture}
\end{center}

%%%%%%%%%%%%%%%%%%%%%%%%%%%%%%%%%%%%%%%%%%%%%%%%%%%%%%%%%%%%%%%%%%%%%%%%%
\subsection{Hierarchical Clustering (Agglomerative Algorithm)}

 Hierarchical clustering consists in building a hierarchy (tree-like structure) of
 clusters.
 There is a set of \(n\) elements
 ~\(A = \{ A_{1},...,A_{i},...,A_{n}\}\)~
 and a corresponding vector estimate of
 \(m\) attributes/parameters
 (\(T_{1},...,T_{j},...,T_{m}\))
  for each element \(i\):
 ~\(z_{i} = ( z_{i,1},...,z_{i,j},...,z_{i,m} )\).
 The basic agglomerative algorithm
 (polynomial, {\it algorithm 1}) is as follows ('Bottom-Up' element pair
 integration process)
 (e.g., \cite{gordon99}, \cite{sneath73}):

~~

 {\it Stage 1.} Computing the matrix of element pair
 \(\forall (A(i_{1}),A(i_{2}))\), ~\(A(i_{1})\in A\), ~\(A(i_{2})\in A \),
 ~\(i_{1} \neq i_{2}\)
 ``distances'' (a simple case,
 Euclidean distance):
 \[ d_{i_{1}i_{2}} = \sqrt{ \sum_{j=1}^{m} ( z_{i_{1},j} - z_{i_{2},j} )^{2} }. \]

 {\it Stage2.} Revelation of the smallest pair ``distance''
 and integration of the corresponding two elements into a resultant
 ``integrated'' element.

 {\it Stage 3.} Stopping process or re-computing the matrix of
 pair ``distances'' and ~{\it Go To}~~ {\it Stage 2}.

~~

 As result, a tree-like
 structure for the element pair integration process ('Bottom-Up') is
 obtained (one element pair integration at each integration step).
 A basic procedure for aggregation of items
 (aggregation as average values) is as follows
 (\(J_{i_{1},i_{2}}= A_{i_{1}} \& A_{i_{2}}\)):
 ~\(\forall j ~~ z_{J_{i_{1},i_{2}},j} =   \frac{z_{i_{1},j} + z_{i_{2},j}} {2}
 \).
 The item pair aggregation process can be based on other
 functions ({\it e.g.}, \(\max \), \(\min \) ).
 Integration of several
 items can be considered analogically.
 An illustrative example
 corresponds to an analysis of \(8\) persons by
 their inclination/interests
 (Table 1).
 Here four attributes are used
 (ordinal scale \([0,5]\)):
 (i) inclination for mathematics or logical thinking
 (\(K_{1}\)),
 (ii) interest to music  (\(K_{2}\)),
 (iii) interest to sport  (\(K_{3}\)),
 and
 (iv) interest to trips   (\(K_{4}\)).
 The result of the basic hierarchical clustering
 (as a hierarchical structure for eight-person team)
 is depicted in Fig. 12.

\begin{center}
\begin{picture}(50,48.5)

\put(0.5,44){\makebox(0,0)[bl]{Table 1. Criteria and estimates}}

\put(00,00){\line(1,0){50}} \put(00,35){\line(1,0){50}}
\put(00,42){\line(1,0){50}}

\put(00,0){\line(0,1){42}} \put(10,0){\line(0,1){42}}
\put(20,0){\line(0,1){42}} \put(30,0){\line(0,1){42}}
\put(40,0){\line(0,1){42}} \put(50,0){\line(0,1){42}}

\put(01,30){\makebox(0,0)[bl]{\(A_{1}\)}}
\put(01,26){\makebox(0,0)[bl]{\(A_{2}\)}}
\put(01,22){\makebox(0,0)[bl]{\(A_{3}\)}}
\put(01,18){\makebox(0,0)[bl]{\(A_{4}\)}}

\put(01,14){\makebox(0,0)[bl]{\(A_{5}\)}}
\put(01,10){\makebox(0,0)[bl]{\(A_{6}\)}}
\put(01,06){\makebox(0,0)[bl]{\(A_{7}\)}}
\put(01,02){\makebox(0,0)[bl]{\(A_{8}\)}}

\put(1.5,39){\makebox(0,0)[bl]{Per-}}
\put(1.5,36){\makebox(0,0)[bl]{son}}

\put(12,37){\makebox(0,0)[bl]{\(K_{1}\)}}

\put(22,37){\makebox(0,0)[bl]{\(K_{2}\)}}

\put(32,37){\makebox(0,0)[bl]{\(K_{3}\)}}

\put(42,37){\makebox(0,0)[bl]{\(K_{4}\)}}

%--------------------

%A1
\put(14,30){\makebox(0,0)[bl]{\(0\)}}
\put(24,30){\makebox(0,0)[bl]{\(5\)}}
\put(34,30){\makebox(0,0)[bl]{\(2\)}}
\put(44,30){\makebox(0,0)[bl]{\(3\)}}

%A2
\put(14,26){\makebox(0,0)[bl]{\(5\)}}
\put(24,26){\makebox(0,0)[bl]{\(2\)}}
\put(34,26){\makebox(0,0)[bl]{\(3\)}}
\put(44,26){\makebox(0,0)[bl]{\(3\)}}

%A3
\put(14,22){\makebox(0,0)[bl]{\(4\)}}
\put(24,22){\makebox(0,0)[bl]{\(3\)}}
\put(34,22){\makebox(0,0)[bl]{\(1\)}}
\put(44,22){\makebox(0,0)[bl]{\(2\)}}

%A4
\put(14,18){\makebox(0,0)[bl]{\(4\)}}
\put(24,18){\makebox(0,0)[bl]{\(3\)}}
\put(34,18){\makebox(0,0)[bl]{\(4\)}}
\put(44,18){\makebox(0,0)[bl]{\(2\)}}

%A5
\put(14,14){\makebox(0,0)[bl]{\(3\)}}
\put(24,14){\makebox(0,0)[bl]{\(5\)}}
\put(34,14){\makebox(0,0)[bl]{\(3\)}}
\put(44,14){\makebox(0,0)[bl]{\(5\)}}

%A6
\put(14,10){\makebox(0,0)[bl]{\(1\)}}
\put(24,10){\makebox(0,0)[bl]{\(5\)}}
\put(34,10){\makebox(0,0)[bl]{\(2\)}}
\put(44,10){\makebox(0,0)[bl]{\(5\)}}

%A7
\put(14,06){\makebox(0,0)[bl]{\(3\)}}
\put(24,06){\makebox(0,0)[bl]{\(3\)}}
\put(34,06){\makebox(0,0)[bl]{\(5\)}}
\put(44,06){\makebox(0,0)[bl]{\(5\)}}

%A8
\put(14,02){\makebox(0,0)[bl]{\(3\)}}
\put(24,02){\makebox(0,0)[bl]{\(3\)}}
\put(34,02){\makebox(0,0)[bl]{\(4\)}}
\put(44,02){\makebox(0,0)[bl]{\(4\)}}

\end{picture}
\end{center}

\begin{center}
\begin{picture}(90,67.5)
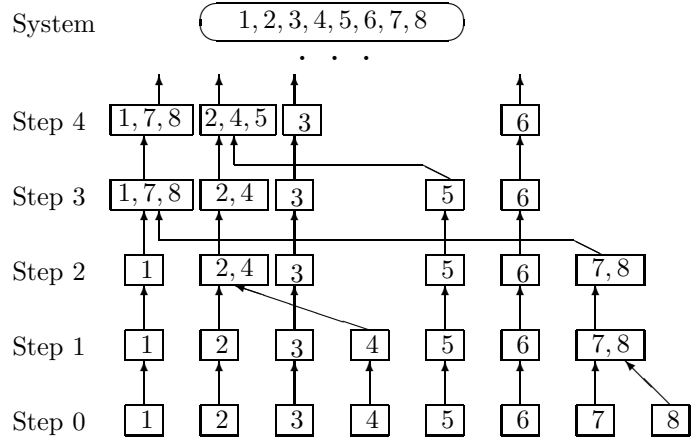
\put(08.5,00){\makebox(0,0)[bl]{Fig. 12. Example of basic
hierarchical clustering}}

\put(00,07){\makebox(0,0)[bl]{Step \(0\)}}
%-----------------------------------------
\put(17,08){\makebox(0,0)[bl]{\(1\)}}

\put(15,11){\line(1,0){5}} \put(15,07){\line(1,0){5}}
\put(15,07){\line(0,1){4}} \put(20,07){\line(0,1){4}}

%--------------------
\put(27,08){\makebox(0,0)[bl]{\(2\)}}

\put(25,11){\line(1,0){5}} \put(25,07){\line(1,0){5}}
\put(25,07){\line(0,1){4}} \put(30,07){\line(0,1){4}}
%-------------------
\put(37,08){\makebox(0,0)[bl]{\(3\)}}

\put(35,11){\line(1,0){5}} \put(35,7){\line(1,0){5}}
\put(35,07){\line(0,1){4}} \put(40,7){\line(0,1){4}}
%-------------------
\put(47,8){\makebox(0,0)[bl]{\(4\)}}

\put(45,11){\line(1,0){5}} \put(45,7){\line(1,0){5}}
\put(45,7){\line(0,1){4}} \put(50,7){\line(0,1){4}}
%--------------------------
\put(57,8){\makebox(0,0)[bl]{\(5\)}}

\put(55,11){\line(1,0){5}} \put(55,7){\line(1,0){5}}
\put(55,7){\line(0,1){4}} \put(60,7){\line(0,1){4}}
%--------------------------
\put(67,8){\makebox(0,0)[bl]{\(6\)}}

\put(65,11){\line(1,0){5}} \put(65,7){\line(1,0){5}}
\put(65,7){\line(0,1){4}} \put(70,7){\line(0,1){4}}
%-----------------------
\put(77,8){\makebox(0,0)[bl]{\(7\)}}

\put(75,11){\line(1,0){5}} \put(75,7){\line(1,0){5}}
\put(75,7){\line(0,1){4}} \put(80,7){\line(0,1){4}}

%-----------------------
\put(87,8){\makebox(0,0)[bl]{\(8\)}}

\put(85,11){\line(1,0){5}} \put(85,7){\line(1,0){5}}
\put(85,7){\line(0,1){4}} \put(90,7){\line(0,1){4}}
%====================================================

\put(00,17){\makebox(0,0)[bl]{Step \(1\)}}
%-----------------------------------------
\put(17,18){\makebox(0,0)[bl]{\(1\)}}

\put(15,21){\line(1,0){5}} \put(15,17){\line(1,0){5}}
\put(15,17){\line(0,1){4}} \put(20,17){\line(0,1){4}}
\put(17.5,11){\vector(0,1){6}}
%--------------------
\put(27,18){\makebox(0,0)[bl]{\(2\)}}

\put(25,21){\line(1,0){5}} \put(25,17){\line(1,0){5}}
\put(25,17){\line(0,1){4}} \put(30,17){\line(0,1){4}}
\put(27.5,11){\vector(0,1){6}}
%-------------------
\put(37,17.5){\makebox(0,0)[bl]{\(3\)}}

\put(35,21){\line(1,0){5}} \put(35,17){\line(1,0){5}}
\put(35,17){\line(0,1){4}} \put(40,17){\line(0,1){4}}
\put(37.5,11){\vector(0,1){6}}
%-------------------
\put(47,18){\makebox(0,0)[bl]{\(4\)}}

\put(45,21){\line(1,0){5}} \put(45,17){\line(1,0){5}}
\put(45,17){\line(0,1){4}} \put(50,17){\line(0,1){4}}
\put(47.5,11){\vector(0,1){6}}
%--------------------------
\put(57,18){\makebox(0,0)[bl]{\(5\)}}

\put(55,21){\line(1,0){5}} \put(55,17){\line(1,0){5}}
\put(55,17){\line(0,1){4}} \put(60,17){\line(0,1){4}}
\put(57.5,11){\vector(0,1){6}}
%--------------------------
\put(67,18){\makebox(0,0)[bl]{\(6\)}}

\put(65,21){\line(1,0){5}} \put(65,17){\line(1,0){5}}
\put(65,17){\line(0,1){4}} \put(70,17){\line(0,1){4}}
\put(67.5,11){\vector(0,1){6}}
%-----------------------
\put(77,17.5){\makebox(0,0)[bl]{\(7,8\)}}

\put(75,21){\line(1,0){9}} \put(75,17){\line(1,0){9}}
\put(75,17){\line(0,1){4}} \put(84,17){\line(0,1){4}}

\put(77.5,11){\vector(0,1){6}} \put(87.5,11){\vector(-1,1){6}}
%====================================================

\put(00,27){\makebox(0,0)[bl]{Step \(2\)}}
%-----------------------------------------
\put(17,28){\makebox(0,0)[bl]{\(1\)}}

\put(15,31){\line(1,0){5}} \put(15,27){\line(1,0){5}}
\put(15,27){\line(0,1){4}} \put(20,27){\line(0,1){4}}
\put(17.5,21){\vector(0,1){6}}
%--------------------
\put(27,27.5){\makebox(0,0)[bl]{\(2,4\)}}

\put(25,31){\line(1,0){9}} \put(25,27){\line(1,0){9}}
\put(25,27){\line(0,1){4}} \put(34,27){\line(0,1){4}}
\put(27.5,21){\vector(0,1){6}}

\put(47.5,21){\vector(-3,1){18}}
%-------------------
\put(37,27.5){\makebox(0,0)[bl]{\(3\)}}

\put(35,31){\line(1,0){5}} \put(35,27){\line(1,0){5}}
\put(35,27){\line(0,1){4}} \put(40,27){\line(0,1){4}}
\put(37.5,21){\vector(0,1){6}}

%--------------------------
\put(57,28){\makebox(0,0)[bl]{\(5\)}}

\put(55,31){\line(1,0){5}} \put(55,27){\line(1,0){5}}
\put(55,27){\line(0,1){4}} \put(60,27){\line(0,1){4}}
\put(57.5,21){\vector(0,1){6}}
%--------------------------
\put(67,27.5){\makebox(0,0)[bl]{\(6\)}}

\put(65,31){\line(1,0){5}} \put(65,27){\line(1,0){5}}
\put(65,27){\line(0,1){4}} \put(70,27){\line(0,1){4}}
\put(67.5,21){\vector(0,1){6}}
%-----------------------
\put(77,27.5){\makebox(0,0)[bl]{\(7,8\)}}

\put(75,31){\line(1,0){9}} \put(75,27){\line(1,0){9}}
\put(75,27){\line(0,1){4}} \put(84,27){\line(0,1){4}}
\put(77.5,21){\vector(0,1){6}}
%====================================================

\put(00,37){\makebox(0,0)[bl]{Step \(3\)}}
%-----------------------------------------
\put(14,37.5){\makebox(0,0)[bl]{\(1,7,8\)}}

\put(13,41){\line(1,0){11}} \put(13,37){\line(1,0){11}}
\put(13,37){\line(0,1){4}} \put(24,37){\line(0,1){4}}
\put(17.5,31){\vector(0,1){6}}

 \put(19.5,33){\vector(0,1){4}}

 \put(74.5,33){\line(-1,0){55}}

 \put(78.5,31){\line(-2,1){4}}

%
%--------------------
\put(27,37.5){\makebox(0,0)[bl]{\(2,4\)}}

\put(25,41){\line(1,0){9}} \put(25,37){\line(1,0){9}}
\put(25,37){\line(0,1){4}} \put(34,37){\line(0,1){4}}
\put(27.5,31){\vector(0,1){6}}
%-------------------
\put(37,37.5){\makebox(0,0)[bl]{\(3\)}}

\put(35,41){\line(1,0){5}} \put(35,37){\line(1,0){5}}
\put(35,37){\line(0,1){4}} \put(40,37){\line(0,1){4}}
\put(37.5,31){\vector(0,1){6}}
%-------------------
%--------------------------
\put(57,38){\makebox(0,0)[bl]{\(5\)}}

\put(55,41){\line(1,0){5}} \put(55,37){\line(1,0){5}}
\put(55,37){\line(0,1){4}} \put(60,37){\line(0,1){4}}
\put(57.5,31){\vector(0,1){6}}
%--------------------------
\put(67,37.5){\makebox(0,0)[bl]{\(6\)}}

\put(65,41){\line(1,0){5}} \put(65,37){\line(1,0){5}}
\put(65,37){\line(0,1){4}} \put(70,37){\line(0,1){4}}
\put(67.5,31){\vector(0,1){6}}
%-----------------------

%====================================================

\put(00,47){\makebox(0,0)[bl]{Step \(4\)}}
%-----------------------------------------
\put(14,47.5){\makebox(0,0)[bl]{\(1,7,8\)}}

\put(13,51){\line(1,0){11}} \put(13,47){\line(1,0){11}}
\put(13,47){\line(0,1){4}} \put(24,47){\line(0,1){4}}

\put(17.5,41){\vector(0,1){6}}

\put(19.5,51){\vector(0,1){4}}
%--------------------
\put(25.5,47.5){\makebox(0,0)[bl]{\(2,4,5\)}}

\put(25,51){\line(1,0){10}} \put(25,47){\line(1,0){10}}
\put(25,47){\line(0,1){4}} \put(35,47){\line(0,1){4}}
\put(27.5,41){\vector(0,1){6}}

\put(29.5,43){\vector(0,1){4}}

 \put(54.5,43){\line(-1,0){25}}

 \put(58.5,41){\line(-2,1){4}}

\put(27.5,51){\vector(0,1){4}}
%-------------------
\put(38,47.5){\makebox(0,0)[bl]{\(3\)}}

\put(36,51){\line(1,0){5}} \put(36,47){\line(1,0){5}}
\put(36,47){\line(0,1){4}} \put(41,47){\line(0,1){4}}
\put(37.5,41){\vector(0,1){6}}

\put(37.5,51){\vector(0,1){4}}
%-------------------

%--------------------------
\put(67,47.5){\makebox(0,0)[bl]{\(6\)}}

\put(65,51){\line(1,0){5}} \put(65,47){\line(1,0){5}}
\put(65,47){\line(0,1){4}} \put(70,47){\line(0,1){4}}
\put(67.5,41){\vector(0,1){6}}

\put(67.5,51){\vector(0,1){4}}

%-----------------------
%====================================================
\put(00,60){\makebox(0,0)[bl]{System}}

%\put(00,67){\makebox(0,0)[bl]{Step \(6\)}}
%-----------------------------------------

\put(42.5,62){\oval(35,5)}

%--------------------
\put(30,60.5){\makebox(0,0)[bl]{\(1,2,3,4,5,6,7,8\)}}

\put(38,57){\makebox(0,0)[bl]{{\bf . ~. ~.}}}

\end{picture}
\end{center}

 Evidently, modifications of the agglomerative algorithms have
 been proposed
 (e.g., \cite{fern08}, \cite{szekely05}).
 Our approach to modifications of the agglomerative
 algorithm above is the following \cite{lev07clust}.

 First, computing the item ``distance'' based on metric \(l_{2}\)
 is a very 'simplified' mathematical approach.
 It is possible to examine an ordinal sale for the item ``proximity''.
 As a result,
  it is possible to select the smallest ordinal item ``proximity''
 and integrate corresponding item pairs.

 Second, multicriteria approach for the item ``proximity'',
 (e.g., Pareto approach \cite{pareto71}) can be used.

 Third, clustering can be organized as
 a series revelation of cliques \cite{gar79}.

 In addition,
 it is necessary to note
 clustering processes can be examined as optimization problems.
 The basic approach to goal functions in clustering is based on the following:
  (1) inter-cluster ``distances'' (i.e., proximity between elements in the same clusters),
  (2) intra-cluster ``distances'' (i.e., proximity between elements in different clusters).
 Some other objectives can be used as well
 (e.g.,
  number of clusters or closeness to a required interval of cluster numbers,
  cluster's cardinalities).
 Some other objective functions can be used as well
 (e.g., \cite{hopp99}, \cite{jon01}).
 Thus multicriteria problem formulations based on criteria above
 may be examined as well.
 Generally, it may be very prospective to integrate
  hierarchical clustering and
   expert-based interactive procedures
  (previous section)
   for the design of system hierarchies.

%%%%%%%%%%%%%%%%%%%%%%%%%%%%%%%%%%%%%
\subsection{Towards Ontology}

 In recent two decades, ontology based approaches have been widely
 used for representation and processing of information
 in various domains, for example
(e.g.,
 \cite{corcho03},\cite{gang05}, \cite{gruber93},
 \cite{hol02},\cite{lopez99},\cite{mussen92},\cite{noy97},\cite{usch96}):
  knowledge-based systems,
  system design,
  systems engineering,
  library science,
  chemistry,
  biomedical informatics,
  conceptual modeling,
  semantic Web.
 Ontologies are artifacts as structures
 (logical, linguistic, ``taxonomical'')
 for description of a domain and/or a ``space'' of tasks
 (e.g., \cite{noy97},\cite{usch96}).
 In fact, hierarchical (multi-layer) structures are used here.
 Basic problems over ontologies are the following
 (e.g., \cite{corcho03},\cite{hol02},\cite{noy97},\cite{noy00},\cite{noy05},\cite{usch96}):
  (a) design,
  (b) comparison,
  (c) integration/merging,
  (d) alignment.
 Ontology is often presented as the following activity over various patterns
 (e.g., logical, reasoning, architectural, naming, content):
 searching, selecting, composing.
 A framework for ontology design has been suggested in
 \cite{noy97}.
 An approach to ontology extraction was presented
  in \cite{gaeta11}:

 ~~

 {\it Stage 1.} Preprocessing (a preliminary work on the available documents is carried out).

 {\it Stage 2.} Creation of the first version of the ontology (as a structure).

 {\it Stage 3.} Creation of concept and relationship.
 (The creation of the whole ontology extracting the concepts and
 their relations from the text documents is carried out).

 {\it Stage 4.} Harmonization.
 The extracted ontology is ``harmonized''
 through the analysis of other domain ontologies and concept
 description from other systems (e.g., Wikipedia).

 {\it Stage 5.} Refinement and validation.
 The resultant
 ontology is refined and validated.

~~

 The following main resources for ontologies design are usually pointed out:
 informal data structures,
 concept schemes (e.g., classifications, thesauri, nomenclatures),
 Web-based resources,
 natural language documents,
 lexical resources (e.g., dictionaries),
 modeling languages (e.g., UML, Petri nets).
 Evidently, modular approaches (based on modular architecture)
 may be widely used in the ontology design and ontology
 ``life cycle'' (i.e., creation, evaluation, testing, utilization, modification).

%%%%%%%%%%%%%%%%%%%%%%%%%%%%%%%%%%%%%%%%%%%%%%%
\subsection{Spanning Trees}

  Evidently, spanning tree problems can be used for
 the design of hierarchical system models,
 in the case when a preliminary network system model over
 system elements exists).
 Many decades, spanning trees problems are
 used in applications (e.g., network design and maintenance,
  communication protocol design, VLSI design)
 and intensively studied in combinatorial optimization
 (e.g.,
 \cite{cormen01}, \cite{du08}, \cite{gar79}, \cite{goemans93},
 \cite{hwang92},
 \cite{knuth68}, \cite{winter87}, \cite{wu03}).
 In this section,
 three basic spanning problems (and their modifications)
 are briefly described:
 (a) minimum spanning tree problem,
 (b) minimum Steiner tree problem, and
 (c) maxim leaf spanning tree problem.

\subsubsection{Minimum Spanning Tree and Steiner Tree Problems}

 In this section, a brief description of spanning trees problems
 is presented (from \cite{lev09}).
  Table 2 contains a list of basic spanning tree problems and
 corresponding literature sources.
 Fig. 13 illustrates two spanning tree problems.

\begin{center}
\begin{picture}(45,53)
\put(26,00){\makebox(0,0)[bl]{Fig. 13. Illustration for spanning
trees  \cite{lev11ADES}}}

\put(08,48){\makebox(0,0)[bl]{Initial graph}}

%--10
\put(05,38){\circle*{1.3}}

%=======================10-9
\put(05,38){\line(1,1){5}}

%=======================9-8
\put(10,43){\line(1,0){15}}

%=======================9-4
\put(10,43){\line(1,-4){5}}

%=======================9-7
\put(10,43){\line(1,-1){10}}

%--9
\put(10,43){\circle*{1.3}}

%--8
\put(25,43){\circle*{1.3}}

%=======================8-7
\put(25,43){\line(-1,-2){5}}

%=======================8-5
\put(25,43){\line(0,-1){15}}

%--7
\put(20,33){\circle*{1.3}}

%=======================7-4
\put(20,33){\line(-1,-2){5}}

%=======================7-5
\put(20,33){\line(1,-1){5}}

%--6
\put(30,18){\circle*{1.3}}

%--6'
\put(30,33){\circle*{1.3}}

%=======================6'-6
\put(30,33){\line(0,-1){15}}

\put(30,33){\line(-1,-1){5}}

\put(30,33){\line(-1,2){5}}

%=======================5-6
\put(25,28){\line(1,-2){5}}

%--5
\put(25,28){\circle*{1.3}}

%%%10:(5,38)
%=======================4-10
\put(15,23){\line(-2,3){10}}

%=======================4-5
\put(15,23){\line(2,1){10}}

%--4
\put(15,23){\circle*{1.3}}

%=======================4-3
\put(15,23){\line(-3,-1){15}}

%=======================4-2
\put(15,23){\line(-1,-2){5}}

%=======================4-1
\put(15,23){\line(1,-3){5}}

%9:10,25
%=======================3-9
\put(00,18){\line(1,2){5}} \put(05,28){\line(1,3){5}}

%--3
\put(00,18){\circle*{1.3}}

%=======================2-6
\put(10,13){\line(4,1){20}}

%=======================2-3
\put(10,13){\line(-2,1){10}}

%--2
\put(10,13){\circle*{1.3}}

%=======================1-5
\put(20,7){\line(1,4){5}}

%--1
 \put(20,7){\circle*{1.3}}

\end{picture}
%\end{center}
%
\begin{picture}(45,53)

\put(07,48){\makebox(0,0)[bl]{Spanning tree}}

%--10
\put(05,38){\circle*{1.3}}

%=======================10-9
\put(05,38){\line(1,1){5}}

%=======================9-8
%\put(10,43){\line(1,0){15}}

%=======================9-4
%\put(10,43){\line(1,-4){5}}

%=======================9-7
\put(10,43){\line(1,-1){10}}

%--9
\put(10,43){\circle*{1.3}}

%--8
\put(25,43){\circle*{1.3}}

%=======================8-7
\put(25,43){\line(-1,-2){5}}

%=======================8-5
%\put(25,43){\line(0,-1){15}}

%--7
\put(20,33){\circle*{1.3}}

%=======================7-4
\put(20,33){\line(-1,-2){5}}

%=======================7-5
\put(20,33){\line(1,-1){5}}

%--6
\put(30,18){\circle*{1.3}}

%--6'
\put(30,33){\circle*{1.3}}

%=======================6'-X
\put(30,33){\line(-1,-1){5}}

%=======================5-6
\put(25,28){\line(1,-2){5}}

%--5
\put(25,28){\circle*{1.3}}

%%%10:(5,38)
%=======================4-10
%\put(15,23){\line(-2,3){10}}

%=======================4-5
%\put(15,23){\line(2,1){10}}

%--4
\put(15,23){\circle*{1.3}}

%=======================4-3
\put(15,23){\line(-3,-1){15}}

%=======================4-2
\put(15,23){\line(-1,-2){5}}

%--3
\put(00,18){\circle*{1.3}}

%--2
\put(10,13){\circle*{1.3}}

%=======================1-5
\put(20,7){\line(1,4){5}}

%--1
 \put(20,7){\circle*{1.3}}

\end{picture}
%\end{center}
%
%
\begin{picture}(30,53)

\put(08.5,48){\makebox(0,0)[bl]{Steiner tree}}
\put(02,44){\makebox(0,0)[bl]{(4 Steiner vertices)}}

%--10
\put(05,38){\circle*{1.3}}

%=======================10-9
\put(05,38){\line(1,1){5}}

%--9
\put(10,43){\circle*{1.3}}

%--8
\put(25,43){\circle*{1.3}}

%%%%%%%%%%%%%%%%%%%%%%%%%%%%%%%%%%%%%%%%%%%Steiner point 7-8-9
\put(20,38){\circle*{1.0}} \put(20,38){\circle{2.0}}

%=======================7-8-9-> 9
\put(20,38){\line(-2,1){10}} \put(20,38){\line(1,1){5}}
\put(20,38){\line(0,-1){5}}
%%%%%%%%%%%%%%%%%%%%%%%%%%%%%%%%%%%%%%%%%%%%%%%%

%--7
\put(20,33){\circle*{1.3}}

%--6
\put(30,18){\circle*{1.3}}

%--6'
\put(30,33){\circle*{1.3}}

%=======================6'-X
\put(30,33){\line(-1,-1){5}}

%%%%%%%%%%%%%%%%%%%%%%%%%%%%%%%%%%%%%%%%%%% Steiner point 1-5-6
\put(25,23){\circle*{1.0}} \put(25,23){\circle{2.0}}

%=======================1-5-6->
\put(25,23){\line(0,1){5}} \put(25,23){\line(1,-1){5}}
\put(25,23){\line(-1,-3){5}}
%%%%%%%%%%%%%%%%%%%%%%%%%%%%%%%%%%%%%%%%%%%%%%%%

%--5
\put(25,28){\circle*{1.3}}

%%%%%%%%%%%%%%%%%%%%%%%%%%%%%%%%%%%%%%%%%%% Steiner point 1-4-5
\put(20,28){\circle*{1.0}} \put(20,28){\circle{2.0}}

%=======================1-4-5->
\put(20,28){\line(0,1){5}} \put(20,28){\line(1,0){5}}
\put(20,28){\line(-1,-1){5}}
%%%%%%%%%%%%%%%%%%%%%%%%%%%%%%%%%%%%%%%%%%%%%%%%

%%%%%%%%%%%%%%%%%%%%%%%%%%%%%%%%%%%%%%%%%%%Steiner point 4-7-5
\put(20,38){\circle*{1.0}} \put(20,38){\circle{2.0}}

%=======================4-7-5->
\put(20,38){\line(-2,1){10}} \put(20,38){\line(1,1){5}}
\put(20,38){\line(0,-1){5}}
%%%%%%%%%%%%%%%%%%%%%%%%%%%%%%%%%%%%%%%%%%%%%%%%

%--4
\put(15,23){\circle*{1.3}}

%%%%%%%%%%%%%%%%%%%%%%%%%%%%%%%%%%%%%%%%%%% Steiner point 2-3-4
\put(10,18){\circle*{1.0}} \put(10,18){\circle{2.0}}

%=======================2-3-4-> 9
\put(10,18){\line(-1,0){10}} \put(10,18){\line(0,-1){5}}
\put(10,18){\line(1,1){5}}
%%%%%%%%%%%%%%%%%%%%%%%%%%%%%%%%%%%%%%%%%%%%%%%%

%--3
\put(00,18){\circle*{1.3}}

%--2
\put(10,13){\circle*{1.3}}

%--1
 \put(20,7){\circle*{1.3}}

\end{picture}
\end{center}

 The basic {\it spanning} problem is {\it the minimum spanning tree}
 or {\it the minimum weight spanning tree}
 problem
 (e.g.,
 \cite{cormen01}, \cite{drortree00},
 \cite{gabow86}, \cite{gar79}, \cite{pettie02}).
 Let \(G = (A,E)\) be a connected graph
 (\(A\) is
 the set of vertices,
 \(E\) is
 the set of edges/arcs)
  with nonnegative weights of
 edges/arcs.

 A spanning tree of the graph  ~\(T=(A,E')\)~ (\( E' \subseteq E \))
  is a subgraph which is a tree and
 connects all the vertices together.
 The total weight (cost) of the spanning tree
 ~\(c(T)\)~
 is the sum of weights
 of all its edges/arcs (i.e.,  \(E'\)).
 A {\it minimum spanning tree} or {\it minimum weight spanning tree}
 \(T^{*}\)
 is a spanning tree with the total weight less or equal to the
 weights of every other spanning tree
 ~~\( c(T^{*}) = \min_{\{T\}}   c(T) \).
 The problem is polynomially solvable
 (e.g.,
  \cite{cormen01}, \cite{gabow86}, \cite{gar79}, \cite{pettie02}).
 The basic well-known
 polynomial algorithms are the following
  (e.g.,  \cite{cormen01}, \cite{gar79}:
 (a) Prim's algorithm, (b) Kruskal's algorithm, (c) Boruvka's algorithm.

 In more general case, a spanning structure corresponds to {\it spanning forest}:
 any graph (not necessary connected) has
 a {\it minimum spanning forest}
 which is a union of minimum spanning trees for its connected components.
 (e.g., \cite{gar79}, \cite{pettie02a}).

\begin{center}
\begin{picture}(108,133)

\put(25,129){\makebox(0,0)[bl]{Table 2. Studies in spanning trees
}}

\put(00,00){\line(1,0){108}} \put(00,121){\line(1,0){108}}
\put(00,127){\line(1,0){108}}

\put(00,0){\line(0,1){127}} \put(74,0){\line(0,1){127}}
\put(108,0){\line(0,1){127}}

\put(01,122){\makebox(0,0)[bl]{Spanning problem}}
\put(75,122){\makebox(0,0)[bl]{Sources}}

%------------------------------------------

\put(01,115){\makebox(0,0)[bl]{1.Spanning tree}}

%--
\put(01,111){\makebox(0,0)[bl]{1.1.Minimum spanning tree}}
\put(75,111){\makebox(0,0)[bl]{\cite{cormen01},\cite{gabow86},\cite{gar79},\cite{pettie02}
}}

%--
\put(01,107){\makebox(0,0)[bl]{1.2.Minimum diameter spanning
tree}} \put(75,107){\makebox(0,0)[bl]{\cite{gfe12} }}

%--
\put(01,103){\makebox(0,0)[bl]{1.3.Minimum spanning forest}}
\put(75,103){\makebox(0,0)[bl]{\cite{gar79},\cite{pettie02a} }}

%--
\put(01,99){\makebox(0,0)[bl]{1.4.Minimum spanning multi-tree}}
\put(75,99){\makebox(0,0)[bl]{\cite{geo12},\cite{itai89},\cite{tar76}
}}

%--
\put(01,95){\makebox(0,0)[bl]{1.5.Multicriteria spanning tree}}
\put(75,95){\makebox(0,0)[bl]{\cite{arroyo08},\cite{chen07},\cite{dell00},\cite{hamaher95}
}}

%--
\put(01,90.7){\makebox(0,0)[bl]{2.Steiner tree}}
%\put(75,90){\makebox(0,0)[bl]{\cite{} }}

%--
\put(01,86){\makebox(0,0)[bl]{2.1.Bottleneck Steiner tree
problem}}

\put(75,86){\makebox(0,0)[bl]{\cite{bae10},\cite{du01},\cite{du08}}}

%--
\put(01,82){\makebox(0,0)[bl]{2.2.Steiner tree problem with
minimum }}

\put(07,78){\makebox(0,0)[bl]{number  of Steiner points}}

\put(75,82){\makebox(0,0)[bl]{\cite{chen00},
\cite{du01},\cite{du08},\cite{lin99} }}

%--
\put(01,74){\makebox(0,0)[bl]{2.3.Terminal Steiner tree problem}}
\put(75,74){\makebox(0,0)[bl]{\cite{drake04},\cite{fuchs03},\cite{lin00}
}}

%--
\put(01,70){\makebox(0,0)[bl]{2.4.Node weighted Steiner tree
problem}}

\put(75,70){\makebox(0,0)[bl]{\cite{guha99},\cite{klein95},\cite{segev87},\cite{zou08}
}}

%--
\put(01,66){\makebox(0,0)[bl]{2.5.Prize-collecting Steiner
problem}}

\put(75,66){\makebox(0,0)[bl]{\cite{gutner08},\cite{hwang92},\cite{lju06},\cite{uchoa06}}}

%--
\put(01,62){\makebox(0,0)[bl]{2.6.Survivable Steiner network
problem }}

\put(75,62){\makebox(0,0)[bl]{\cite{du08} }}

%--
\put(01,58){\makebox(0,0)[bl]{2.7.Steiner tree coloring problem}}
\put(75,58){\makebox(0,0)[bl]{\cite{du08} }}

%--
\put(01,54){\makebox(0,0)[bl]{2.8.Steiner tree scheduling
problem}}

\put(75,54){\makebox(0,0)[bl]{\cite{du08} }}

%--
\put(01,50){\makebox(0,0)[bl]{2.9.Constrained Steiner tree
problem}}

\put(75,50){\makebox(0,0)[bl]{\cite{costa08},\cite{rosen95} }}

%--
\put(01,46){\makebox(0,0)[bl]{2.10.Steiner tree problem with hop
constraints}}

\put(75,46){\makebox(0,0)[bl]{\cite{costa08},\cite{voss99} }}

%--
\put(01,42){\makebox(0,0)[bl]{2.11.Steiner tree problems with
profits}}

\put(75,42){\makebox(0,0)[bl]{\cite{costa06} }}

%--
\put(01,38){\makebox(0,0)[bl]{2.12.Generalized Steiner problem}}
\put(75,38){\makebox(0,0)[bl]{\cite{agra95},\cite{awer04},\cite{dror00},\cite{drortree00},\cite{west93}}}

%--
\put(01,34){\makebox(0,0)[bl]{2.13.Generalized Steiner star
problem}}

\put(75,34){\makebox(0,0)[bl]{\cite{khu02} }}

%--
\put(01,30){\makebox(0,0)[bl]{2.14.Stochastic Steiner tree
problem}}

\put(75,30){\makebox(0,0)[bl]{\cite{gupta05},\cite{gupta07a} }}

%--
\put(01,26){\makebox(0,0)[bl]{2.15.Dynamic Steiner tree problem}}
\put(75,26){\makebox(0,0)[bl]{\cite{imase91} }}

%--
\put(01,22){\makebox(0,0)[bl]{2.16.On-line Steiner tree problem}}
\put(75,22){\makebox(0,0)[bl]{\cite{alon06},\cite{awer04},\cite{west93}}}

%--
\put(01,18){\makebox(0,0)[bl]{2.17.Group Steiner tree problem}}
\put(75,18){\makebox(0,0)[bl]{\cite{chekuri06},\cite{duin04}}}

%--
\put(01,14){\makebox(0,0)[bl]{2.18.Steiner forest problem }}

\put(75,14){\makebox(0,0)[bl]{\cite{feld09} }}

%--
\put(01,10){\makebox(0,0)[bl]{2.19.Multicriteria Steiner tree
problem}}

\put(75,10){\makebox(0,0)[bl]{\cite{levnur11},\cite{levzam11},\cite{vuj03}
}}

%--
\put(01,06){\makebox(0,0)[bl]{2.20.Multicriteria Steiner tree
problem }}

\put(08,02.5){\makebox(0,0)[bl]{with the cost of Steiner
vertices}}

\put(75,06){\makebox(0,0)[bl]{\cite{levzam11}}}

\end{picture}
\end{center}

   In recent years,  {\it multicriteria spanning tree} problems
 (or {\it multi-objective spanning tree} problems)
 are examined
 (e.g., \cite{arroyo08},\cite{chen07}, \cite{dell00},
  \cite{hamaher95}).
 Here a vector-weight corresponds to each graph edge/arc
 and the objective function is a vector-like one where
 each component of the vector above equals a sum of corresponding components of
 the edges/arcs weights:
 ~\( \overline{c} (T)\).
 In this case, it is reasonable to search for Pareto-efficient
 (by ~\( \overline{c} (T)\))
 solutions.
 This problem is NP-hard (even in bicriteria case)
 (e.g., \cite{arroyo08},\cite{chen07}).
 Here the following algorithms are used (e.g.,  \cite{arroyo08}, \cite{chen07}):
 (i) multicriteria Prim's spanning tree algorithm,
 (ii) multicriteria Kruskal's spanning tree algorithm,
 (iii) genetic algorithms,
 (iv) multiobjective evolutionary optimization algorithms,
 (v) heuristics (e.g., GRASP algorithm, simulated annealing
 algorithm), and
  (vi) knowledge-based approaches.

 In the case of Steiner tree problem
 (e.g.,
 \cite{agra95},\cite{awer04},\cite{dror00},\cite{drortree00},
 \cite{du08},\cite{gar79},\cite{goemans93}, \cite{hwang92},
 \cite{klein95},\cite{segev87},\cite{west93}\cite{winter87}),
 the spanning tree can
 include additional vertices (i.e., Steiner vertices).
 Thus, the total weight (cost) of
  the resultant spanning structure can be less than in the case of
  the basic spanning tree problem.
%-------------------------------------------------------------
%
  The basic formulation of  {\it Steiner tree problem} is the following.
  Given an undirected connected graph
  \( G=(A,E) \) with
  node/vertices set \(A\),
  edge set \(E\),
  and nonnegative weights associated with the edges.
  Given a set
 ~\( Q \)
  of specified vertices
 (terminals/basic vertices, i.e., Steiner points).
 The problem is:

~~

 {\it Find a minimum cost subgraph}
 ~\( T_{s} = (Q,E') \subseteq  G \)
  ~(\( E' \subseteq E \))
  {\it such that there exists a path in the subgraph}
  ~\( T_{s} \)
 {\it between every pair of basic vertices.}

~~

  Here the optimal solution
  ~\( T_{s} \)
  is a tree  (Steiner tree).
  The Steiner tree problem is NP-hard.
 Descriptions of many Steiner tree problem formulations are presented in
 \cite{goemans93}.
 Mainly, the following solving approaches are used for the Steiner tree problems:
 (1) exact algorithms (enumerative algorithms, for example, branch-and-bound, branch-and-cut)
  (e.g., \cite{chorpa92}, \cite{lju06})
  and
  (2) various heuristics
 (e.g., \cite{costa08}, \cite{voss92}, \cite{voss00})
 including the following:
  (a) fast (greedy) algorithms
  (e.g., \cite{chekuri06},  \cite{costa08});
  (b) approximation algorithms
      (e.g.,
  \cite{agra95}, \cite{du08},
  \cite{guha99}, \cite{gupta07a}, \cite{gutner08}, \cite{klein95},
  \cite{zou08});
  (c) genetic algorithms
 (e.g., \cite{kapsalis93});
  (d) AI-based methods (e.g., \cite{joob86});
  (e) dual heuristics
  (e.g., \cite{de01});
 (f) distributed primal-dual heuristic
  (e.g., \cite{santos07});
  and
  (g) local search techniques (e.g., \cite{canuto01}).

 In  {\it multicriteria Steiner tree} problems
 (e.g., \cite{levnur11},\cite{levzam11},\cite{vuj03})
 a vector weight for each edge/arc is under consideration, and
 the vector-like objective function (as a vector of corresponding sums)
 can be used (in a simplest case):
 ~\( \overline{c} (T_{s}) \).
 Thus,
 Pareto-efficient solutions
 (by ~\( \overline{c} (T_{s})\)) are searched for
 (e.g., \cite{levnur11},\cite{levzam11}).
 Here, the following algorithmic approaches can be pointed out:
  ~(i) basic heuristics
   (e.g., \cite{vuj03});
 ~(ii) special multi-layer macro-heuristics:
   (a) partition-synthesis heuristic (e.g., \cite{lev98}),
   (b) spanning-tree based heuristic (e.g., \cite{levnur11}),
 (c) composite multistage solving scheme
 (e.g., \cite{levzam11}).

\subsubsection{Maximum Leaf Spanning Tree Problem}

 The ``maximum leaf spanning tree'' problem is the following
 (e.g.,
 \cite{alon09},
 \cite{gar79}, \cite{klei91}):

~~~

 {\it Find a spanning tree of an input graph
 so that the number of the tree leafs is maximal.}

~~

 Generally,
 the spanning tree of a graph contains the following types of nodes:
 (a) root,
 (b) internal nodes
  (the internal nodes may be considered as a virtual
 ``bus'' in networking),
  and
 (c) leaf nodes.
 Thus, the problem consists in maximizing the number of leaf nodes
 or
 minimizing the number of internal nodes.
 The problem is one of the basic NP-hard problems \cite{gar79}.
 Fig. 14 depicts an illustrative numerical example.

\begin{center}
\begin{picture}(45,52)
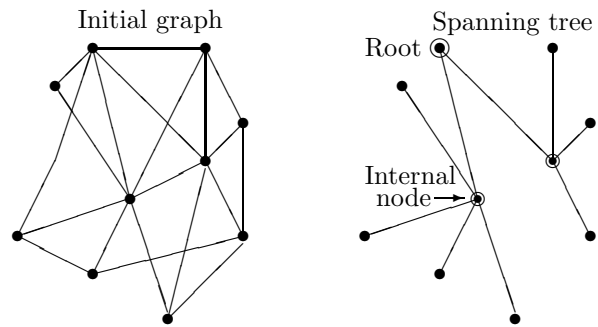
\put(00,00){\makebox(0,0)[bl]{Fig. 14. Illustration for maximum
 leaf spanning tree}}

\put(08,45){\makebox(0,0)[bl]{Initial graph}}

%--10
\put(05,38){\circle*{1.3}}

%=======================10-9
\put(05,38){\line(1,1){5}}

%=======================9-8
\put(10,43){\line(1,0){15}}

%=======================9-4
\put(10,43){\line(1,-4){5}}

%=======================9-7
\put(10,43){\line(1,-1){10}}

%--9
\put(10,43){\circle*{1.3}}

%--8
\put(25,43){\circle*{1.3}}

%=======================8-7
\put(25,43){\line(-1,-2){5}}

%=======================8-5
\put(25,43){\line(0,-1){15}}

%--7
%\put(20,33){\circle*{1.3}}

%=======================7-4
\put(20,33){\line(-1,-2){5}}

%=======================7-5
\put(20,33){\line(1,-1){5}}

%--6
\put(30,18){\circle*{1.3}}

%--6'
\put(30,33){\circle*{1.3}}

%=======================6'-6
\put(30,33){\line(0,-1){15}}

\put(30,33){\line(-1,-1){5}}

\put(30,33){\line(-1,2){5}}

%=======================5-6
\put(25,28){\line(1,-2){5}}

%--5
\put(25,28){\circle*{1.3}}

%%%10:(5,38)
%=======================4-10
\put(15,23){\line(-2,3){10}}

%=======================4-5
\put(15,23){\line(2,1){10}}

%--4
\put(15,23){\circle*{1.3}}

%=======================4-3
\put(15,23){\line(-3,-1){15}}

%=======================4-2
\put(15,23){\line(-1,-2){5}}

%=======================4-1
\put(15,23){\line(1,-3){5}}

%9:10,25
%=======================3-9
\put(00,18){\line(1,2){5}} \put(05,28){\line(1,3){5}}

%--3
\put(00,18){\circle*{1.3}}

%=======================2-6
\put(10,13){\line(4,1){20}}

%=======================2-3
\put(10,13){\line(-2,1){10}}

%--2
\put(10,13){\circle*{1.3}}

%=======================1-5
\put(20,7){\line(1,4){5}}

%--1
\put(20,7){\circle*{1.3}}

%=======================1-6
\put(20,7){\line(1,1){10}}

\end{picture}
%\end{center}
%%%%%%%%%%%%%%%%%%%%%%%%%%%%%%%%%%%%%%%%%%%%%%%%%%
\begin{picture}(35,50)

\put(09,44.8){\makebox(0,0)[bl]{Spanning tree}}

%--10
\put(05,38){\circle*{1.3}}

%=======================9-4
\put(10,43){\line(1,-4){5}}

%=======================9-7
\put(10,43){\line(1,-1){10}}

\put(00,42){\makebox(0,0)[bl]{Root}}
%--9
\put(10,43){\circle*{1.4}} \put(10,43){\circle{2.5}}

%--8
\put(25,43){\circle*{1.3}}

%=======================8-5
\put(25,43){\line(0,-1){15}}

%=======================7-5
\put(20,33){\line(1,-1){5}}

%--6
\put(30,18){\circle*{1.3}}

%--6'
\put(30,33){\circle*{1.3}}

\put(30,33){\line(-1,-1){5}}

%=======================5-6
\put(25,28){\line(1,-2){5}}

%--5
\put(25,28){\circle*{1.0}} \put(25,28){\circle{1.8}}

%%%10:(5,38)
%=======================4-10
\put(15,23){\line(-2,3){10}}

%--4
\put(15,23){\circle*{1.0}} \put(15,23){\circle{1.8}}

\put(00,25){\makebox(0,0)[bl]{Internal}}
\put(01.5,22){\makebox(0,0)[bl]{node}}

\put(9.2,23.1){\vector(1,0){4}}

%=======================4-3
\put(15,23){\line(-3,-1){15}}

%=======================4-2
\put(15,23){\line(-1,-2){5}}

%=======================4-1
\put(15,23){\line(1,-3){5}}

%--3
\put(00,18){\circle*{1.3}}

%--2
\put(10,13){\circle*{1.3}}

%--1
\put(20,7){\circle*{1.3}}

\end{picture}
\end{center}

 The problem is  applied in networking
  (e.g., minimum-energy broadcast trees in ad hoc wireless
  networks),
 (e.g,
 \cite{liang02}, \cite{thai07}),
 in circuit layout \cite{storer81}.
 The following main types of algorithms have been suggested:
 (1) exact algorithms (branch-and bound)
 (\cite{fernau11}, \cite{fuj03}),
 (2) local optimization \cite{lu92},
 (3) greedy 3-approximation algorithm (linear time) \cite{lu98},
 (4) 2-approximation algorithm \cite{solis98},
 and
%
% () approximation algorithms
%
 (5) heuristics
 (e.g., Bee Colony algorithm)
 (e.g., \cite{drescher10}, \cite{singh09}).
%
%%%%%%%%%%%%%%%%%%%%%%%%%%%%%%%%%%%%%%%%%%%%%%%%%%%%%%%%%%
%
 Some prospective generalizations of the
 ``maximum leaf spanning tree'' problem may be considered
 while taking into account the following:
 (a) multicriteria descriptions of leaf nodes,
 (b) uncertainty,
 (c) dynamics.

 In sense of exact algorithms, this problem is equivalent to
 ``connected dominating set'' problem
 (e.g.,
 \cite{blum05}, \cite{caro00}, \cite{cheng03}, \cite{gar79}):

~~

 Find a minimum set of vertices \(D \subseteq A\)
 of  input graph \(G = (A,E)\)
 that the induced by \(D\) subgraph \(G'=(D,E')\) (\( E' \subseteq E \))
 is connected dominated set and \(D\) is a dominating set of \(G\).

~~

  A recent survey on the problem is presented in \cite{blum05}.
 The  ``connected dominating set'' problem
 plays the central role
 in sensor wireless networks,
 in mobile ad-hoc networks (MANETs),
 in network testing
 (e.g., \cite{blum05},\cite{dai04}, \cite{liang02}, \cite{thai07}).
 This problem is used in communication protocols
 including the following \cite{blum05}:
 (i) media access coordination,
 (ii) unicast,
 (iii) multicast/broadcast,
 (iv) location-based routing,
 (v)  energy conservation,
 (vi) topology control, and
 (vii) resource discovery in MANET.
 Mainly, the following algorithms
 are used for the problem:
 (i) approximation algorithms
 (e.g., \cite{cheng03},\cite{guha98}),
 (ii) heuristics
 (e.g., \cite{butenko04},\cite{sanchis02}).
 Analogically,
 prospective generalizations of the
 problem may involve the following:
 (a) multicriteria descriptions of
 the elements of the dominating set,
 (b) uncertainty,
 (c) dynamics.

\subsection{Towards Optimal Organizational Hierarchies}

 Organizational hierarchies have been investigated many decades
 (e.g.,
  \cite{baligh06}, \cite{dean92},  \cite{gubko06},
  \cite{march93}, \cite{mishin07}, \cite{tannen74},
  \cite{voronin02}, \cite{waller76}).
 The following basic approaches to organizational structures are
 used in firms (manufacturing, sales, etc.):
 bureaucratic structures,
 functional structures,
 division structure (or product structure),
 matrix structure (integration of functions and products).
 Generally,
 the following hierarchical layers are considered for organizations:
 (a) organization,
 (b) branches,
 (c) departments (divisions),
 (d) work groups,
 (e) individuals.
 Fig. 15 depicts an example of hierarchical structure for universities.

\begin{center}
\begin{picture}(120,83)

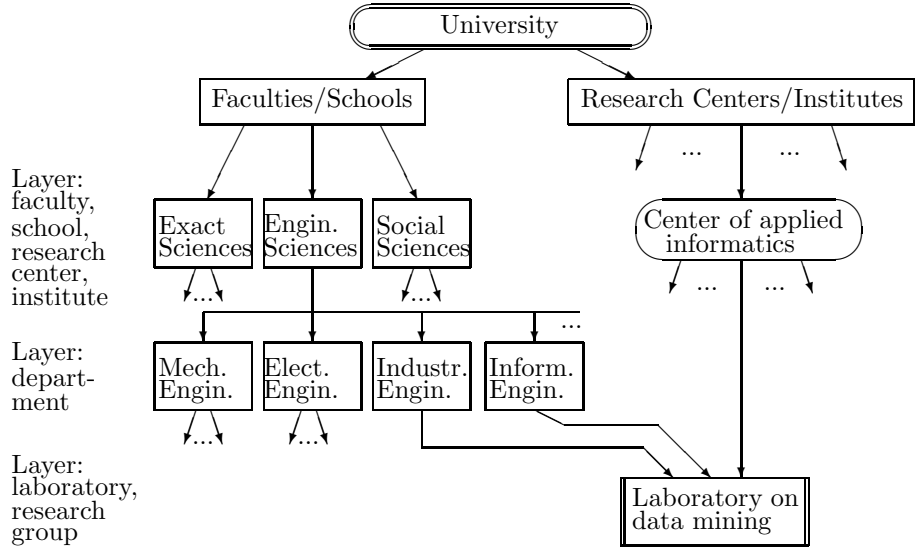
\put(16,00){\makebox(0,0)[bl]{Fig. 15. Example of hierarchical
 structure (university) }}

% 74-80
\put(65,77){\oval(40,6)} \put(65,77){\oval(39,5)}

\put(57,75.5){\makebox(0,0)[bl]{University}}

\put(55,74){\vector(-2,-1){08}} \put(75,74){\vector(2,-1){08}}

%%%%%%%%%%%%%%%%%%%%%%%%%%%%%%%%%%%%%%%%%%% 64-70
\put(26.6,65.5){\makebox(0,0)[bl]{Faculties/Schools}}

\put(25,64){\line(1,0){30}} \put(25,70){\line(1,0){30}}
\put(25,64){\line(0,1){06}} \put(55,64){\line(0,1){06}}

\put(40,64){\vector(0,-1){10}}

\put(31,64){\vector(-1,-2){05}} \put(49,64){\vector(1,-2){05}}
%--

\put(19.5,49.5){\makebox(0,0)[bl]{Exact}}
\put(19.5,46.5){\makebox(0,0)[bl]{Sciences}}

\put(19,45){\line(1,0){13}} \put(19,54){\line(1,0){13}}
\put(19,45){\line(0,1){09}} \put(32,45){\line(0,1){09}}

\put(24.5,45){\vector(-1,-3){1.5}}
\put(26.5,45){\vector(1,-3){1.5}}
\put(24,40.5){\makebox(0,0)[bl]{...}}

%--

\put(33.5,49){\makebox(0,0)[bl]{Engin.}}
\put(33.5,46.5){\makebox(0,0)[bl]{Sciences}}

\put(33.5,45){\line(1,0){13}} \put(33.5,54){\line(1,0){13}}
\put(33.5,45){\line(0,1){09}} \put(46.5,45){\line(0,1){09}}

%--

\put(48.5,49.5){\makebox(0,0)[bl]{Social}}
\put(48.5,46.5){\makebox(0,0)[bl]{Sciences}}

\put(48,45){\line(1,0){13}} \put(48,54){\line(1,0){13}}
\put(48,45){\line(0,1){09}} \put(61,45){\line(0,1){09}}

\put(53.5,45){\vector(-1,-3){1.5}}
\put(55.5,45){\vector(1,-3){1.5}}
\put(53,40.5){\makebox(0,0)[bl]{...}}

%%%%%%%%%%%%%%%%%%%% Departments 24-37

%--

\put(25.5,39){\line(1,0){50}}

\put(73,37){\makebox(0,0)[bl]{...}}

\put(25.5,39){\vector(0,-1){4}}

\put(19.5,30.5){\makebox(0,0)[bl]{Mech.}}
\put(19.5,27){\makebox(0,0)[bl]{Engin.}}

\put(19,26){\line(1,0){13}} \put(19,35){\line(1,0){13}}
\put(19,26){\line(0,1){09}} \put(32,26){\line(0,1){09}}

\put(24.5,26){\vector(-1,-3){1.5}}
\put(26.5,26){\vector(1,-3){1.5}}
\put(24,21.5){\makebox(0,0)[bl]{...}}

%--
\put(40,39){\line(0,1){6}}

\put(40,39){\vector(0,-1){4}}

\put(33.5,30.5){\makebox(0,0)[bl]{Elect.}}
\put(33.5,27){\makebox(0,0)[bl]{Engin.}}

\put(33.5,26){\line(1,0){13}} \put(33.5,35){\line(1,0){13}}
\put(33.5,26){\line(0,1){09}} \put(46.5,26){\line(0,1){09}}

\put(38.5,26){\vector(-1,-3){1.5}}
\put(40.5,26){\vector(1,-3){1.5}}
\put(38,21.5){\makebox(0,0)[bl]{...}}

%--

\put(54.5,39){\vector(0,-1){4}}

\put(48.5,30.5){\makebox(0,0)[bl]{Industr.}}
\put(48.5,27){\makebox(0,0)[bl]{Engin.}}

\put(48,26){\line(1,0){13}} \put(48,35){\line(1,0){13}}
\put(48,26){\line(0,1){09}} \put(61,26){\line(0,1){09}}

\put(54.5,26){\line(0,-1){05}} \put(54.5,21){\line(1,0){29.5}}
\put(84,21){\vector(1,-1){4}}

%--

\put(69.5,39){\vector(0,-1){4}}

\put(63.5,30.5){\makebox(0,0)[bl]{Inform.}}
\put(63.5,27){\makebox(0,0)[bl]{Engin.}}

\put(63,26){\line(1,0){13}} \put(63,35){\line(1,0){13}}
\put(63,26){\line(0,1){09}} \put(76,26){\line(0,1){09}}

\put(69.5,26){\line(2,-1){04}} \put(73.5,24){\line(1,0){12.5}}
\put(86,24){\vector(1,-1){7}}

%%%%%%%%%%%%%%%%%%%%%%%%%%%%%%%%%%%%%%%%%%%%64-70

\put(75.7,65.5){\makebox(0,0)[bl]{Research Centers/Institutes}}

\put(74,64){\line(1,0){46}} \put(74,70){\line(1,0){46}}
\put(74,64){\line(0,1){06}} \put(120,64){\line(0,1){06}}

%--

\put(85,64){\vector(-1,-3){2}} \put(109,64){\vector(1,-3){2}}

\put(89,60){\makebox(0,0)[bl]{...}}
\put(102,60){\makebox(0,0)[bl]{...}}

\put(97,64){\vector(0,-1){10}}

%---

\put(98,50){\oval(30,8)}

\put(84,49.5){\makebox(0,0)[bl]{Center of applied}}
\put(87,47){\makebox(0,0)[bl]{informatics}}

%--

\put(89,46){\vector(-1,-3){1.5}} \put(105,46){\vector(1,-3){1.5}}

\put(91,42){\makebox(0,0)[bl]{...}}
\put(100,42){\makebox(0,0)[bl]{...}}

\put(97,46){\vector(0,-1){29}}

%%%%%%%%%%%%%%%%%%--41-60

\put(00,55){\makebox(0,0)[bl]{Layer:}}
\put(00,52){\makebox(0,0)[bl]{faculty,}}
\put(00,48.5){\makebox(0,0)[bl]{school,}}
\put(00,46){\makebox(0,0)[bl]{research}}
\put(00,42.5){\makebox(0,0)[bl]{center,}}
\put(00,40){\makebox(0,0)[bl]{institute}}

%--24-37

\put(00,32){\makebox(0,0)[bl]{Layer:}}
\put(00,28.5){\makebox(0,0)[bl]{depart-}}
\put(00,26){\makebox(0,0)[bl]{ment}}

%7-20

\put(00,17){\makebox(0,0)[bl]{Layer:}}
\put(00,14){\makebox(0,0)[bl]{laboratory,}}
\put(00,11.5){\makebox(0,0)[bl]{research}}
\put(00,08){\makebox(0,0)[bl]{group}}

%----

\put(82.5,12.5){\makebox(0,0)[bl]{Laboratory on}}
\put(82.5,9.5){\makebox(0,0)[bl]{data mining}}

\put(81,08){\line(1,0){25}} \put(81,17){\line(1,0){25}}
\put(81,08){\line(0,1){09}} \put(106,08){\line(0,1){09}}
\put(81.5,08){\line(0,1){09}} \put(105.5,08){\line(0,1){09}}

\end{picture}
\end{center}

 In recent years, interests to design of 'optimal'
 organizational hierarchy
 have been increased
 (e.g., \cite{gubko06}, \cite{mishin07}, \cite{voronin02}).
 Here, approaches can be based on spanning structures
 (see previous section) or direct design.
 Two basic kinds of direct design problems may be examined
 (e.g., \cite{lev89}, \cite{lev09}):

 ~{\it 1.} to build a hierarchy that has the best value(s) of its
 property(ies);

 ~{\it 2.} to build a hierarchy that is the most 'close'
 to an ideal one (or a set of ideal structures).

 Let ~\(\{ H_{i}=(A_{i},E_{i}), i=\overline{1,n} \}\)
 be a set of initial structures (hierarchies).
 Let \( H^{g}(A_{g},E_{g}) \)
 ~(\(  H^{g}  \in  \{ G_{i} \}\))
  be a goal (i.e., ideal) structure,
 let ~ \( \rho ( H_{i_{1}},H_{i_{2}} ) \)
 be a proximity (or 'distance') between two structures (hierarchies)
 \(H_{i_{1}}\) and \(H_{i_{2}}\)
 (\(  \forall H_{i_{1}},G_{i_{2}}  \in  \{ H_{i} \}\)).
 In addition, it is reasonable to consider
 properties of a structure
 (e.g., vertex degree, vertex connectivity)
 \(H_{i_{1}} \in \{ H_{i} \}\):
 ~\( \overline{\psi(H_{i_{1}}) } =
 ( \psi_{1}(H_{i_{1}}),...,\psi_{\mu}(H_{i_{1}}),..., \psi_{m} (H_{i_{1}}) )\) .
 Here hierarchy \(H^{*}  \in \{ H_{i} \}\) is searched for.
 Thus,
 in the 1st kind of the problem the objective function is
 (a case of maximization for each property):
  \[\max_{ H^{*} \in   \{ H_{i} \} } ~ \psi_{\mu } ( H^{*} ) ~~~ \forall
  ~\mu=\overline{1,m}.\]
 This problem is close to the well-known class of ``graph augmentation''  problems,
 i.e., modification of a graph to get some required properties of the modified graph
  (e.g., \cite{esw76},\cite{khu97}).

 In the second kind of the problem the objective function is:
 ~\(  \min_{H^{*} \in H_{i}} ~\rho ( H^{*},H^{g} )\).
 In general, multicriteria problem formulations may be examined as
 well.
 Clearly, it is possible to consider integrated problem formulations,
 for example:
 \[  \min_{H^{*} \in H_{i}} ~\rho ( H^{*},H^{g} ) ~~~~~~
  s.t.~~~  \psi_{\mu } ( H^{*} ) \geq d^{\mu}   ~~ \forall ~\mu=\overline{1,m},\]
 where  ~\( ( d^{1},...,d^{\mu},...,d^{m} ) \) is a vector
 constraint for properties.
 The optimization models above usually correspond to complicated integer
  (or mixed integer) optimization problems and
  here various solving approaches are used
  (e.g., enumerative methods, heuristics, AI techniques).

 In addition, it is necessary to point out
 that an attempt to build  generalized approach for
 'optimal' organization hierarchy
 is presented in
 (\cite{gubko06}, \cite{mishin07}, \cite{voronin02}).
 This approach is based on the usage of
 a  functional \(P\)
  of hierarchical structure \(G \in \Omega\)
 (\( \Omega \) is the set of hierarchical structures as direct
 acyclic 'layered' graphs;
 \( P:~ \Omega \longrightarrow ~[0,+ \infty ) \))
  and
%  general
  problem is:
 ~~\( \arg \min_{G \in \Omega} P(G) \).
 Mainly, structure \(G\) is considered as tree and
 functional \(P\) is considered as convex
 (\cite{gubko06}, \cite{mishin07}, \cite{voronin02}).
 The set of applications involves
 technological process,
 supply chain network, etc.

\subsection{Multi-layer Structures}

\subsubsection{Multi-layer Approach}

 Multi-layer (or multi-level) approach is a basic one for
 representation of complex systems
 (e.g., \cite{lev11ADES},\cite{magnani11},
 \cite{mes70},\cite{tanen02},\cite{tanen06}).
 In fact, this approach is a basic methodological method for
 decreasing the system complexity
 (i.e., multi-layer partitioning a systems and corresponding
 partitioning the system problems set).
 Here, the following main methodological steps can be pointed out.
 First, levels for main properties of complex systems have to be
 considered, for example:
 (i)  stability,
 (ii) controlability,
 (iii) adaptability, and
 (iv) self-organization \cite{mes70}.
 Second, seven-layer structure for computer systems was
 suggested as follows (e.g., \cite{tanen06}):
 (1) layer of hardware,
 (2) layer of microprogramming,
 (3) layer of operation system,
 (4) assembler-based layer,
 (5) layer of algorithmic languages,
 (6) layer of applied support systems (e.g., DBMS, DSS), and
 (7) layer of applications.
 Third, seven basic layers
 (OSI model)
 for data transmission in communication networks
 were suggested (e.g., \cite{tanen02},\cite{zimmer80}):
%\cite{sisco},
 (1) physical layer (media, signals and binary transmission),
 (2) data link layer (physical addressing),
 (3) network layer (path determination and logical addressing),
 (4) transport layer (end-to-end connections, reliability and flow control),
 (5) session layer (interhost communication, managing sessions between applications),
 (6) presentation layer (data presentation, encryption and decryption, convert machine dependent data
 to machine independent data), and
 (7) application layer (network process to application).

 Evidently, the multi-layer approach can be applied in many domains.
 For example, four-layer structure was suggested for
 representation of combinatorial ``optimization problems domain''
 in \cite{lev11ADES}:
 (i) layer of basic combinatorial optimization problems,
 (ii) layer of multicriteria combinatorial problems,
 (iii) layer of typical composite problem frameworks, and
 (iv) layer of typical applications.

 In general, a multi-layer system hierarchy consists of the following:
 (a) hierarchical layers;
 (b) set of elements for each hierarchical layer,
  description of the elements (i.e., attributes);
 (c) interconnections (some relations) over the set
 of elements for each hierarchical layer;
 (d) connections between  elements of neighbor hierarchical
 layers.
 Thus,
 the following design framework for a multi-layer system hierarchy
 can be considered:

 {\it 1.} Generation of multi-layer structure
 (i.e., the layer), for example:
  (i) dividing the initial elements/nodes into parts corresponding to
 layers (levels),
 (ii) description of layer elements (nodes, arcs) (e.g., traffic),
 (iii) building a structure for each part,
 e.g., path, multiple paths, tree, ring, complete graph
 (i.e., clique) or their combinations.

 {\it 2.} Definition (searching for)
 connections between elements of neighbor hierarchical layers.

\subsubsection{Typical Hierarchical Layers in Communication
Network}

 First, it is reasonable to
 describe  a typical hierarchy of communication networks
 (e.g.,  \cite{current86},
  \cite{gavish92}, \cite{kuz05}, \cite{mur99}).
 A traditional network hierarchy consists of the following basic layers:
 ~(a) international (multi-country, continent) network  GAN;
 ~(b) metropolitan network MN;
 ~(c) wide area network WAN; and
 ~(d) local area network LAN.
 ~{\it IBM Red Book} contains an interesting dimensional classes of
 networks by node numbers as follows:
 ~(i) large size communication network ( \(> 500\) nodes );
 ~(ii) medium size communication network (\(< 500\) nodes); and
 ~(iii) small size communication network (\(< 80\) nodes) \cite{mur99}.
 From the ``engineering'' viewpoint,
 the following hierarchical layers can be considered:
% involve the following:

 {\bf 1.} Backbone network.

 {\bf 2.} Global network as a set of
 interconnected network segments including the following:
 ~(a) additional centers,
 ~(b) cross-connections, and ~(c) bridges.

 {\bf 3.} Access network / network segment (cluster):
 (e.g., bi-connected topology, about 20 nodes).

 {\bf 4.} Distributed network: a simple hard topology
 (e.g., bus, star, tree, ring).

 As a result, a class of small-dimensional networks is added to the
 above-mentioned classification: \(20 ... 25\) nodes.
 Here a simplified network hierarchy is examined as follows (Fig 16):

 {\it 1. TOP LAYER:}  ~nodes and links for connection of clusters:
 {\it 1.1.} basic node clusters (network segments),
 {\it 1.2.} communication centers,
 {\it 1.3.} cross-connection links (center's connection), and 1.4. links
 for connection of neighbor clusters.

 {\it 2. MEDIUM LAYER:} ~basic
 clusters (network segment/ access networks, bi-connected
 topological modules).

 {\it 3. BOTTOM LAYER:} ~distribution networks (e.g., tree,
 ring, bus).

\begin{center}
\begin{picture}(92,101)

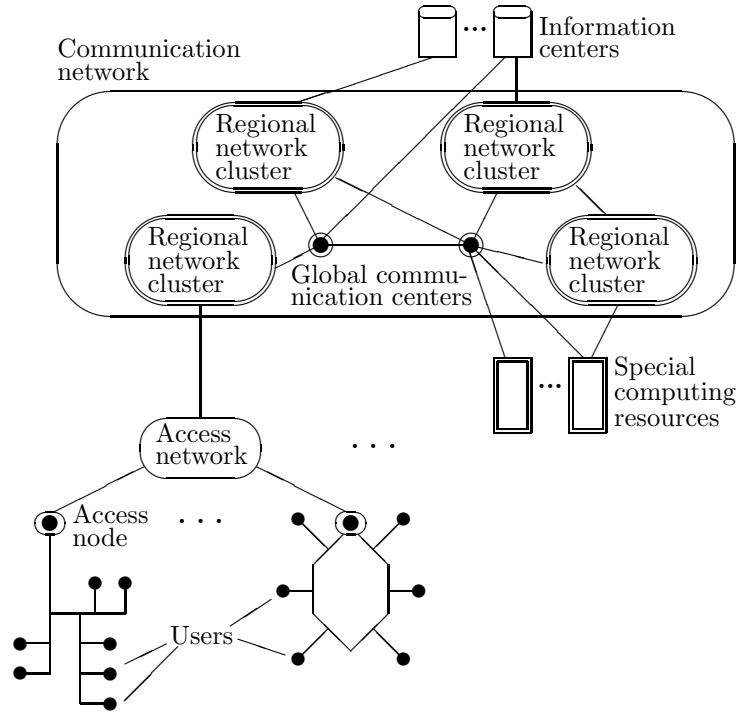
\put(00,00){\makebox(0,0)[bl]{Fig. 16. Illustration for
multi-layer communication network}}

%---- Top- layer

\put(05,92){\makebox(0,0)[bl]{Communication}}
\put(05,89){\makebox(0,0)[bl]{network}}

\put(50,72.5){\oval(90,30)}

%-- Information centers

\put(64.5,92){\line(-1,-1){24.5}}

\put(66,92){\line(0,-1){06}}

\put(55,92){\line(-3,-1){18}}

\put(55.5,98){\oval(5,1.5)}

\put(53,92){\line(1,0){05}}  \put(53,92){\line(0,1){06}}
\put(58,92){\line(0,1){06}}

\put(58.8,95.5){\makebox(0,0)[bl]{{\bf ...}}}

\put(65.5,98){\oval(5,1.5)}

\put(63,92){\line(1,0){05}} \put(63,92){\line(0,1){06}}
\put(68,92){\line(0,1){06}}

\put(69,95){\makebox(0,0)[bl]{Information }}
\put(69,92){\makebox(0,0)[bl]{centers}}

%--

\put(76,52){\line(1,2){03.5}} \put(75,52){\line(-1,1){15}}

\put(64.5,52){\line(-1,3){05}}

\put(63,42){\line(1,0){05}} \put(63,52){\line(1,0){5}}
\put(63,42){\line(0,1){10}} \put(68,42){\line(0,1){10}}

\put(63.5,42.5){\line(1,0){04}} \put(63.5,51.5){\line(1,0){4}}
\put(63.5,42.5){\line(0,1){09}} \put(67.5,42.5){\line(0,1){09}}

\put(68.8,48){\makebox(0,0)[bl]{{\bf ...}}}

\put(73,42){\line(1,0){05}} \put(73,52){\line(1,0){5}}
\put(73,42){\line(0,1){10}} \put(78,42){\line(0,1){10}}

\put(73.5,42.5){\line(1,0){04}} \put(73.5,51.5){\line(1,0){4}}
\put(73.5,42.5){\line(0,1){09}} \put(77.5,42.5){\line(0,1){09}}

\put(79,49){\makebox(0,0)[bl]{Special}}
\put(79,45.7){\makebox(0,0)[bl]{computing}}
\put(79,43){\makebox(0,0)[bl]{resources}}

%--

\put(26,81.4){\makebox(0,0)[bl]{Regional}}
\put(26,78.8){\makebox(0,0)[bl]{network}}
\put(26,75.8){\makebox(0,0)[bl]{cluster}}

\put(33,80){\oval(20,12)} \put(33,80){\oval(19,11)}

%--

\put(59,81.4){\makebox(0,0)[bl]{Regional}}
\put(59,78.8){\makebox(0,0)[bl]{network}}
\put(59,75.8){\makebox(0,0)[bl]{cluster}}

\put(66,80){\oval(20,12)} \put(66,80){\oval(19,11)}

%%%%%%%%%%%%%%%%

\put(78,71){\line(-1,1){04}}

\put(73,66.4){\makebox(0,0)[bl]{Regional}}
\put(73,63.8){\makebox(0,0)[bl]{network}}
\put(73,60.8){\makebox(0,0)[bl]{cluster}}

\put(80,65){\oval(20,12)} \put(80,65){\oval(19,11)}

%%%%%%%%%%%%%

\put(60,67){\circle{3.2}} \put(60,67){\circle*{2}}
\put(60,67){\line(-1,0){20}} \put(60,67){\line(4,-1){9.5}}

\put(60,67){\line(1,2){3.5}} \put(60,67){\line(-2,1){18}}

%%%%%%%%%%%%%

\put(40,67){\circle{3.2}} \put(40,67){\circle*{2}}
\put(40,67){\line(-2,-1){6}} \put(40,67){\line(-1,2){3.5}}

\put(36,62){\makebox(0,0)[bl]{Global commu-}}
\put(36,59){\makebox(0,0)[bl]{nication centers}}

%--

\put(17,66.4){\makebox(0,0)[bl]{Regional}}
\put(17,63.8){\makebox(0,0)[bl]{network}}
\put(17,60.8){\makebox(0,0)[bl]{cluster}}

\put(24,65){\oval(20,12)} \put(24,65){\oval(19,11)}

%%%%%%%%%%%%%%%%%%%%%%%%%%%%%%%%%%%%%%%%%

\put(44,40){\makebox(0,0)[bl]{{\bf .~.~.}}}

\put(18,40.8){\makebox(0,0)[bl]{Access}}
\put(18,37.8){\makebox(0,0)[bl]{network}}

\put(24,44){\line(0,1){15}}

\put(24,40){\oval(16,8)}

\put(21,30){\makebox(0,0)[bl]{{\bf .~.~.}}}

%%%%%%%%%%%%%%%%%%%%%%%%%%%%%%%

\put(20,14){\makebox(0,0)[bl]{Users}}

\put(19.5,14){\line(-2,-1){5.5}}

\put(21,13.5){\line(-1,-1){7.1}}

\put(28,16.3){\line(3,2){5.4}} \put(29,14.5){\line(3,-1){6.3}}

%%%%%%%%%%%%%%%%%%%%

\put(43.5,31.5){\line(-2,1){12}}

\put(44,30){\oval(4,3)} \put(44,30){\circle*{2}}

%--

\put(39,24.5){\line(1,1){4}} \put(49,24.5){\line(-1,1){4}}

\put(37,30.5){\circle*{1.7}} \put(51,30.5){\circle*{1.7}}
\put(41,26.5){\line(-1,1){04}} \put(47,26.5){\line(1,1){04}}

%--

\put(39,18){\line(0,1){6.5}} \put(49,18){\line(0,1){6.5}}

\put(35,21){\circle*{1.7}} \put(53,21){\circle*{1.7}}
\put(39,21){\line(-1,0){04}} \put(49,21){\line(1,0){04}}

%--

\put(39,18){\line(1,-1){5}} \put(49,18){\line(-1,-1){5}}

\put(37,12){\circle*{1.7}} \put(51,12){\circle*{1.7}}
\put(41,16){\line(-1,-1){04}} \put(47,16){\line(1,-1){04}}

%%%%%%%%%%%%%%%%%%%%%

\put(07,30){\makebox(0,0)[bl]{Access}}
\put(07,27){\makebox(0,0)[bl]{node}}

\put(04.5,31.5){\line(2,1){12}}

\put(04,30){\oval(4,3)} \put(04,30){\circle*{2}}

%--

\put(00,10){\circle*{1.7}} \put(00,14){\circle*{1.7}}
\put(00,10){\line(1,0){04}} \put(00,14){\line(1,0){04}}
\put(04,10){\line(0,1){18.5}}

%--

\put(12,10){\circle*{1.7}}  \put(12,14){\circle*{1.7}}
\put(12,06){\circle*{1.7}}

\put(12,06){\line(-1,0){04}} \put(12,10){\line(-1,0){04}}
\put(12,14){\line(-1,0){04}}

\put(08,06){\line(0,1){12}}

%--

\put(04,18){\line(1,0){10}}

\put(10,22){\circle*{1.7}}  \put(14,22){\circle*{1.7}}
\put(10,18){\line(0,1){04}} \put(14,18){\line(0,1){04}}

\end{picture}
\end{center}

\subsubsection{Layered  k-connected Network}

 A special version of \(2\)-connected network was suggested
 in \cite{bel97}.
 The generalized \(k\)-connected network of this kind
 and its design scheme were briefly  described in
  \cite{lev11ADES}.
 Let ~\(G=(A,E)\)~ be an initial graph (network)
 where \(A\) is a vertex (node) set,
 \(E\) is an edge set and ~\( |A| \geq k\times(k+1) + k \).

 Our kind of \(k\)-connected networks is:
 ~(a) \(k\) ``centers'' where each ``center'' is
 ~(\(k+1\))-vertex clique (the ``centers'' have not intersections),
  ~(b) node set \(M \subseteq A \) ( \(|M| \geq k\) ) includes all nodes which do
 not belong to the ``centers'';
  there is a connection (i.e., edge)
  between
  \(\forall \alpha \in M\) and
 each ``center'', i.e.,
 \(k\) edges  (one edge for a connection to a ``center'').
 As a result, the following three-layer network is obtained (Fig. 17):
 (i) layer of end nodes,
 (ii) layer of special ``central'' communication nodes:
  selected nodes and/or special additional nodes
 (the nodes are useful for location of key communication
 equipment); and
 (iii) communication ``center'' (the  ``center'' can correspond
  to ``communication providers'' / ``communication operators'' or their branches).

 Fig. 18 depicts an example of  \(4\)-connected structure.
 Proof of \(k\)-connectivity for the defined kind of networks is based on
 four  special cases.
 Consider two nodes \(a,b \in A\).
 The cases are as follows:

 {\it Case 1:} \(a\) and \(b\) belong to the same ``center''
 (i.e., ~(\(k+1\))-clique).
 Nodes \(a\) and \(b\) have connections as follows:~
 (i) direct connection \((a,b)\),
 (ii) (\(k-1\)) two-edge connection by other nodes of this
 ``center''.
 Thus, nodes \(a\) and \(b\) have \(k\) connections (without
 intersection).

 {\it Case 2:} \(a\) belongs to  ``center'' \(A\)
 and \(b\) belongs to  ``center'' \(B\)
 (i.e., another one).
 Nodes \(a\) and \(b\) have connections as follows:~
  two-edge connection by another node  \(\alpha \in M\)
 (i.e., ~(\(k\))-connections of this kind).
 Thus, nodes \(a\) and \(b\) have \(k\) connections (without
 intersection).

 {\it Case 3:} \(a\)  belongs to a ``center''
  and \(b \in B\).
  Nodes \(a\) and \(b\) have connections as follows:~
 (i) direct connection \((a,b)\),
 (ii) (\(k-1\)) three-edge connection:
  by a node of each other ``center'' and
  by another node from \(M\) ((\(k-1\)) different ways).
 Thus, nodes \(a\) and \(b\) have \(k\) connections (without
 intersection).

 {\it Case 4:} \(a\),\(b\) \( \in M\).
 Nodes \(a\) and \(b\) have \(k\) different connections
 as follows:~
 node \(a\), a ``center'', node \(b\).

\begin{center}
\begin{picture}(60,34)
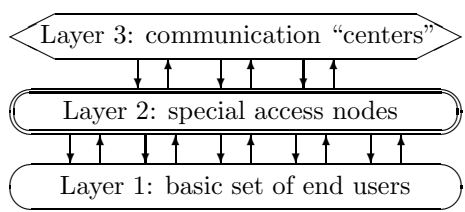
\put(0.5,0){\makebox(0,0)[bl]{Fig. 17. Three-layer model of
network}}

%-----------------------

\put(00,29){\line(2,1){06}} \put(00,29){\line(2,-1){06}}

\put(60,29){\line(-2,1){06}} \put(60,29){\line(-2,-1){06}}

\put(06,26){\line(1,0){48}} \put(06,32){\line(1,0){48}}

\put(4,27.5){\makebox(0,0)[bl]{Layer \(3\): communication
``centers''}}

%\put(30,41){\oval(56,10)}

\put(17,26){\vector(0,-1){4}} \put(21,22){\vector(0,1){4}}
\put(28,26){\vector(0,-1){4}} \put(32,22){\vector(0,1){4}}
\put(39,26){\vector(0,-1){4}} \put(43,22){\vector(0,1){4}}

%-----------------------

\put(07,17.3){\makebox(0,0)[bl]{Layer \(2\): special access
nodes}}

\put(30,19){\oval(60,6)} \put(30,19){\oval(59,05)}

\put(8,16){\vector(0,-1){4}} \put(12,12){\vector(0,1){4}}
\put(18,16){\vector(0,-1){4}} \put(22,12){\vector(0,1){4}}
\put(28,16){\vector(0,-1){4}} \put(32,12){\vector(0,1){4}}
\put(38,16){\vector(0,-1){4}} \put(42,12){\vector(0,1){4}}
\put(48,16){\vector(0,-1){4}} \put(52,12){\vector(0,1){4}}

%-----------------------End Users

\put(06.6,07.4){\makebox(0,0)[bl]{Layer \(1\):  basic set of end
users}}

\put(30,09){\oval(60,06)}

%------------------------------------

\end{picture}
\end{center}

\begin{center}
\begin{picture}(110,51)
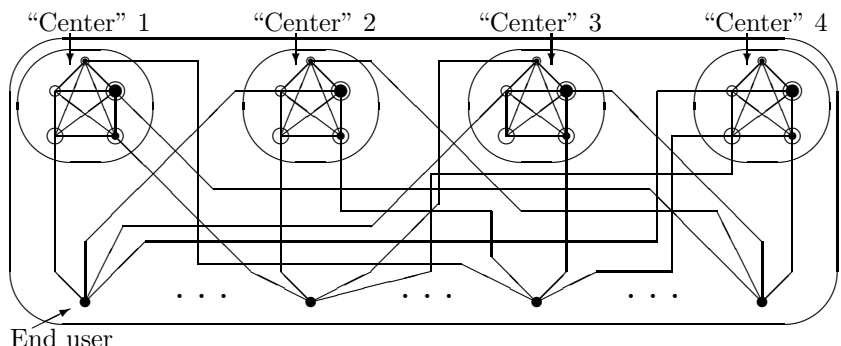
\put(18,0){\makebox(0,0)[bl]{Fig. 18. Illustration for 4-connected
network}}

\put(00,4){\makebox(0,0)[bl]{End user}}
\put(03,6.5){\vector(2,1){5}}

\put(02,46){\makebox(0,0)[bl]{``Center'' 1}}
\put(32,46){\makebox(0,0)[bl]{``Center'' 2}}
\put(62,46){\makebox(0,0)[bl]{``Center'' 3}}
\put(92,46){\makebox(0,0)[bl]{``Center'' 4}}

\put(08,45.6){\vector(0,-1){4}} \put(38,45.6){\vector(0,-1){4}}
\put(72,45.6){\vector(0,-1){4}} \put(98,45.6){\vector(0,-1){4}}
%------------------------------------------------------
\put(55,26){\oval(110,38)}

%------------------------------------------------------1

\put(10,36){\oval(18,15)}

\put(14,38){\circle*{1.8}} \put(14,38){\circle{2.6}}

\put(14,32){\circle*{1.2}} \put(14,32){\circle{2}}

\put(06,38){\circle{1.5}}

\put(06,32){\circle{2}}

\put(10,42){\circle*{0.6}} \put(10,42){\circle{1.2}}

%---

\put(06,32){\line(0,1){6}} \put(06,32){\line(1,0){8}}
\put(06,32){\line(4,3){8}} \put(14,32){\line(-4,3){8}}
\put(14,32){\line(0,1){6}} \put(06,38){\line(1,0){8}}
\put(06,38){\line(1,1){4}} \put(14,38){\line(-1,1){4}}
\put(14,32){\line(-2,5){4}} \put(06,32){\line(2,5){4}}

%------------------------------------------------------2

\put(40,36){\oval(18,15)}

\put(44,38){\circle*{1.8}} \put(44,38){\circle{2.6}}

\put(44,32){\circle*{1.2}} \put(44,32){\circle{2}}

\put(36,38){\circle{1.5}}

\put(36,32){\circle{2}}

\put(40,42){\circle*{0.6}} \put(40,42){\circle{1.2}}

%---

\put(36,32){\line(0,1){6}} \put(36,32){\line(1,0){8}}
\put(36,32){\line(4,3){8}} \put(44,32){\line(-4,3){8}}
\put(44,32){\line(0,1){6}} \put(36,38){\line(1,0){8}}
\put(36,38){\line(1,1){4}} \put(44,38){\line(-1,1){4}}
\put(44,32){\line(-2,5){4}} \put(36,32){\line(2,5){4}}

%------------------------------------------------------3

\put(70,36){\oval(18,15)}

\put(74,38){\circle*{1.8}} \put(74,38){\circle{2.6}}

\put(74,32){\circle*{1.2}} \put(74,32){\circle{2}}

\put(66,38){\circle{1.5}}

\put(66,32){\circle{2}}

\put(70,42){\circle*{0.6}} \put(70,42){\circle{1.2}}

%---

\put(66,32){\line(0,1){6}} \put(66,32){\line(1,0){8}}
\put(66,32){\line(4,3){8}} \put(74,32){\line(-4,3){8}}
\put(74,32){\line(0,1){6}} \put(66,38){\line(1,0){8}}
\put(66,38){\line(1,1){4}} \put(74,38){\line(-1,1){4}}
\put(74,32){\line(-2,5){4}} \put(66,32){\line(2,5){4}}

%------------------------------------------------------4

\put(100,36){\oval(18,15)}

\put(104,38){\circle*{1.8}} \put(104,38){\circle{2.6}}

\put(104,32){\circle*{1.2}} \put(104,32){\circle{2}}

\put(96,38){\circle{1.5}}

\put(96,32){\circle{2}}

\put(100,42){\circle*{0.6}} \put(100,42){\circle{1.2}}

%---

\put(96,32){\line(0,1){6}} \put(96,32){\line(1,0){8}}
\put(96,32){\line(4,3){8}} \put(104,32){\line(-4,3){8}}
\put(104,32){\line(0,1){6}} \put(96,38){\line(1,0){8}}
\put(96,38){\line(1,1){4}} \put(104,38){\line(-1,1){4}}
\put(104,32){\line(-2,5){4}} \put(96,32){\line(2,5){4}}

%=====================================

%----------------------1

\put(10,10){\circle*{1.3}}

\put(10,10){\line(-1,1){4}} \put(06,14){\line(0,1){18}}

\put(10,10){\line(0,1){08}} \put(10,18){\line(1,1){20}}
\put(30,38){\line(1,0){06}}

\put(10,10){\line(1,2){5}} \put(15,20){\line(1,0){33}}
\put(48,20){\line(1,1){18}}

\put(10,10){\line(1,1){8}} \put(18,18){\line(1,0){68}}
\put(86,18){\line(0,1){20}}

\put(86,38){\line(1,0){10}}

%----------------------2

\put(40,10){\circle*{1.3}}

\put(40,10){\line(-1,1){4}} \put(36,14){\line(0,1){18}}

\put(40,10){\line(-2,1){8}} \put(32,14){\line(-1,1){18}}

\put(40,10){\line(2,1){8}} \put(48,14){\line(1,1){9}}

\put(57,23){\line(0,1){19}} \put(57,42){\line(1,0){13}}

\put(40,10){\line(4,1){16}} \put(56,14){\line(0,1){13}}
\put(56,27){\line(1,0){40}} \put(96,27){\line(0,1){5}}

%----------------------3

\put(70,10){\circle*{1.3}}

\put(70,10){\line(1,1){4}} \put(74,14){\line(0,1){18}}

\put(70,10){\line(2,1){8}} \put(78,14){\line(1,0){10}}
\put(88,14){\line(0,1){18}} \put(88,32){\line(1,0){08}}

\put(70,10){\line(-1,1){6}} \put(64,16){\line(0,1){6}}
\put(64,22){\line(-1,0){20}} \put(44,22){\line(0,1){10}}

\put(70,10){\line(-2,1){10}} \put(60,15){\line(-1,0){35}}
\put(25,15){\line(0,1){27}} \put(25,42){\line(-1,0){15}}

%----------------------4

\put(100,10){\circle*{1.3}}

\put(100,10){\line(1,1){4}} \put(104,14){\line(0,1){18}}

\put(100,10){\line(0,1){08}} \put(100,18){\line(-1,1){20}}
\put(80,38){\line(-1,0){06}}

\put(100,10){\line(-1,2){6}} \put(94,22){\line(-1,0){26}}
\put(68,22){\line(-1,1){20}} \put(48,42){\line(-1,0){8}}

\put(100,10){\line(-1,1){15}} \put(85,25){\line(-1,0){58}}
\put(27,25){\line(-1,1){13}}

%===========================================

\put(22,10.5){\makebox(0,0)[bl]{{\bf . . .}}}

%---

\put(52,10.5){\makebox(0,0)[bl]{{\bf . . .}}}

%---

\put(82,10.5){\makebox(0,0)[bl]{{\bf . . .}}}

\end{picture}
\end{center}

 Further, a solving scheme ('Bottom-Up')
 to build \(k\)-connected structure is:

%~~

 {\it Stage 1.} Selection of ~\(k \times (k+1)\)~ vertices
 for ~\(k\)  ``centers''
 (here multicriteria ranking problem can be used).

 {\it Stage 2.} Clustering of the selected vertices to get
 \(k\) clusters as  ``centers''
 (each cluster consists of ~(\(k+1\)) vertices,
 the clusters have not intersections).

 {\it Stage 3.} Building of connections for end users:
 connection of each vertex
 (that does not belong to a ``center'')
  with each ``centers''
  (i.e., with the only one node of each ``centers'').
 Here multiple choice problem or its modifications can be used.

%~~

 Now an illustrative example of illustrate the examined
 type of bi-connected communication network is presented.
 Initial nodes are depicted in Fig. 19.

\begin{center}
\begin{picture}(110,47)
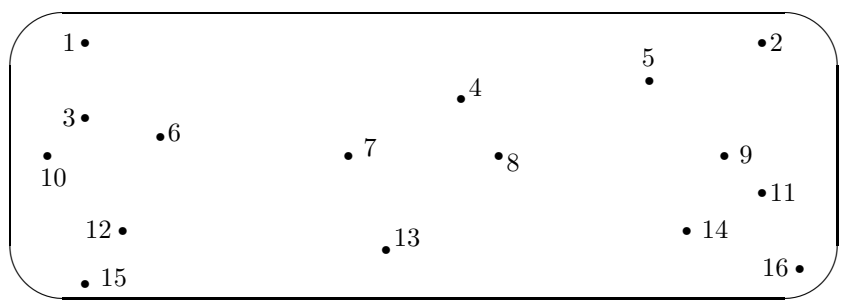
\put(30,0){\makebox(0,0)[bl]{Fig. 19. Initial network nodes}}

%------------------------------------------------------
\put(55,25){\oval(110,38)}

%----------------------------------Center 1

\put(05,25){\circle*{1}} \put(04,21){\makebox(0,0)[bl]{\(10\)}}

\put(10,30){\circle*{1}} \put(07,29){\makebox(0,0)[bl]{\(3\)}}

\put(15,15){\circle*{1}} \put(10,14){\makebox(0,0)[bl]{\(12\)}}

%----------------------------------Center 2

\put(95,25){\circle*{1}} \put(97,24){\makebox(0,0)[bl]{\(9\)}}

\put(85,35){\circle*{1}} \put(84,37){\makebox(0,0)[bl]{\(5\)}}

\put(90,15){\circle*{1}} \put(92,14){\makebox(0,0)[bl]{\(14\)}}

%------------------------------------------------------1

\put(50,12.5){\circle*{1}} \put(51,13){\makebox(0,0)[bl]{\(13\)}}

%------------------------------------------------------2

\put(10,08){\circle*{1}} \put(12,07.6){\makebox(0,0)[bl]{\(15\)}}

%------------------------------------------------------3

\put(45,25){\circle*{1}} \put(47,25){\makebox(0,0)[bl]{\(7\)}}

%------------------------------------------------------4

\put(10,40){\circle*{1}} \put(7,39){\makebox(0,0)[bl]{\(1\)}}

%------------------------------------------------------5

\put(100,40){\circle*{1}} \put(101,39){\makebox(0,0)[bl]{\(2\)}}

%------------------------------------------------------6

\put(100,20){\circle*{1}} \put(101,19){\makebox(0,0)[bl]{\(11\)}}

%------------------------------------------------------7

\put(20,27.5){\circle*{1}} \put(21,27){\makebox(0,0)[bl]{\(6\)}}

%------------------------------------------------------8

\put(65,25){\circle*{1}} \put(66,23){\makebox(0,0)[bl]{\(8\)}}

%------------------------------------------------------9

\put(60,32.5){\circle*{1}} \put(61,033){\makebox(0,0)[bl]{\(4\)}}

%------------------------------------------------------10

\put(105,10){\circle*{1}} \put(100,09){\makebox(0,0)[bl]{\(16\)}}

\end{picture}
\end{center}

 Here the following nodes are selected for ``centers'':
 ~\(\{3,5,9,10,12,14\}\).
 Other nodes correspond to users.

 The corresponding  version ({\it version 1}) of the resultant bi-connected networks is presented
 in Fig. 20:

 {\it Version 1}. ``Centers'' are  locally-allocated
 (i.e., regional ``centers''):
``center'' \(1\): nodes \(3\),  \(10\),  and \(12\);
  ``center'' \(2\): nodes \(5\),  \(9\),  and \(14\).

\begin{center}
\begin{picture}(110,51)
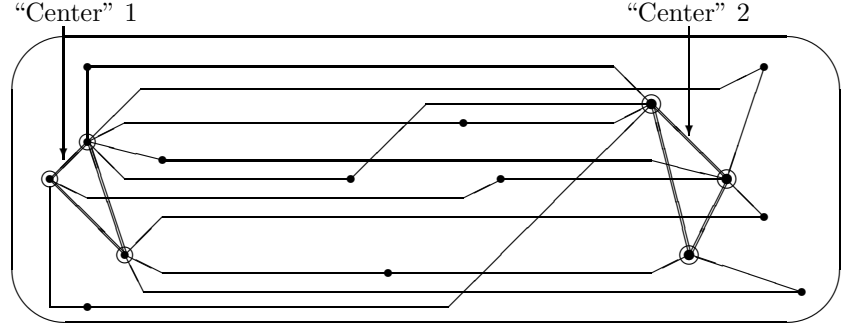
\put(11,0){\makebox(0,0)[bl]{Fig. 20. Bi-connected network:
locally-allocated  ``centers }}

%------------------------------------------------------
\put(55,25){\oval(110,38)}

%-------------------------------------------

\put(00,46){\makebox(0,0)[bl]{``Center'' 1}}
\put(81.4,46){\makebox(0,0)[bl]{``Center'' 2}}

\put(06.8,45.3){\vector(0,-1){17.6}}
\put(90,45.3){\vector(0,-1){14.5}}

%----------------------------------Center 1

\put(05,25){\circle*{1.2}} \put(05,25){\circle{2}}
\put(10,30){\circle*{1.2}} \put(10,30){\circle{2}}
\put(15,15){\circle*{1.2}} \put(15,15){\circle{2}}

\put(05,25){\line(1,1){5}} \put(05,25){\line(1,-1){10}}
\put(15,15){\line(-1,3){5}}

\put(05.3,25){\line(1,1){5}} \put(05.3,25){\line(1,-1){10}}
\put(15.3,15){\line(-1,3){5}}

%----------------------------------Center 2

\put(95,25){\circle*{1.5}} \put(95,25){\circle{2.4}}
\put(85,35){\circle*{1.5}} \put(85,35){\circle{2.4}}
\put(90,15){\circle*{1.5}} \put(90,15){\circle{2.4}}

\put(95,25){\line(-1,1){10}} \put(95,25){\line(-1,-2){5}}
\put(90,15){\line(-1,4){5}}

\put(95.3,25){\line(-1,1){10}} \put(95.3,25){\line(-1,-2){5}}
\put(90.3,15){\line(-1,4){5}}

%------------------------------------------------------1

\put(50,12.5){\circle*{1}}

\put(50,12.5){\line(1,0){35}} \put(85,12.5){\line(2,1){5}}
\put(50,12.5){\line(-1,0){30}} \put(20,12.5){\line(-2,1){5}}

%------------------------------------------------------2

\put(10,08){\circle*{1}}

\put(05,08){\line(0,1){16}}

\put(05,08){\line(1,0){53}} \put(58,08){\line(1,1){27}}

%------------------------------------------------------3

\put(45,25){\circle*{1}}

\put(45,25){\line(1,1){10}} \put(55,35){\line(1,0){30}}

\put(15,25){\line(-1,1){5}} \put(45,25){\line(-1,0){30}}

%------------------------------------------------------4

\put(10,40){\circle*{1}}

\put(10,40){\line(0,-1){10}}

\put(10,40){\line(1,0){70}} \put(80,40){\line(1,-1){5}}

%------------------------------------------------------5

\put(100,40){\circle*{1}}

\put(100,40){\line(-1,-3){5}}

\put(100,40){\line(-2,-1){6}} \put(94,37){\line(-1,0){77}}
\put(17,37){\line(-1,-1){7}}

%------------------------------------------------------6

\put(100,20){\circle*{1}}

\put(100,20){\line(-1,1){5}} \put(100,20){\line(-1,0){80}}
\put(20,20){\line(-1,-1){5}}

%------------------------------------------------------7

\put(20,27.5){\circle*{1}}

\put(20,27.5){\line(-4,1){10}}

\put(20,27.5){\line(1,0){65}} \put(85,27.5){\line(4,-1){10}}

%------------------------------------------------------8

\put(65,25){\circle*{1}}

\put(65,25){\line(1,0){30}} \put(65,25){\line(-2,-1){5}}
\put(60,22.5){\line(-1,0){50}} \put(10,22.5){\line(-2,1){5}}

%------------------------------------------------------9

\put(60,32.5){\circle*{1}}

\put(60,32.5){\line(-1,0){45}}  \put(15,32.5){\line(-2,-1){5}}

\put(60,32.5){\line(1,0){20}}  \put(80,32.5){\line(2,1){5}}

%------------------------------------------------------10

\put(105,10){\circle*{1}}

\put(105,10){\line(-3,1){15}}

\put(105,10){\line(-1,0){87.5}}

\put(17.5,10){\line(-1,2){2.5}}

\end{picture}
\end{center}

 On the other hand, it may  be reasonable
 to consider another situation with
 ``simplification'' of  connections for end users
 (e.g., minimization of connection distances).
 In this case, it is necessary to use the following algorithmic rule
 at {\it stage 2} of the above-mentioned solving scheme
 (i.e., grouping of the selected nodes to centers):

~~

  {\it Nodes of each ``center'' have to 'cover' the network
 (e.g., to be very close to all nodes of end users).}

~~

 Thus,  the second solving scheme to build \(k\)-connected structure is:

~~

 {\it Stage 1.} Clustering of of the initial set of nodes to get
 \(k\) clusters.

 {\it Stage 2.} Building of \(k\) ``centers'':
 selection  of a vertex in each obtained cluster
 as an element for each ``center''
 (a representative of the ``center''),
  connection of the vertices in each ``center''.

 {\it Stage 3.} Building of connections for end users:
 connection of each vertex
 (that does not belong to a ``center'')
  with each ``centers''
  (i.e., with the only one node of each ``centers'').

~~

 Fig. 21 depicts an example of this type of the resultant
 bi-connected network:

  {\it Version 2}. ``Centers'' are distributed over the network:
 ``center'' \(1\): nodes \(3\),  \(9\),  and \(12\);
  ``center'' \(2\): nodes \(5\),  \(10\),  and \(14\).

\begin{center}
\begin{picture}(110,51)
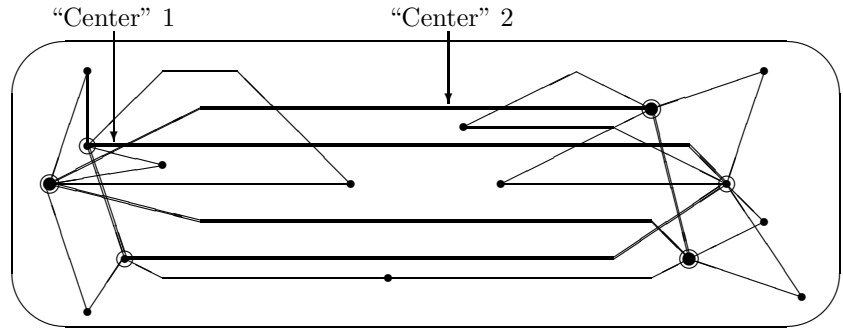
\put(04,0){\makebox(0,0)[bl]{Fig. 21. Bi-connected network:
 distributed over network  ``centers''}}

%------------------------------------------------------
\put(55,25){\oval(110,38)}

%-------------------------------------------

\put(05,46){\makebox(0,0)[bl]{``Center'' 1}}
\put(50,46){\makebox(0,0)[bl]{``Center'' 2}}

\put(13.5,45.3){\vector(0,-1){14.66}}
\put(58,45.3){\vector(0,-1){9.6}}

%----------------------------------Center 1

\put(10,30){\circle*{1.2}} \put(10,30){\circle{2}}
%\put(07,29){\makebox(0,0)[bl]{\(3\)}}

\put(10,30){\line(1,0){80}} \put(90,30){\line(1,-1){5}}
\put(10,30.3){\line(1,0){80}} \put(90,30.3){\line(1,-1){5}}

\put(10,30){\line(1,-3){5}} \put(10.3,30){\line(1,-3){5}}

\put(95,25){\circle*{1.2}} \put(95,25){\circle{2}}
%\put(97,24){\makebox(0,0)[bl]{\(9\)}}

\put(15,15){\circle*{1.2}} \put(15,15){\circle{2}}
%\put(10,14){\makebox(0,0)[bl]{\(12\)}}

\put(15,15){\line(1,0){65}} \put(80,15){\line(3,2){15}}
\put(15,15.3){\line(1,0){65}} \put(80,15.3){\line(3,2){15}}

%----------------------------------Center 2

\put(85,35){\circle*{1.6}} \put(85,35){\circle{2.5}}
%\put(84,37){\makebox(0,0)[bl]{\(5\)}}

\put(85,35){\line(1,-4){5}} \put(85.2,35){\line(1,-4){5}}

\put(05,25){\circle*{1.6}} \put(05,25){\circle{2.5}}
%\put(04,21){\makebox(0,0)[bl]{\(10\)}}

\put(05,25){\line(4,-1){20}} \put(25,20){\line(1,0){60}}
\put(85,20){\line(1,-1){5}}

\put(05,25.2){\line(4,-1){20}} \put(25,20.2){\line(1,0){60}}
\put(85,20.2){\line(1,-1){5}}

\put(05,25){\line(2,1){20}} \put(25,35){\line(1,0){60}}
\put(05,25.2){\line(2,1){20}} \put(25,35.2){\line(1,0){60}}

\put(90,15){\circle*{1.6}} \put(90,15){\circle{2.5}}
%\put(92,14){\makebox(0,0)[bl]{\(14\)}}

%------------------------------------------------------1

\put(50,12.5){\circle*{1}}
%\put(51,13){\makebox(0,0)[bl]{\(13\)}}

\put(50,12.5){\line(1,0){35}} \put(85,12.5){\line(2,1){5}}
\put(50,12.5){\line(-1,0){30}} \put(20,12.5){\line(-2,1){5}}

%------------------------------------------------------2

\put(10,08){\circle*{1}}
%\put(12,07.6){\makebox(0,0)[bl]{\(15\)}}

\put(10,08){\line(2,3){5}} \put(10,08){\line(-1,3){05.45}}

%------------------------------------------------------3

\put(45,25){\circle*{1}}
%\put(47,25){\makebox(0,0)[bl]{\(7\)}}

\put(45,25){\line(-1,0){40}}

\put(45,25){\line(-1,1){15}} \put(30,40){\line(-1,0){10}}
\put(20,40){\line(-1,-1){10}}

%------------------------------------------------------4

\put(10,40){\circle*{1}}
%\put(7,39){\makebox(0,0)[bl]{\(1\)}}

\put(10,40){\line(-1,-3){5}} \put(10,40){\line(0,-1){10}}

%------------------------------------------------------5

\put(100,40){\circle*{1}}
%\put(101,39){\makebox(0,0)[bl]{\(2\)}}

\put(100,40){\line(-3,-1){15}} \put(100,40){\line(-1,-3){5}}

%------------------------------------------------------6

\put(100,20){\circle*{1}}
%\put(101,19){\makebox(0,0)[bl]{\(11\)}}

\put(100,20){\line(-1,1){5}} \put(100,20){\line(-2,-1){10}}

%------------------------------------------------------7

\put(20,27.5){\circle*{1}}
%\put(21,27){\makebox(0,0)[bl]{\(6\)}}

\put(20,27.5){\line(-4,1){10}} \put(20,27.5){\line(-6,-1){15}}

%------------------------------------------------------8

\put(65,25){\circle*{1}}
%\put(66,23){\makebox(0,0)[bl]{\(8\)}}

\put(65,25){\line(1,0){30}} \put(65,25){\line(2,1){20}}

%------------------------------------------------------9

\put(60,32.5){\circle*{1}}
%\put(61,33){\makebox(0,0)[bl]{\(4\)}}

\put(60,32.5){\line(2,1){15}} \put(75,40){\line(2,-1){10}}

\put(60,32.5){\line(1,0){20}} \put(80,32.5){\line(2,-1){15}}

%------------------------------------------------------10

\put(105,10){\circle*{1}}
%\put(100,09){\makebox(0,0)[bl]{\(16\)}}

\put(105,10){\line(-3,1){15}} \put(105,10){\line(-2,3){10}}

\end{picture}
\end{center}

\subsubsection{Towards Hierarchical Network Design Problems}

 The basic hierarchical two-level network design problem consists in
 finding  a minimum cost two-level spanning network,
 consisting of two parts:
 (i) main path (or several paths, tree, ring)
 (ii) secondary trees
 (e.g.,
 \cite{current86}, \cite{obr08},  \cite{obr10}, \cite{pirkul91}).
 Thus, the initial network is divided into two parts:

 {\it 1.} The higher level part (primary nodes, main path):
  a path (or several paths, tree, ring)
 composed of primary arcs,
 which visits some of the nodes of the network
 (i.e., primary nodes).

 {\it 2.} The lower level part (secondary nodes, secondary trees):
 the part is composed of one or more
 trees whose arcs, termed secondary,
 are less expensive to build than the primary arcs.

 Here, each arc has a cost (\(d_{ij},~ \forall i,j \in A\), \(A\) is the set of
 nodes).
 The total cost of the selected arcs in the spanning structure
 is used as the minimized objective function.
 The problem is formulated as combinatorial optimization model
 (e.g., \cite{current86}),
 it is NP-hard \cite{bala94}.
 Various approaches have been suggested for the problem,
 for example:
 (a) exact (enumerative) methods
 (e.g., branch-and-cut approach, dynamic programming),
 (b) heuristics (e.g., Lagrangian relaxation),
 (c) evolutionary algorithms.
 Often, a preliminary minimum spanning tree is used
 to construct the solution
 (leaf nodes of a spanning tree are used as the secondary nodes).
 Further, three illustrative examples are presented:

 (i) higher part is path over
 primary nodes  ~\(<1,6,5,7,2,3,4>\)~ (Fig. 22),

 (ii) higher part is tree over
 primary nodes  ~\(\{1,2,3,4,5,6,7\}\)~ (Fig. 23),

 (iii) higher part is ring over
 primary nodes  ~\(\{1,2,3,4,5,6,7\}\)~
 (two-connected case) (Fig. 24).

 Main applications of the problem involve communication networks,
 computer networks,
 transportation networks,
 power line distributed systems.
 Evidently, more general multi-level networks are examined as well
(e.g.,
 \cite{bala94}, \cite{chorpa02}).

\begin{center}
\begin{picture}(110,47)
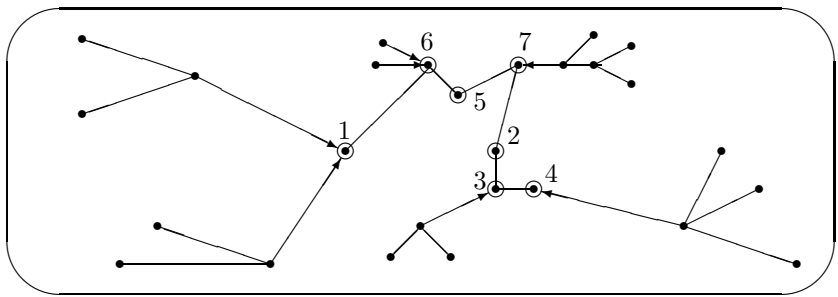
\put(09,0){\makebox(0,0)[bl]{Fig. 22. Example of two-level network
 (higher level: path)}}

%------------------------------------------------------
\put(55,25){\oval(110,38)}

%%%%%%%%%%%%%%%%%%%%%%%% Primary PATH
%----------------------------------------------------7->1

\put(45,25){\circle*{1}} \put(45,25){\circle{2}}

\put(44,26.7){\makebox(0,0)[bl]{\(1\)}}

\put(45,25){\line(1,1){11.3}}

%\put(45,25){\line(1,0){20}}

%------------------------------------------Root 15
%\put(45,25){\line(-2,-3){10}}

\put(35,10){\vector(2,3){09.4}}

\put(35,10){\circle*{1}}

%\put(36,09){\makebox(0,0)[bl]{\(15\)}}

%--

\put(35,10){\line(-1,0){20}} \put(15,10){\circle*{1}}

%\put(05,11){\makebox(0,0)[bl]{\(12\)}}

%--
\put(35,10){\line(-3,1){15}}

\put(20,15){\circle*{1}}

%\put(04,21){\makebox(0,0)[bl]{\(10\)}}

%--------------------------------------------ROOT 6
%\put(45,25){\line(-2,1){20}}

\put(25,35){\vector(2,-1){19}}

\put(25,35){\circle*{1}}

%\put(26,37){\makebox(0,0)[bl]{\(6\)}}

%--1

\put(25,35){\line(-3,1){15}}

\put(10,40){\circle*{1}}

%\put(7,39){\makebox(0,0)[bl]{\(1\)}}

%--

\put(25,35){\line(-3,-1){15}}

\put(10,30){\circle*{1}}

%\put(07,29){\makebox(0,0)[bl]{\(3\)}}

%------------------------------------------------------8->2
%----2-7
\put(65,25){\line(1,4){3}}

\put(65,25){\circle*{1}} \put(65,25){\circle{2}}

\put(66.5,26){\makebox(0,0)[bl]{\(2\)}}

\put(62,20){\makebox(0,0)[bl]{\(3\)}}

\put(71.5,21){\makebox(0,0)[bl]{\(4\)}}

\put(65,25){\line(0,-1){5}} \put(65,20){\line(1,0){5}}

%\put(65,25){\line(1,-1){5}}

%-----------------------------8-1
\put(65,20){\circle*{1}} \put(65,20){\circle{2}}

%---------------------- Root  13
%\put(65,20){\line(-2,-1){10}}

\put(55,15){\vector(2,1){09}}

\put(55,15){\circle*{1}}

%\put(56,13){\makebox(0,0)[bl]{\(13\)}}

%--
\put(55,15){\line(-1,-1){04}} \put(51,11){\circle*{1}}

%--
\put(55,15){\line(1,-1){04}} \put(59,11){\circle*{1}}

%-----------------------------8-2
\put(70,20){\circle*{1}} \put(70,20){\circle{2}}

%----------------------- Root 14
%\put(70,20){\line(4,-1){20}}

\put(90,15){\vector(-4,1){19}}

\put(90,15){\circle*{1}}

%\put(92,14){\makebox(0,0)[bl]{\(14\)}}

%--
\put(90,15){\line(1,2){5}}

\put(95,25){\circle*{1}}

%\put(97,24){\makebox(0,0)[bl]{\(9\)}}

%--

\put(90,15){\line(2,1){10}}

\put(100,20){\circle*{1}}

%\put(101,19){\makebox(0,0)[bl]{\(11\)}}

%--
\put(90,15){\line(3,-1){15}}

\put(105,10){\circle*{1}}

%\put(100,09){\makebox(0,0)[bl]{\(16\)}}

%-----------------------------------4->5
\put(60,32.5){\circle*{1}} \put(60,32.5){\circle{2}}

\put(62,30.5){\makebox(0,0)[bl]{\(5\)}}

\put(55,38.5){\makebox(0,0)[bl]{\(6\)}}

\put(60,32.5){\line(-1,1){4}} \put(60,32.5){\line(2,1){8}}

%-----------------4-1
\put(56,36.5){\circle*{1}} \put(56,36.5){\circle{2}}

%\put(56,36.5){\line(-2,1){6}}

\put(50,39.5){\vector(2,-1){5}}

\put(50,39.5){\circle*{1}}

%\put(56,36.5){\line(-1,0){7}}

\put(49,36.5){\vector(1,0){6.5}}

\put(49,36.5){\circle*{1}}

%-----------------4-2
\put(68,36.5){\circle*{1}} \put(68,36.5){\circle{2}}

\put(68,38.5){\makebox(0,0)[bl]{\(7\)}}

%\put(64,36.5){\line(1,0){6}}

\put(74,36.5){\vector(-1,0){5.4}}

%%%%%%%%%%%%%%%%%%%%%%%%%%%%%%%%%%%%%%%%%%%%%%%%%%%%%%%

\put(74,36.5){\circle*{1}}

\put(74,36.5){\line(1,1){4}} \put(78,40.5){\circle*{1}}

\put(74,36.5){\line(1,0){5}} \put(78,36.5){\circle*{1}}

\put(78,36.5){\line(2,1){5}} \put(83,39){\circle*{1}}

\put(78,36.5){\line(2,-1){5}} \put(83,34){\circle*{1}}

\end{picture}
\end{center}

\begin{center}
\begin{picture}(110,47)
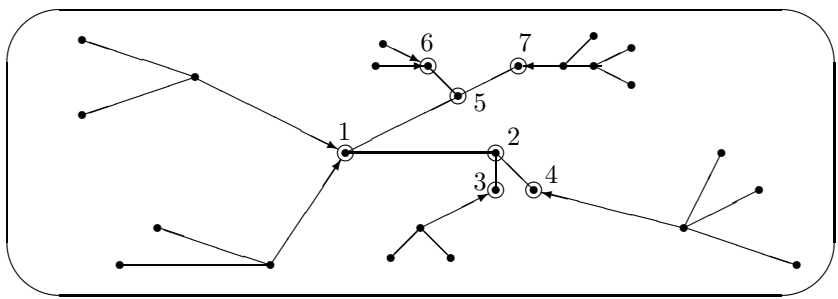
\put(11,0){\makebox(0,0)[bl]{Fig. 23. Example of two-level network
 (higher level: tree)}}

%------------------------------------------------------
\put(55,25){\oval(110,38)}

%%%%%%%%%%%%%%%%%%%%%%%% Primary TREE
%----------------------------------------------------7->1

\put(45,25){\circle*{1}} \put(45,25){\circle{2}}

\put(44,26.7){\makebox(0,0)[bl]{\(1\)}}

\put(45,25){\line(2,1){15}}

\put(45,25){\line(1,0){20}}

%------------------------------------------Root 15
%\put(45,25){\line(-2,-3){10}}

\put(35,10){\vector(2,3){09.4}}

\put(35,10){\circle*{1}}

%\put(36,09){\makebox(0,0)[bl]{\(15\)}}

%--

\put(35,10){\line(-1,0){20}} \put(15,10){\circle*{1}}

%\put(05,11){\makebox(0,0)[bl]{\(12\)}}

%--
\put(35,10){\line(-3,1){15}}

\put(20,15){\circle*{1}}

%\put(04,21){\makebox(0,0)[bl]{\(10\)}}

%--------------------------------------------ROOT 6
%\put(45,25){\line(-2,1){20}}

\put(25,35){\vector(2,-1){19}}

\put(25,35){\circle*{1}}

%\put(26,37){\makebox(0,0)[bl]{\(6\)}}

%--1

\put(25,35){\line(-3,1){15}}

\put(10,40){\circle*{1}}

%\put(7,39){\makebox(0,0)[bl]{\(1\)}}

%--

\put(25,35){\line(-3,-1){15}}

\put(10,30){\circle*{1}}

%\put(07,29){\makebox(0,0)[bl]{\(3\)}}

%------------------------------------------------------8->2

\put(65,25){\circle*{1}} \put(65,25){\circle{2}}

\put(66.5,26){\makebox(0,0)[bl]{\(2\)}}

\put(62,20){\makebox(0,0)[bl]{\(3\)}}

\put(71.5,21){\makebox(0,0)[bl]{\(4\)}}

\put(65,25){\line(0,-1){5}} \put(65,25){\line(1,-1){5}}

%-----------------------------8-1
\put(65,20){\circle*{1}} \put(65,20){\circle{2}}

%---------------------- Root  13
%\put(65,20){\line(-2,-1){10}}

\put(55,15){\vector(2,1){09}}

\put(55,15){\circle*{1}}

%\put(56,13){\makebox(0,0)[bl]{\(13\)}}

%--
\put(55,15){\line(-1,-1){04}} \put(51,11){\circle*{1}}

%--
\put(55,15){\line(1,-1){04}} \put(59,11){\circle*{1}}

%-----------------------------8-2
\put(70,20){\circle*{1}} \put(70,20){\circle{2}}

%----------------------- Root 14
%\put(70,20){\line(4,-1){20}}

\put(90,15){\vector(-4,1){19}}

\put(90,15){\circle*{1}}

%\put(92,14){\makebox(0,0)[bl]{\(14\)}}

%--
\put(90,15){\line(1,2){5}}

\put(95,25){\circle*{1}}

%\put(97,24){\makebox(0,0)[bl]{\(9\)}}

%--

\put(90,15){\line(2,1){10}}

\put(100,20){\circle*{1}}

%\put(101,19){\makebox(0,0)[bl]{\(11\)}}

%--
\put(90,15){\line(3,-1){15}}

\put(105,10){\circle*{1}}

%\put(100,09){\makebox(0,0)[bl]{\(16\)}}

%-----------------------------------4->5
\put(60,32.5){\circle*{1}} \put(60,32.5){\circle{2}}

\put(62,30.5){\makebox(0,0)[bl]{\(5\)}}

\put(55,38.5){\makebox(0,0)[bl]{\(6\)}}

\put(60,32.5){\line(-1,1){4}} \put(60,32.5){\line(2,1){8}}

%-----------------4-1
\put(56,36.5){\circle*{1}} \put(56,36.5){\circle{2}}

%\put(56,36.5){\line(-2,1){6}}

\put(50,39.5){\vector(2,-1){5}}

\put(50,39.5){\circle*{1}}

%\put(56,36.5){\line(-1,0){7}}

\put(49,36.5){\vector(1,0){6.5}}

\put(49,36.5){\circle*{1}}

%-----------------4-2
\put(68,36.5){\circle*{1}} \put(68,36.5){\circle{2}}

\put(68,38.5){\makebox(0,0)[bl]{\(7\)}}

%\put(64,36.5){\line(1,0){6}}

\put(74,36.5){\vector(-1,0){5.4}}

%%%%%%%%%%%%%%%%%%%%%%%%%%%%%%%%%%%%%%%%%%%%%%%%%%%%%%%

\put(74,36.5){\circle*{1}}

\put(74,36.5){\line(1,1){4}} \put(78,40.5){\circle*{1}}

\put(74,36.5){\line(1,0){5}} \put(78,36.5){\circle*{1}}

\put(78,36.5){\line(2,1){5}} \put(83,39){\circle*{1}}

\put(78,36.5){\line(2,-1){5}} \put(83,34){\circle*{1}}

\end{picture}
\end{center}

\begin{center}
\begin{picture}(110,47)
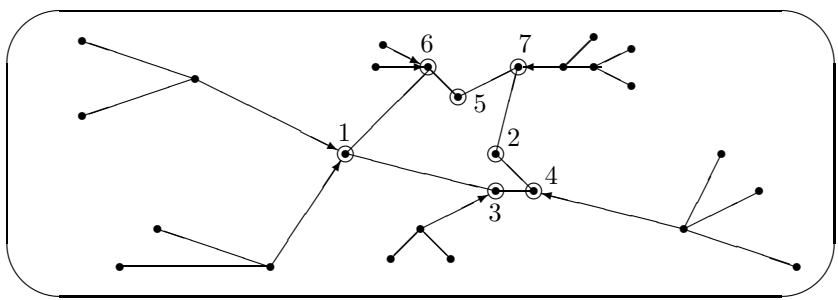
\put(10,0){\makebox(0,0)[bl]{Fig. 24. Example of two-level network
 (higher level: ring)}}

%------------------------------------------------------
\put(55,25){\oval(110,38)}

%%%%%%%%%%%%%%%%%%%%%%%% Primary RING
%----------------------------------------------------7->1

\put(45,25){\circle*{1}} \put(45,25){\circle{2}}

\put(44,26.7){\makebox(0,0)[bl]{\(1\)}}

\put(45,25){\line(1,1){11.3}}

\put(45,25){\line(4,-1){20}}

%------------------------------------------Root 15
%\put(45,25){\line(-2,-3){10}}

\put(35,10){\vector(2,3){09.4}}

\put(35,10){\circle*{1}}

%\put(36,09){\makebox(0,0)[bl]{\(15\)}}

%--

\put(35,10){\line(-1,0){20}} \put(15,10){\circle*{1}}

%\put(05,11){\makebox(0,0)[bl]{\(12\)}}

%--
\put(35,10){\line(-3,1){15}}

\put(20,15){\circle*{1}}

%\put(04,21){\makebox(0,0)[bl]{\(10\)}}

%--------------------------------------------ROOT 6
%\put(45,25){\line(-2,1){20}}

\put(25,35){\vector(2,-1){19}}

\put(25,35){\circle*{1}}

%\put(26,37){\makebox(0,0)[bl]{\(6\)}}

%--1

\put(25,35){\line(-3,1){15}}

\put(10,40){\circle*{1}}

%\put(7,39){\makebox(0,0)[bl]{\(1\)}}

%--

\put(25,35){\line(-3,-1){15}}

\put(10,30){\circle*{1}}

%\put(07,29){\makebox(0,0)[bl]{\(3\)}}

%------------------------------------------------------8->2
%----2-7
\put(65,25){\line(1,4){3}}

\put(65,25){\circle*{1}} \put(65,25){\circle{2}}

\put(66.5,26){\makebox(0,0)[bl]{\(2\)}}

\put(64,16){\makebox(0,0)[bl]{\(3\)}}

\put(71.5,21){\makebox(0,0)[bl]{\(4\)}}

%\put(65,25){\line(0,-1){5}}

\put(70,20){\line(-1,0){5}}

\put(65,25){\line(1,-1){5}}

%-----------------------------8-1
\put(65,20){\circle*{1}} \put(65,20){\circle{2}}

%---------------------- Root  13
%\put(65,20){\line(-2,-1){10}}

\put(55,15){\vector(2,1){09}}

\put(55,15){\circle*{1}}

%\put(56,13){\makebox(0,0)[bl]{\(13\)}}

%--
\put(55,15){\line(-1,-1){04}} \put(51,11){\circle*{1}}

%--
\put(55,15){\line(1,-1){04}} \put(59,11){\circle*{1}}

%-----------------------------8-2
\put(70,20){\circle*{1}} \put(70,20){\circle{2}}

%----------------------- Root 14
%\put(70,20){\line(4,-1){20}}

\put(90,15){\vector(-4,1){19}}

\put(90,15){\circle*{1}}

%\put(92,14){\makebox(0,0)[bl]{\(14\)}}

%--
\put(90,15){\line(1,2){5}}

\put(95,25){\circle*{1}}

%\put(97,24){\makebox(0,0)[bl]{\(9\)}}

%--

\put(90,15){\line(2,1){10}}

\put(100,20){\circle*{1}}

%\put(101,19){\makebox(0,0)[bl]{\(11\)}}

%--
\put(90,15){\line(3,-1){15}}

\put(105,10){\circle*{1}}

%\put(100,09){\makebox(0,0)[bl]{\(16\)}}

%-----------------------------------4->5
\put(60,32.5){\circle*{1}} \put(60,32.5){\circle{2}}

\put(62,30.5){\makebox(0,0)[bl]{\(5\)}}

\put(55,38.5){\makebox(0,0)[bl]{\(6\)}}

\put(60,32.5){\line(-1,1){4}} \put(60,32.5){\line(2,1){8}}

%-----------------4-1
\put(56,36.5){\circle*{1}} \put(56,36.5){\circle{2}}

%\put(56,36.5){\line(-2,1){6}}

\put(50,39.5){\vector(2,-1){5}}

\put(50,39.5){\circle*{1}}

%\put(56,36.5){\line(-1,0){7}}

\put(49,36.5){\vector(1,0){6.5}}

\put(49,36.5){\circle*{1}}

%-----------------4-2
\put(68,36.5){\circle*{1}} \put(68,36.5){\circle{2}}

\put(68,38.5){\makebox(0,0)[bl]{\(7\)}}

%\put(64,36.5){\line(1,0){6}}

\put(74,36.5){\vector(-1,0){5.4}}

%%%%%%%%%%%%%%%%%%%%%%%%%%%%%%%%%%%%%%%%%%%%%%%%%%%%%%%

\put(74,36.5){\circle*{1}}

\put(74,36.5){\line(1,1){4}} \put(78,40.5){\circle*{1}}

\put(74,36.5){\line(1,0){5}} \put(78,36.5){\circle*{1}}

\put(78,36.5){\line(2,1){5}} \put(83,39){\circle*{1}}

\put(78,36.5){\line(2,-1){5}} \put(83,34){\circle*{1}}

\end{picture}
\end{center}

\subsubsection{Connection Assignment in Two-layer Network (Access Points - Users)}

 This section
  \footnote{
 From (with amendments):\\
 M.Sh. Levin,
 Towards communication network development (structural systems
 issues, combinatorial problems).
 {\it IEEE Region 8 Int. Conf. Sibircon 2010},
  vol. 1, (2010) 204--208.
  }
%%%%%%%%%%%%%%%%%
 focuses on assignment design for two-layer network (access points - users)
  (Fig. 25).
 The design problem consists in connection (assignment)
 of each user to access point.

\begin{center}
\begin{picture}(80,32)
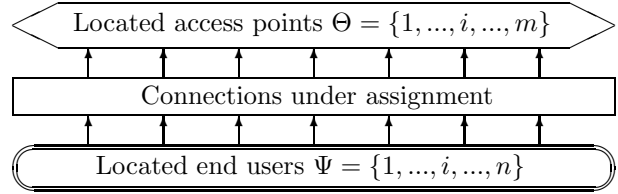
\put(10,00){\makebox(0,0)[bl]{Fig. 25. Layers: access points -
users}}

%------------------------------
%\put(30,28){\oval(60,6)} \put(30,28){\oval(59,5)}

\put(06,25){\line(1,0){68}} \put(06,31){\line(1,0){68}}

\put(00,28){\line(2,1){6}}  \put(80,28){\line(-2,1){6}}
\put(00,28){\line(2,-1){6}}  \put(80,28){\line(-2,-1){6}}

\put(08,26.4){\makebox(0,0)[bl]{Located access points  \( \Theta =
 \{1,...,i,...,m  \} \)
}}

%%%%%%%%%%%%%%%%%%%%%%%%%%%
\put(10,21){\vector(0,1){4}} \put(20,21){\vector(0,1){4}}
\put(30,21){\vector(0,1){4}} \put(40,21){\vector(0,1){4}}
\put(50,21){\vector(0,1){4}} \put(60,21){\vector(0,1){4}}
\put(70,21){\vector(0,1){4}}

%------------------------------
\put(00,16){\line(1,0){80}} \put(00,21){\line(1,0){80}}
\put(00,16){\line(0,1){5}}  \put(80,16){\line(0,1){5}}

\put(17,16.9){\makebox(0,0)[bl]{Connections under assignment}}

%%%%%%%%%%%%%%%%%%%%%%%%%%%
\put(10,12){\vector(0,1){4}} \put(20,12){\vector(0,1){4}}
\put(30,12){\vector(0,1){4}} \put(40,12){\vector(0,1){4}}
\put(50,12){\vector(0,1){4}} \put(60,12){\vector(0,1){4}}
\put(70,12){\vector(0,1){4}}
%------------------------------
\put(40,09){\oval(80,6)} \put(40,09){\oval(79,5)}

\put(11,7.3){\makebox(0,0)[bl]{Located end users
 \(\Psi = \{1,...,i,...,n  \}\)
}}

%------------------------------

%\put(83,34){\circle*{1}}

\end{picture}
\end{center}

 The problem and initial data are partially based on \cite{levpet10a}.
 The following parameters are used:
 set of users \( \Psi = \{1,...,i,...,n  \} \) (\(n=25\)),
 set of access points  \( \Theta = \{1,...,i,...,m  \} \)
 (\(m=6\)).
 Each user is described by  parameter vector
 \((x_{i},y_{i},z_{i},f_{i},p_{i} )\),
 where vector components are as follows
 (Table 3):
 coordinates of user \((x_{i},y_{i},z_{i}\),
 required frequency bandwidth \(f_{i}\)
 (scale: 1 Mbit/s ... 10 Mbit/s),
 priority \(p_{i}\)
 (ordinal scale [1,2,3], all user requirements are satisfied in case
 \(p_{i} = 1\)),
 required reliability \(r_{i}\)
  (ordinal scale [1,10],
 \(10\) corresponds to maximum reliability).
%
% required reliability \(d_{i}\)
%  (ordinal scale [1,10],
% \(10\) corresponds to the highest level of data protection).
%
%
 Analogically,
 parameters of access points are considered
 (by index \(j\), Table 4) including parameter \(n_{j}\)
 (maximal possible number of users under service).

\begin{center}
\begin{picture}(48,68)

\put(26,63){\makebox(0,0)[bl]{Table 3. Data on users
 \cite{levsib10}}}

\put(00,00){\line(1,0){48}} \put(00,54){\line(1,0){48}}
\put(00,61){\line(1,0){48}}

\put(00,00){\line(0,1){61}} \put(06,00){\line(0,1){61}}
\put(13,00){\line(0,1){61}} \put(20,00){\line(0,1){61}}
\put(27,00){\line(0,1){61}} \put(34,00){\line(0,1){61}}
\put(41,00){\line(0,1){61}} \put(48,00){\line(0,1){61}}

\put(02,56){\makebox(0,0)[bl]{\(i\)}}

\put(08,56){\makebox(0,0)[bl]{\(x_{i}\)}}
\put(15,56){\makebox(0,0)[bl]{\(y_{i}\)}}
\put(22,56){\makebox(0,0)[bl]{\(z_{i}\)}}

\put(29,56){\makebox(0,0)[bl]{\(f_{i}\)}}

\put(36,56){\makebox(0,0)[bl]{\(p_{i}\)}}
\put(43,56){\makebox(0,0)[bl]{\(r_{i}\)}}
%

%----------------------------------------
%%98 - 54 = 44
%--1
\put(02,50){\makebox(0,0)[bl]{\(1\)}}
\put(08,50){\makebox(0,0)[bl]{\(30\)}}
\put(14,50){\makebox(0,0)[bl]{\(165\)}}
\put(22.5,50){\makebox(0,0)[bl]{\(5\)}}

\put(28.5,50){\makebox(0,0)[bl]{\(10\)}}

\put(36.5,50){\makebox(0,0)[bl]{\(2\)}}
\put(43.5,50){\makebox(0,0)[bl]{\(5\)}}

%--2
\put(02,46){\makebox(0,0)[bl]{\(2\)}}
\put(08,46){\makebox(0,0)[bl]{\(58\)}}
\put(14,46){\makebox(0,0)[bl]{\(174\)}}
\put(22.5,46){\makebox(0,0)[bl]{\(5\)}}

\put(29.5,46){\makebox(0,0)[bl]{\(5\)}}

\put(36.5,46){\makebox(0,0)[bl]{\(1\)}}
\put(43.5,46){\makebox(0,0)[bl]{\(9\)}}

%--3
\put(02,42){\makebox(0,0)[bl]{\(3\)}}
\put(08,42){\makebox(0,0)[bl]{\(95\)}}
\put(14,42){\makebox(0,0)[bl]{\(156\)}}
\put(22.5,42){\makebox(0,0)[bl]{\(0\)}}
\put(29.5,42){\makebox(0,0)[bl]{\(6\)}}

%\put(36.5,46){\makebox(0,0)[bl]{\(1\)}}

\put(36.5,42){\makebox(0,0)[bl]{\(1\)}}
\put(43.5,42){\makebox(0,0)[bl]{\(6\)}}

%--8->7-4
\put(02,38){\makebox(0,0)[bl]{\(4\)}}
\put(08,38){\makebox(0,0)[bl]{\(52\)}}
\put(14,38){\makebox(0,0)[bl]{\(134\)}}
\put(22.5,38){\makebox(0,0)[bl]{\(5\)}}
\put(29.5,38){\makebox(0,0)[bl]{\(6\)}}

\put(36.5,38){\makebox(0,0)[bl]{\(1\)}}
\put(43.5,38){\makebox(0,0)[bl]{\(8\)}}

%--9->8
\put(02,34){\makebox(0,0)[bl]{\(5\)}}
\put(08,34){\makebox(0,0)[bl]{\(85\)}}
\put(14,34){\makebox(0,0)[bl]{\(134\)}}
\put(22.5,34){\makebox(0,0)[bl]{\(3\)}}
\put(29.5,34){\makebox(0,0)[bl]{\(6\)}}

\put(36.5,34){\makebox(0,0)[bl]{\(1\)}}
\put(43.5,34){\makebox(0,0)[bl]{\(7\)}}

%--13->12
\put(02,30){\makebox(0,0)[bl]{\(6\)}}
\put(08,30){\makebox(0,0)[bl]{\(27\)}}
\put(14,30){\makebox(0,0)[bl]{\(109\)}}
\put(22.5,30){\makebox(0,0)[bl]{\(7\)}}
\put(29.5,30){\makebox(0,0)[bl]{\(8\)}}

\put(36.5,30){\makebox(0,0)[bl]{\(3\)}}
\put(43.5,30){\makebox(0,0)[bl]{\(5\)}}

%--14->13
\put(02,26){\makebox(0,0)[bl]{\(7\)}}
\put(08,26){\makebox(0,0)[bl]{\(55\)}}
\put(14,26){\makebox(0,0)[bl]{\(105\)}}
\put(22.5,26){\makebox(0,0)[bl]{\(2\)}}
\put(29.5,26){\makebox(0,0)[bl]{\(7\)}}

\put(36.5,26){\makebox(0,0)[bl]{\(2\)}}
\put(42.5,26){\makebox(0,0)[bl]{\(10\)}}

%--15->14
\put(02,22){\makebox(0,0)[bl]{\(8\)}}
\put(08,22){\makebox(0,0)[bl]{\(98\)}}
\put(15,22){\makebox(0,0)[bl]{\(89\)}}
\put(22.5,22){\makebox(0,0)[bl]{\(3\)}}
\put(28.5,22){\makebox(0,0)[bl]{\(10\)}}

\put(36.5,22){\makebox(0,0)[bl]{\(1\)}}
\put(42.5,22){\makebox(0,0)[bl]{\(10\)}}

%--20->16
\put(02,18){\makebox(0,0)[bl]{\(9\)}}
\put(08,18){\makebox(0,0)[bl]{\(25\)}}
\put(15,18){\makebox(0,0)[bl]{\(65\)}}
\put(22.5,18){\makebox(0,0)[bl]{\(2\)}}
\put(29.5,18){\makebox(0,0)[bl]{\(7\)}}

\put(36.5,18){\makebox(0,0)[bl]{\(3\)}}
\put(43.5,18){\makebox(0,0)[bl]{\(5\)}}

%--21->17
\put(01,14){\makebox(0,0)[bl]{\(10\)}}
\put(08,14){\makebox(0,0)[bl]{\(52\)}}
\put(15,14){\makebox(0,0)[bl]{\(81\)}}
\put(22.5,14){\makebox(0,0)[bl]{\(1\)}}
\put(28.5,14){\makebox(0,0)[bl]{\(10\)}}

\put(36.5,14){\makebox(0,0)[bl]{\(1\)}}
\put(43.5,14){\makebox(0,0)[bl]{\(8\)}}

%--28->20
\put(01,10){\makebox(0,0)[bl]{\(11\)}}
\put(08,10){\makebox(0,0)[bl]{\(65\)}}
\put(15,10){\makebox(0,0)[bl]{\(25\)}}
\put(22.5,10){\makebox(0,0)[bl]{\(7\)}}
\put(29.5,10){\makebox(0,0)[bl]{\(6\)}}

\put(36.5,10){\makebox(0,0)[bl]{\(2\)}}
\put(43.5,10){\makebox(0,0)[bl]{\(9\)}}

%--29->21
\put(01,06){\makebox(0,0)[bl]{\(12\)}}
\put(08,06){\makebox(0,0)[bl]{\(93\)}}
\put(15,06){\makebox(0,0)[bl]{\(39\)}}
\put(22.5,06){\makebox(0,0)[bl]{\(1\)}}
\put(28.5,06){\makebox(0,0)[bl]{\(10\)}}

\put(36.5,06){\makebox(0,0)[bl]{\(1\)}}
\put(42.5,06){\makebox(0,0)[bl]{\(10\)}}
%\put(57.5,06){\makebox(0,0)[bl]{\(9\)}}

%--new 23
\put(01,02){\makebox(0,0)[bl]{\(13\)}}
\put(07,02){\makebox(0,0)[bl]{\(172\)}}
\put(15,02){\makebox(0,0)[bl]{\(26\)}}
\put(22.5,02){\makebox(0,0)[bl]{\(2\)}}
\put(28.5,02){\makebox(0,0)[bl]{\(10\)}}

\put(36.5,02){\makebox(0,0)[bl]{\(2\)}}
\put(43.5,02){\makebox(0,0)[bl]{\(7\)}}

\end{picture}
%\end{center}
%
%
%\begin{center}
\begin{picture}(48,67)

\put(00,00){\line(1,0){48}} \put(00,50){\line(1,0){48}}
\put(00,57){\line(1,0){48}}

\put(00,00){\line(0,1){57}} \put(06,00){\line(0,1){57}}
\put(13,00){\line(0,1){57}} \put(20,00){\line(0,1){57}}
\put(27,00){\line(0,1){57}} \put(34,00){\line(0,1){57}}
\put(41,00){\line(0,1){57}} \put(48,00){\line(0,1){57}}
%\put(55,00){\line(0,1){53}} \put(62,00){\line(0,1){53}}

\put(02,52){\makebox(0,0)[bl]{\(i\)}}

\put(08,52){\makebox(0,0)[bl]{\(x_{i}\)}}
\put(15,52){\makebox(0,0)[bl]{\(y_{i}\)}}
\put(22,52){\makebox(0,0)[bl]{\(z_{i}\)}}

\put(29,52){\makebox(0,0)[bl]{\(f_{i}\)}}

\put(36,52){\makebox(0,0)[bl]{\(p_{i}\)}}
\put(43,52){\makebox(0,0)[bl]{\(r_{i}\)}}
%

%----------------------------------------
%--4
\put(01,46){\makebox(0,0)[bl]{\(14\)}}
\put(07,46){\makebox(0,0)[bl]{\(110\)}}
\put(14,46){\makebox(0,0)[bl]{\(169\)}}
\put(22.5,46){\makebox(0,0)[bl]{\(5\)}}
\put(29.5,46){\makebox(0,0)[bl]{\(7\)}}

%\put(36.5,42){\makebox(0,0)[bl]{\(1\)}}

\put(36.5,46){\makebox(0,0)[bl]{\(2\)}}
\put(43.5,46){\makebox(0,0)[bl]{\(5\)}}
%\put(57.5,42){\makebox(0,0)[bl]{\(6\)}}

%--5
\put(01,42){\makebox(0,0)[bl]{\(15\)}}
\put(07,42){\makebox(0,0)[bl]{\(145\)}}
\put(14,42){\makebox(0,0)[bl]{\(181\)}}
\put(22.5,42){\makebox(0,0)[bl]{\(3\)}}
\put(29.5,42){\makebox(0,0)[bl]{\(5\)}}

%\put(36.5,38){\makebox(0,0)[bl]{\(1\)}}

\put(36.5,42){\makebox(0,0)[bl]{\(2\)}}
\put(43.5,42){\makebox(0,0)[bl]{\(4\)}}
%\put(57.5,38){\makebox(0,0)[bl]{\(6\)}}

%--6
\put(01,38){\makebox(0,0)[bl]{\(16\)}}
\put(07,38){\makebox(0,0)[bl]{\(170\)}}
\put(14,38){\makebox(0,0)[bl]{\(161\)}}
\put(22.5,38){\makebox(0,0)[bl]{\(5\)}}
\put(29.5,38){\makebox(0,0)[bl]{\(7\)}}

%\put(36.5,34){\makebox(0,0)[bl]{\(1\)}}

\put(36.5,38){\makebox(0,0)[bl]{\(2\)}}
\put(43.5,38){\makebox(0,0)[bl]{\(4\)}}
%\put(57.5,34){\makebox(0,0)[bl]{\(7\)}}

%--10->9
\put(01,34){\makebox(0,0)[bl]{\(17\)}}
\put(07,34){\makebox(0,0)[bl]{\(120\)}}
\put(14,34){\makebox(0,0)[bl]{\(140\)}}
\put(22.5,34){\makebox(0,0)[bl]{\(6\)}}
\put(29.5,34){\makebox(0,0)[bl]{\(4\)}}

%\put(36.5,30){\makebox(0,0)[bl]{\(1\)}}

\put(36.5,34){\makebox(0,0)[bl]{\(2\)}}
\put(43.5,34){\makebox(0,0)[bl]{\(6\)}}
%\put(57.5,30){\makebox(0,0)[bl]{\(8\)}}

%--11->10
\put(01,30){\makebox(0,0)[bl]{\(18\)}}
\put(07,30){\makebox(0,0)[bl]{\(150\)}}
\put(14,30){\makebox(0,0)[bl]{\(136\)}}
\put(22.5,30){\makebox(0,0)[bl]{\(3\)}}
\put(29.5,30){\makebox(0,0)[bl]{\(6\)}}

%\put(36.5,26){\makebox(0,0)[bl]{\(1\)}}

\put(36.5,30){\makebox(0,0)[bl]{\(2\)}}
\put(43.5,30){\makebox(0,0)[bl]{\(7\)}}
%\put(57.5,26){\makebox(0,0)[bl]{\(8\)}}

%--12->11
\put(01,26){\makebox(0,0)[bl]{\(19\)}}
\put(07,26){\makebox(0,0)[bl]{\(175\)}}
\put(14,26){\makebox(0,0)[bl]{\(125\)}}
\put(22.5,26){\makebox(0,0)[bl]{\(1\)}}
\put(29.5,26){\makebox(0,0)[bl]{\(8\)}}

%\put(36.5,22){\makebox(0,0)[bl]{\(1\)}}

\put(36.5,26){\makebox(0,0)[bl]{\(3\)}}
\put(43.5,26){\makebox(0,0)[bl]{\(5\)}}
%\put(57.5,22){\makebox(0,0)[bl]{\(6\)}}

%--19->15
\put(01,22){\makebox(0,0)[bl]{\(20\)}}
\put(07,22){\makebox(0,0)[bl]{\(183\)}}
\put(15,22){\makebox(0,0)[bl]{\(91\)}}
\put(22.5,22){\makebox(0,0)[bl]{\(4\)}}
\put(29.5,22){\makebox(0,0)[bl]{\(4\)}}

%\put(36.5,22){\makebox(0,0)[bl]{\(1\)}}

\put(36.5,22){\makebox(0,0)[bl]{\(3\)}}
\put(43.5,22){\makebox(0,0)[bl]{\(5\)}}
%\put(57.5,22){\makebox(0,0)[bl]{\(7\)}}

%--24->18
\put(01,18){\makebox(0,0)[bl]{\(21\)}}
\put(07,18){\makebox(0,0)[bl]{\(135\)}}
\put(15,18){\makebox(0,0)[bl]{\(59\)}}
\put(22.5,18){\makebox(0,0)[bl]{\(4\)}}
\put(28.5,18){\makebox(0,0)[bl]{\(13\)}}

%\put(36.5,18){\makebox(0,0)[bl]{\(1\)}}

\put(36.5,18){\makebox(0,0)[bl]{\(3\)}}
\put(43.5,18){\makebox(0,0)[bl]{\(4\)}}
%\put(57.5,18){\makebox(0,0)[bl]{\(6\)}}

%--25->19
\put(01,14){\makebox(0,0)[bl]{\(22\)}}
\put(07,14){\makebox(0,0)[bl]{\(147\)}}
\put(15,14){\makebox(0,0)[bl]{\(79\)}}
\put(22.5,14){\makebox(0,0)[bl]{\(5\)}}
\put(29.5,14){\makebox(0,0)[bl]{\(7\)}}

%\put(36.5,14){\makebox(0,0)[bl]{\(1\)}}

\put(36.5,14){\makebox(0,0)[bl]{\(3\)}}
\put(42.5,14){\makebox(0,0)[bl]{\(16\)}}
%\put(57.5,14){\makebox(0,0)[bl]{\(8\)}}

%--30->22
\put(01,10){\makebox(0,0)[bl]{\(23\)}}
\put(07,10){\makebox(0,0)[bl]{\(172\)}}
\put(15,10){\makebox(0,0)[bl]{\(26\)}}
\put(22.5,10){\makebox(0,0)[bl]{\(2\)}}
\put(28.5,10){\makebox(0,0)[bl]{\(10\)}}

%\put(36.5,10){\makebox(0,0)[bl]{\(1\)}}

\put(36.5,10){\makebox(0,0)[bl]{\(2\)}}
\put(43.5,10){\makebox(0,0)[bl]{\(7\)}}
%\put(57.5,10){\makebox(0,0)[bl]{\(6\)}}

%--new 24
\put(01,06){\makebox(0,0)[bl]{\(24\)}}
\put(07,06){\makebox(0,0)[bl]{\(165\)}}
\put(15,06){\makebox(0,0)[bl]{\(50\)}}
\put(22.5,06){\makebox(0,0)[bl]{\(3\)}}
\put(28.5,06){\makebox(0,0)[bl]{\(7\)}}

%\put(36.5,06){\makebox(0,0)[bl]{\(1\)}}

\put(36.5,06){\makebox(0,0)[bl]{\(3\)}}
\put(43.5,06){\makebox(0,0)[bl]{\(3\)}}
%\put(57.5,06){\makebox(0,0)[bl]{\(7\)}}

%--new 25
\put(01,02){\makebox(0,0)[bl]{\(25\)}}
\put(07,02){\makebox(0,0)[bl]{\(127\)}}
\put(15,02){\makebox(0,0)[bl]{\(95\)}}
\put(22.5,02){\makebox(0,0)[bl]{\(5\)}}
\put(28.5,02){\makebox(0,0)[bl]{\(7\)}}

%\put(36.5,02){\makebox(0,0)[bl]{\(1\)}}

\put(36.5,02){\makebox(0,0)[bl]{\(2\)}}
\put(43.5,02){\makebox(0,0)[bl]{\(5\)}}
%\put(57.5,02){\makebox(0,0)[bl]{\(7\)}}

%------------------------------------------------

\end{picture}
\end{center}

\begin{center}
\begin{picture}(42,28)

\put(15.5,23){\makebox(0,0)[bl]{Table 4. Data on access points
 \cite{levsib10}}}

\put(00,00){\line(1,0){42}} \put(00,14){\line(1,0){42}}
\put(00,21){\line(1,0){42}}

\put(00,00){\line(0,1){21}} \put(04,00){\line(0,1){21}}
\put(11,00){\line(0,1){21}} \put(18,00){\line(0,1){21}}
\put(24,00){\line(0,1){21}} \put(30,00){\line(0,1){21}}
\put(36,00){\line(0,1){21}} \put(42,00){\line(0,1){21}}
%\put(48,00){\line(0,1){21}}

\put(01,16){\makebox(0,0)[bl]{\(j\)}}
\put(05.5,16){\makebox(0,0)[bl]{\(x_{j}\)}}
\put(12.5,16){\makebox(0,0)[bl]{\(y_{j}\)}}
\put(19,16){\makebox(0,0)[bl]{\(z_{j}\)}}
\put(25,16){\makebox(0,0)[bl]{\(f_{j}\)}}
\put(31,16){\makebox(0,0)[bl]{\(n_{j}\)}}
\put(37,16){\makebox(0,0)[bl]{\(r_{j}\)}}
%

%--

\put(01,10){\makebox(0,0)[bl]{\(1\)}}
\put(05.5,10){\makebox(0,0)[bl]{\(50\)}}
\put(11.5,10){\makebox(0,0)[bl]{\(157\)}}
\put(19,10){\makebox(0,0)[bl]{\(10\)}}
\put(25,10){\makebox(0,0)[bl]{\(30\)}}
\put(32,10){\makebox(0,0)[bl]{\(4\)}}
\put(37,10){\makebox(0,0)[bl]{\(10\)}}
%\put(43,10){\makebox(0,0)[bl]{\(10\)}}

%--

\put(01,06){\makebox(0,0)[bl]{\(2\)}}
\put(05.5,06){\makebox(0,0)[bl]{\(72\)}}
\put(11.5,06){\makebox(0,0)[bl]{\(102\)}}
\put(19,06){\makebox(0,0)[bl]{\(10\)}}
\put(25,06){\makebox(0,0)[bl]{\(42\)}}
\put(32,06){\makebox(0,0)[bl]{\(6\)}}
\put(37,06){\makebox(0,0)[bl]{\(10\)}}
%\put(44,06){\makebox(0,0)[bl]{\(9\)}}

%--=================================================

\put(01,02){\makebox(0,0)[bl]{\(3\)}}
\put(05.5,02){\makebox(0,0)[bl]{\(45\)}}
\put(12.5,02){\makebox(0,0)[bl]{\(52\)}}
\put(19,02){\makebox(0,0)[bl]{\(10\)}}
\put(25,02){\makebox(0,0)[bl]{\(45\)}}
\put(31,02){\makebox(0,0)[bl]{\(10\)}}
\put(37,02){\makebox(0,0)[bl]{\(10\)}}
%\put(43,02){\makebox(0,0)[bl]{\(10\)}}

%------------------------------------------------

\end{picture}
%\end{center}
%
%\begin{center}
\begin{picture}(42,28)

\put(00,00){\line(1,0){42}} \put(00,14){\line(1,0){42}}
\put(00,21){\line(1,0){42}}

\put(00,00){\line(0,1){21}} \put(04,00){\line(0,1){21}}
\put(11,00){\line(0,1){21}} \put(18,00){\line(0,1){21}}
\put(24,00){\line(0,1){21}} \put(30,00){\line(0,1){21}}
\put(36,00){\line(0,1){21}} \put(42,00){\line(0,1){21}}
%\put(48,00){\line(0,1){21}}

\put(01,16){\makebox(0,0)[bl]{\(j\)}}
\put(05.5,16){\makebox(0,0)[bl]{\(x_{j}\)}}
\put(12.5,16){\makebox(0,0)[bl]{\(y_{j}\)}}
\put(19,16){\makebox(0,0)[bl]{\(z_{j}\)}}
\put(25,16){\makebox(0,0)[bl]{\(f_{j}\)}}
\put(31,16){\makebox(0,0)[bl]{\(n_{j}\)}}
\put(37,16){\makebox(0,0)[bl]{\(r_{j}\)}}
%

%--

\put(01,10){\makebox(0,0)[bl]{\(4\)}}
\put(04.5,10){\makebox(0,0)[bl]{\(150\)}}
\put(11.5,10){\makebox(0,0)[bl]{\(165\)}}
\put(19,10){\makebox(0,0)[bl]{\(10\)}}
\put(25,10){\makebox(0,0)[bl]{\(30\)}}
\put(32,10){\makebox(0,0)[bl]{\(5\)}}
\put(37,10){\makebox(0,0)[bl]{\(15\)}}
%\put(44,10){\makebox(0,0)[bl]{\(8\)}}

%--

\put(01,06){\makebox(0,0)[bl]{\(5\)}}
\put(04.5,06){\makebox(0,0)[bl]{\(140\)}}
\put(11.5,06){\makebox(0,0)[bl]{\(112\)}}
\put(19,06){\makebox(0,0)[bl]{\(10\)}}
\put(25,06){\makebox(0,0)[bl]{\(32\)}}
\put(32,06){\makebox(0,0)[bl]{\(5\)}}
\put(38,06){\makebox(0,0)[bl]{\(8\)}}
%\put(44,06){\makebox(0,0)[bl]{\(8\)}}

%--=================================================

\put(01,02){\makebox(0,0)[bl]{\(6\)}}
\put(04.5,02){\makebox(0,0)[bl]{\(147\)}}
\put(12.5,02){\makebox(0,0)[bl]{\(47\)}}
\put(19,02){\makebox(0,0)[bl]{\(10\)}}
\put(25,02){\makebox(0,0)[bl]{\(30\)}}
\put(32,02){\makebox(0,0)[bl]{\(5\)}}
\put(37,02){\makebox(0,0)[bl]{\(15\)}}
%\put(44,02){\makebox(0,0)[bl]{\(7\)}}
%------------------------------------------------

\end{picture}
\end{center}

 Further, each pair ``user-access point'' (i.e.,
 \((i,j),  i \in \Psi, j \in \Theta \)) can be described:
 (1) reliability ~\(r_{ij} = min  \{r_{i},r_{j}\}\),
 (2) distance ~\(l_{ij}\),
 (3) priority ~\(p_{ij} = p_{i}\), and
 (4) required bandwidth  ~\(f_{ij} = f_{i}\).
 In addition, a ``connection'' parameter is considered:~
   \(\beta_{ij}\) equals \(1\) if \(l_{ij} \leq l\) and \(0\)
   otherwise (\(L\) corresponds to distance constraint).
 This parameter defines
 \( \forall i \in \Psi \)
 a subset of possible access points
 \(\Theta_{i} \subseteq \Theta\).
%
%
%%%%%%%%%%%%%%%%%%%%%%%%%%%%%%%%%%%%%%%%%%%%%%%%%%%%%%%
%
 The assignment of user \( i\) to access point \( j\)
 is defined by Boolean variable
 \(x_{ij}\) (\(x_{ij} =1\) in the case of assignment
  \( i\) to  \( j\)  and   \(x_{ij} =0\) otherwise).
 The assignment solution
  ~(\( \Psi \) \(\Rightarrow \)  ~\( \Theta \))~
 is defined by Boolean matrix
 ~\(X = || x_{ij} ||, ~ i=\overline{1,n}, ~j=\overline{1,m}\).
%
%%%%%%%%%%%%%%%%%%%%%%%%%%%%%%%%%%%%%%%%%%%%
%
 The problem formulation is:
 \[\max \sum_{i=1}^{n}   \sum_{j \in \Theta_{i}}  r_{ij} x_{ij},
 ~ \max \sum_{i=1}^{n}   \sum_{j \in \Theta_{i}}  f_{ij}
 x_{ij},
%%%%%%%%%%%%%%%%%%%%%%%%%%%%%%%%%%%%%%%%%%%%%%%%%%%%%%%
%
 ~ \max \sum_{i=1}^{n}   \sum_{j \in \Theta_{i}}  p_{ij}
 x_{ij} \]
%%%%%%%%%%%%%%%%%%%%%%%%%%%%%%%%%%%%%%%%%%%%%%%%%%%%%%%
%
 \[s.t.~
  \sum_{i=1}^{n}  f_{ij} x_{ij}  \leq  f_{j}
 ~ \forall j \in \Theta,
   ~ ~\sum_{i=1}^{n}  x_{ij}  \leq  n_{j}
 ~ \forall j \in \Theta,
%%%%%%%%%%%%%%%%%%%%%%%%%%%%%%%%%%%%%%%%%%%%%%%%%%%%%%%%
%
 ~ \sum_{j \in \Theta_{i}} x_{ij}\leq 1 ~\forall i \in \Psi ,\]
  \[ x_{ij} \in \{0,1\}, ~~\forall ~i=\overline{1,n},  ~\forall
  ~j=\overline{1,m};
%%%%%%%%%%%%%%%%%%%%%%%%%%%%%%%%%%%%%%%%%%%%%%%%%%%%%%%%%%%%%%%%%
%
  ~~ x_{ij} = 0, ~~\forall ~i=\overline{1,n},   ~j \in  \{\Theta \setminus \Theta_{i}
  \}. \]
 This multiple criteria assignment problem is NP-hard
 (e.g., \cite{gar79}).
 Here, a simplified two-stage heuristic is used:
 (i) transformation of vector estimate for each pair \((i,j)\)
 into an ordinal estimate (by multicriteria ranking, ELECTRE-like
 technique
  (\cite{lev12b}, \cite{roy96}),
 (ii) solving the obtained one-criterion assignment problem
 (by greedy algorithm).
 Fig. 26 depicts  the obtained solution:  assignment of users to access points.

\begin{center}
\begin{picture}(74,80)
\put(00,00){\makebox(0,0)[bl]{Fig. 26.
 Assignment:  users - access points
 \cite{levsib10}}}

%-----------------------------------------------

\put(14,32){\oval(3.2,4)} \put(14,32){\circle*{2}}
\put(12,27){\makebox(0,0)[bl]{\(10\)}}
%21->17

\put(15,42){\oval(3.2,4)} \put(15,42){\circle*{2}}
\put(13,37){\makebox(0,0)[bl]{\(7\)}}
%14->13

%------------------------------------------------Point 3-2
\put(19,39){\line(1,0){6}} \put(19,39){\line(1,2){3}}
\put(25,39){\line(-1,2){3}} \put(22,45){\line(0,1){3}}
\put(22,48){\circle*{1}} \put(21,39.5){\makebox(0,0)[bl]{\(2\)}}

\put(22,48){\circle{2.5}} \put(22,48){\circle{3.5}}
%---
%--
\put(24.5,40){\line(2,-1){6}}
%--
%
\put(23,42.3){\line(1,2){3}} \put(26,48.3){\line(0,1){4.8}}
%--
%
\put(20.5,42){\line(-1,0){4.5}}

\put(05,43){\oval(3.2,4)} \put(05,43){\circle*{2}}
\put(03,38){\makebox(0,0)[bl]{\(6\)}}
%13->12

%--
%connection 1-12
\put(11,60){\line(-1,-3){6}}
%--

\put(05,67){\oval(3.2,4)} \put(05,67){\circle*{2}}
\put(03.8,62){\makebox(0,0)[bl]{\(1\)}}

\put(15,74){\oval(3.2,4)} \put(15,74){\circle*{2}}
\put(13.8,69){\makebox(0,0)[bl]{\(2\)}}

%-------------------------------------------------------Point 1
\put(09,60){\line(1,0){6}} \put(09,60){\line(1,2){3}}
\put(15,60){\line(-1,2){3}} \put(12,66){\line(0,1){3}}
\put(12,69){\circle*{1}} \put(11,60.5){\makebox(0,0)[bl]{\(1\)}}

\put(12,69){\circle{2.5}} \put(12,69){\circle{3.5}}
%---
%
%--
%
\put(10,62){\line(-1,1){5}}
%--
%
\put(13.5,60){\line(0,-1){6}}
%--
%
\put(14,62.3){\line(1,2){2}} \put(16,66.3){\line(0,1){5.8}}

\put(14,53){\oval(3.2,4)} \put(14,53){\circle*{2}}
\put(12.8,48){\makebox(0,0)[bl]{\(4\)}}
%8->7

\put(04,27){\oval(3.2,4)} \put(04,27){\circle*{2}}
\put(02,22){\makebox(0,0)[bl]{\(9\)}}
%20->16

\put(18,10){\oval(3.2,4)} \put(18,10){\circle*{2}}
\put(16,5){\makebox(0,0)[bl]{\(11\)}}
%28->20

\put(28,14){\oval(3.2,4)} \put(28,14){\circle*{2}}
\put(26,9){\makebox(0,0)[bl]{\(12\)}}
%29->21

\put(30,38){\oval(3.2,4)} \put(30,38){\circle*{2}}
\put(28,33){\makebox(0,0)[bl]{\(8\)}}
%15->14

\put(38,60){\oval(3.2,4)} \put(38,60){\circle*{2}}
\put(36,55){\makebox(0,0)[bl]{\(17\)}}
%10->9

\put(50,58){\oval(3.2,4)} \put(50,58){\circle*{2}}
\put(49,53){\makebox(0,0)[bl]{\(18\)}}
%11->10

%----------------------------------------------------------Point 4-5
\put(42,46){\line(1,0){6}} \put(42,46){\line(1,2){3}}
\put(48,46){\line(-1,2){3}} \put(45,52){\line(0,1){3}}
\put(45,55){\circle*{1}} \put(44,46.5){\makebox(0,0)[bl]{\(5\)}}

\put(45,55){\circle{2.5}} \put(45,55){\circle{3.5}}
%---
%Connection 4-9, 4-11, 4-12
%--
%connection 4-9
%\put(43.5,48.5){\line(-3,1){4.5}} \put(39,50){\line(0,1){8}}

\put(43,48){\line(-1,3){4}}

%--
%connection 4-11
\put(46,49.4){\line(1,2){4.5}}
%--
%connection 4-12
\put(47,48){\line(3,1){13}}
%--

\put(60,52){\oval(3.2,4)} \put(60,52){\circle*{2}}
\put(58,47){\makebox(0,0)[bl]{\(19\)}}
%12->11

\put(63,35){\oval(3.2,4)} \put(63,35){\circle*{2}}
\put(61,30){\makebox(0,0)[bl]{\(20\)}}
%19->15

\put(49,30){\oval(3.2,4)} \put(49,30){\circle*{2}}
\put(47,25){\makebox(0,0)[bl]{\(22\)}}
%25->19

\put(44,20){\oval(3.2,4)} \put(44,20){\circle*{2}}
\put(42,15){\makebox(0,0)[bl]{\(21\)}}
%24->18

\put(61,10){\oval(3.2,4)} \put(61,10){\circle*{2}}
\put(59,05){\makebox(0,0)[bl]{\(23\)}}
%30->22

\put(61,20){\oval(3.2,4)} \put(61,20){\circle*{2}}
\put(59,15){\makebox(0,0)[bl]{\(24\)}}

\put(42,37){\oval(3.2,4)} \put(42,37){\circle*{2}}
\put(40,32){\makebox(0,0)[bl]{\(25\)}}

%--

\put(42,37){\line(-1,0){6}} \put(36,37){\line(-1,1){4}}
\put(32,41){\line(-1,0){8}}

%---------------------------------------------------------Point 6
\put(47,13){\line(1,0){6}} \put(47,13){\line(1,2){3}}
\put(53,13){\line(-1,2){3}} \put(50,19){\line(0,1){3}}
\put(50,22){\circle*{1}} \put(49,13.5){\makebox(0,0)[bl]{\(6\)}}

\put(50,22){\circle{2.5}} \put(50,22){\circle{3.5}}
%---
%
%--
%
\put(61,20){\line(-2,-1){9}}
%--
%
\put(48,15.5){\line(-1,1){4.5}}
%--
%
\put(48.7,17){\line(-1,2){4.2}} \put(44.5,25.4){\line(1,1){4}}
%--
%
\put(51,16){\line(1,2){07.5}} \put(58.5,31){\line(1,1){04}}
%--
%
\put(52,15){\line(2,-1){8.6}}
%--

\put(25,55){\oval(3.2,4)} \put(25,55){\circle*{2}}
\put(24,50){\makebox(0,0)[bl]{\(5\)}}
%9->8
%%%%%%%%%%%%%%%%%%%%%%%%%%%%%%%%%%%%%%%%%%%%%%%%%%%%%%%%%%%%%%%%
\put(29,66){\oval(3.2,4)} \put(29,66){\circle*{2}}
\put(28,61){\makebox(0,0)[bl]{\(3\)}}
%--
%
\put(29,66){\line(1,0){15}} \put(44,66){\line(2,1){4}}

\put(35,72){\oval(3.2,4)} \put(35,72){\circle*{2}}
\put(34,67){\makebox(0,0)[bl]{\(14\)}}

\put(45,77){\oval(3.2,4)} \put(45,77){\circle*{2}}
\put(44,72){\makebox(0,0)[bl]{\(15\)}}

\put(58,67){\oval(3.2,4)} \put(58,67){\circle*{2}}
\put(57,62){\makebox(0,0)[bl]{\(16\)}}

%----------------------------------------------------------Point 2-4
\put(47,68){\line(1,0){6}} \put(47,68){\line(1,2){3}}
\put(53,68){\line(-1,2){3}} \put(50,74){\line(0,1){3}}
\put(50,77){\circle*{1}} \put(49,68.5){\makebox(0,0)[bl]{\(4\)}}

\put(50,77){\circle{2.5}} \put(50,77){\circle{3.5}}
%---
%Connection 2-20, 2-21, 2-29, 2-28
%--
%connection 2-2
%\put(48,70.3){\line(-3,1){14}} \put(34,75){\line(-1,0){18}}
%--
%connection 2-3
%\put(48,68){\line(-2,-1){04}} \put(44,66){\line(-1,0){17}}
%--
%connection 2-4
\put(48,70){\line(-1,0){09}} \put(39,70){\line(-2,1){5}}
%--
%connection 2-5
\put(48.4,71.2){\line(-1,2){3}}
%--
%connection 2-6
\put(52,70){\line(2,-1){5}}
%--

%---------------------------------------------------------Point 5-3
\put(09,16){\line(1,0){6}} \put(09,16){\line(1,2){3}}
\put(15,16){\line(-1,2){3}} \put(12,22){\line(0,1){3}}
\put(12,25){\circle*{1}} \put(11,16.5){\makebox(0,0)[bl]{\(3\)}}

\put(12,25){\circle{2.5}} \put(12,25){\circle{3.5}}
%------------------------------------------------
%Connection 5-20, 5-21, 5-29, 5-28
%--
%connection 5-20
\put(10,18){\line(-2,3){6.5}}
%--
%connection 5-21
\put(13.5,18.5){\line(1,2){05}} \put(18.5,28.5){\line(-1,1){04}}
%--
%connection 5-29
\put(14,18){\line(4,-1){14}}
%--
%connection 5-28
\put(14,16){\line(1,-2){3}}
%--

%%%%%%%%%%%%%%%%%%%%%%%%%%%%%%%%%%%%%%%%%%%%%%%%%%%%%%%%%%%%%%%%%%
%%%%%%%%%%%%%%%%%%%%%%%%%%%%%%%%%%%%%%%%%%%%%%%%%%%%%%%% Additions

\put(36,19){\oval(3.2,4)} \put(36,19){\circle*{2}}
\put(34,13){\makebox(0,0)[bl]{\(13\)}}
%connection 5-23 6-23
%\put(34.4,18.6){\line(-1,0){20.6}}

\put(36,19){\line(1,0){3}} \put(39,19){\line(1,-2){2.5}}
\put(41.5,14){\line(1,0){6}}

%--

\end{picture}
\end{center}

%%%%%%%%%%%%%%%%%%%%%%%%%%%%%%%%%%%%%%%%%%%%%%%%%%%%%
\subsection{Morphological Hierarchy}

 Morphological hierarchy generalizes
 system model from morphological analysis
 \cite{zwi69} by the following ways
 (\cite{lev98},\cite{lev06},\cite{lev09},\cite{lev11agg},\cite{lev12morph},\cite{lev12a}):
 (i) hierarchical (tree-like) structure of the examined system;
 (ii) design alternative set for each leaf node
 (system component)
 of the system
 model;
 (iii) assessment of the design alternatives (DAs) on the basis of
 an ordinal scale (e.g., \cite{lev98},\cite{lev06},\cite{lev09})
 or a special interval multiset based scale \cite{lev12a};
 and
 (iv) ordinal compatibility between the design alternatives for
 different system components (instead of binary compatibility that
 was used in morphological analysis).

 An example of morphological hierarchy
 for a compressed plan of car repair
 from  \cite{lev12shop} is the following
 (Fig. 27, ordinal estimates of DAs are shown in parentheses; Table 5, Table 6):

 ~~

 {\bf 0.} Plan ~ \( S = A*B*C \)

 {\bf 1.} Payment ~ \( A = X*F \)

 {\it 1.1.} payment scheme ~\( X \):
  100 \% payment \(X_{0}\),
  prepayment of 50...80 percent for parts \(X_{1}\);
  bank loan \(X_{2}\);

 {\it 1.2.} version ~\( F \):
  cash \(F_{1}\), credit card \(F_{2}\),
  bank transfer \(F_{3}\).

 {\bf 2.} Body ~ \( B = R*Z*M \):

 {\it 2.1.} frame ~\( R \):
  None \(R_{0}\), technical diagnostics \(R_{1}\),
 follow-up assembly \(R_{2}\);

 {\it 2.2.} hardware ~\( Z \):
  None \(Z_{0}\),
  replacement of defect parts \(Z_{1}\),
  repair of body-defects \(Z_{2}\),
  fitting \(Z_{3}\),
  \(Z_{4}=Z_{1} \& Z_{2} \), \(Z_{5} = Z_{1} \& Z_{3}\),
  \(Z_{6} = Z_{2} \& Z_{3}\), \(Z_{7} = Z_{1} \& Z_{2} \& Z_{3}\);

  {\it 2.3.} finishing ~\( M = U * V\):

  {\it 2.3.1.} painting ~\( U \):
   None \(U_{0}\), partial painting \(U_{1}\), painting \(U_{2}\);

  {\it 2.3.2.} appearance restoration ~\( V \):
   None \(V_{0}\), Yes \(V_{1}\).

 {\bf 3.} Electric \& electronic subsystem ~ \( C = P * Q \):

 {\it 3.1.} Computer \& navigation subsystem ~ \( P = K * G \)

 {\it 3.1.1.} Computer ~\( K \):
  None \(K_{0}\), upgrade \(K_{1}\),
  additional or new computer \(K_{2}\);

 {\it 3.1.2.} system GPS ~\( G \):
  None \(G_{0}\), GPS system \(G_{1}\);

 {\it 3.2.} wiring \& lighting  ~ \( Q = O * L \)

 {\it 3.2.1.} wiring ~\( O \):
  None \(O_{0}\), repair \(O_{1}\);

 {\it 3.2.2.} lighting  ~\( L \):
  None \(L_{0}\), partial replacement \(L_{1}\),
  replacement \(L_{2}\).

%~~

\begin{center}
\begin{picture}(72,92)
\put(05,0){\makebox(0,0)[bl] {Fig. 27. Structure of repair plan
\cite{lev12shop} }}

\put(26,87){\makebox(0,8)[bl]{\(S=A\star B\star C\)}}

\put(23,88){\circle*{3}}

\put(23,84){\line(0,1){4}}

\put(0,84){\line(1,0){55}}

%==============================
%-----------------------------B

 \put(23,80){\line(0,1){4}}

\put(24.5,79){\makebox(0,8)[bl]{\(B=R\star Z\star M\)}}

\put(23,79){\circle*{2.3}}

\put(23,44){\line(0,1){34}}

%-----------------------------A

\put(0,80){\line(0,1){4}}

\put(1.5,79){\makebox(0,8)[bl]{\(A=X\star F\)}}

\put(0,79){\circle*{2.3}}

\put(0,74){\line(0,1){4}}

%---------------- X

\put(0,74){\line(1,0){10}}

\put(0,70){\line(0,1){4}}

\put(0,69){\circle*{1.5}}

\put(2,69){\makebox(0,8)[bl]{\(X\)}}

\put(0,64){\makebox(0,8)[bl]{\(X_{0}(2)\)}}
\put(0,60){\makebox(0,8)[bl]{\(X_{1}(1)\)}}
\put(0,56){\makebox(0,8)[bl]{\(X_{2}(3)\)}}

%---------------- F

\put(10,70){\line(0,1){4}}

\put(10,69){\circle*{1.5}}

\put(12,69){\makebox(0,8)[bl]{\(F\)}}

\put(10,64){\makebox(0,8)[bl]{\(F_{1}(2)\)}}
\put(10,60){\makebox(0,8)[bl]{\(F_{2}(1)\)}}
\put(10,56){\makebox(0,8)[bl]{\(F_{3}(3)\)}}

%===============================
%---------------- R

\put(0,44){\line(1,0){25}}

\put(0,40){\line(0,1){4}}

\put(0,39){\circle*{1.5}}

\put(2,39){\makebox(0,8)[bl]{\(R\)}}

\put(0,34){\makebox(0,8)[bl]{\(R_{0}(2)\)}}
\put(0,30){\makebox(0,8)[bl]{\(R_{1}(1)\)}}
\put(0,26){\makebox(0,8)[bl]{\(R_{2}(3)\)}}

%---------------- Z

\put(10,40){\line(0,1){4}}

\put(10,39){\circle*{1.5}}

\put(12,39){\makebox(0,8)[bl]{\(Z\)}}

\put(10,34){\makebox(0,8)[bl]{\(Z_{0}(2)\)}}
\put(10,30){\makebox(0,8)[bl]{\(Z_{1}(1)\)}}
\put(10,26){\makebox(0,8)[bl]{\(Z_{2}(3)\)}}
\put(10,22){\makebox(0,8)[bl]{\(Z_{3}(2)\)}}
\put(07,18){\makebox(0,8)[bl]{\(Z_{4}=Z_{1}\&Z_{2}(2)\)}}
\put(07,14){\makebox(0,8)[bl]{\(Z_{5}=Z_{2}\&Z_{3}(2)\)}}
\put(07,10){\makebox(0,8)[bl]{\(Z_{6}=Z_{1}\&Z_{3}(1)\)}}
\put(07,6){\makebox(0,8)[bl]{\(Z_{7}=Z_{1}\&Z_{2}\&Z_{3}(3)\)}}

%-------------------------M
\put(25,40){\line(0,1){4}}

\put(25,39){\circle*{2}}

\put(27,38){\makebox(0,8)[bl]{\(M=U\star V\)}}

%-------------------------U
\put(25,31){\line(0,1){7}}

\put(25,31){\line(1,0){20}}

\put(35,27){\line(0,1){4}}

\put(35,26){\circle*{1.5}}

\put(37,26){\makebox(0,8)[bl]{\(U\)}}

\put(33,21){\makebox(0,8)[bl]{\(U_{0}(2)\)}}
\put(33,17){\makebox(0,8)[bl]{\(U_{1}(1)\)}}
\put(33,13){\makebox(0,8)[bl]{\(U_{2}(3)\)}}

%---------------- V

\put(45,27){\line(0,1){4}}

\put(45,26){\circle*{1.5}}

\put(47,26){\makebox(0,8)[bl]{\(V\)}}

\put(43,21){\makebox(0,8)[bl]{\(V_{0}(2)\)}}
\put(43,17){\makebox(0,8)[bl]{\(V_{1}(1)\)}}

%======================================================

\put(55,80){\line(0,1){4}}

\put(56.5,78){\makebox(0,8)[bl]{\(C=P\star Q\)}}

\put(55,79){\circle*{2.3}}

\put(55,74){\line(0,1){4}}

\put(25,74){\line(1,0){30}}

%============================================

\put(25,70){\line(0,1){4}}

\put(25,69){\circle*{2}}

\put(26.5,69){\makebox(0,8)[bl]{\(P=K\star G\)}}

%---------------------------------- K
\put(25,64){\line(0,1){4}}

\put(25,64){\line(1,0){10}}

\put(25,60){\line(0,1){4}}

\put(25,59){\circle*{1.5}}

\put(27,59){\makebox(0,8)[bl]{\(K\)}}

\put(25,54){\makebox(0,8)[bl]{\(K_{0}(2)\)}}
\put(25,50){\makebox(0,8)[bl]{\(K_{1}(1)\)}}
\put(25,46){\makebox(0,8)[bl]{\(K_{2}(3)\)}}

%---------------- G

\put(35,60){\line(0,1){4}}

\put(35,59){\circle*{1.5}}

\put(37,59){\makebox(0,8)[bl]{\(G\)}}

\put(35,54){\makebox(0,8)[bl]{\(G_{0}(2)\)}}
\put(35,50){\makebox(0,8)[bl]{\(G_{1}(1)\)}}

%======================================================

\put(55,70){\line(0,1){4}}

\put(55,69){\circle*{2}}

\put(57,68.5){\makebox(0,8)[bl]{\(Q=O\star L\)}}

%---------------------------------- O
\put(55,64){\line(0,1){4}}

\put(55,64){\line(1,0){10}}

\put(55,60){\line(0,1){4}}

\put(55,59){\circle*{1.5}}

\put(57,59){\makebox(0,8)[bl]{\(O\)}}

\put(54,54){\makebox(0,8)[bl]{\(O_{0}(2)\)}}
\put(54,50){\makebox(0,8)[bl]{\(O_{1}(1)\)}}

%---------------- L

\put(65,60){\line(0,1){4}}

\put(65,59){\circle*{1.5}}

\put(67,59){\makebox(0,8)[bl]{\(L\)}}

\put(65,54){\makebox(0,8)[bl]{\(L_{0}(2)\)}}
\put(65,50){\makebox(0,8)[bl]{\(L_{1}(1)\)}}
\put(65,46){\makebox(0,8)[bl]{\(L_{2}(1)\)}}

\end{picture}
\end{center}

\begin{center}
\begin{picture}(37,53)

\put(09.5,49){\makebox(0,0)[bl]{Table 5. Compatibility
\cite{lev12shop}}}

\put(00,0){\line(1,0){37}} \put(00,42){\line(1,0){37}}
\put(00,48){\line(1,0){37}}

\put(00,0){\line(0,1){48}} \put(07,0){\line(0,1){48}}
\put(37,0){\line(0,1){48}}

\put(01,38){\makebox(0,0)[bl]{\(Z_{0}\)}}
\put(01,34){\makebox(0,0)[bl]{\(Z_{1}\)}}
\put(01,30){\makebox(0,0)[bl]{\(Z_{2}\)}}
\put(01,26){\makebox(0,0)[bl]{\(Z_{3}\)}}
\put(01,22){\makebox(0,0)[bl]{\(Z_{4}\)}}
\put(01,18){\makebox(0,0)[bl]{\(Z_{5}\)}}
\put(01,14){\makebox(0,0)[bl]{\(Z_{6}\)}}
\put(01,10){\makebox(0,0)[bl]{\(Z_{7}\)}}
\put(01,06){\makebox(0,0)[bl]{\(M_{1}\)}}
\put(01,02){\makebox(0,0)[bl]{\(M_{2}\)}}

\put(13,42){\line(0,1){6}} \put(19,42){\line(0,1){6}}
\put(25,42){\line(0,1){6}} \put(31,42){\line(0,1){6}}

\put(07.5,44){\makebox(0,0)[bl]{\(M_{1}\)}}
\put(13.5,44){\makebox(0,0)[bl]{\(M_{2}\)}}
\put(19.5,44){\makebox(0,0)[bl]{\(R_{0}\)}}
\put(25.5,44){\makebox(0,0)[bl]{\(R_{1}\)}}
\put(31.5,44){\makebox(0,0)[bl]{\(R_{2}\)}}

\put(09,38){\makebox(0,0)[bl]{\(2\)}}
\put(15,38){\makebox(0,0)[bl]{\(3\)}}
\put(21,38){\makebox(0,0)[bl]{\(3\)}}
\put(27,38){\makebox(0,0)[bl]{\(3\)}}
\put(33,38){\makebox(0,0)[bl]{\(0\)}}

\put(09,34){\makebox(0,0)[bl]{\(3\)}}
\put(15,34){\makebox(0,0)[bl]{\(2\)}}
\put(21,34){\makebox(0,0)[bl]{\(2\)}}
\put(27,34){\makebox(0,0)[bl]{\(3\)}}
\put(33,34){\makebox(0,0)[bl]{\(3\)}}

\put(09,30){\makebox(0,0)[bl]{\(3\)}}
\put(15,30){\makebox(0,0)[bl]{\(2\)}}
\put(21,30){\makebox(0,0)[bl]{\(0\)}}
\put(27,30){\makebox(0,0)[bl]{\(3\)}}
\put(33,30){\makebox(0,0)[bl]{\(3\)}}

\put(09,26){\makebox(0,0)[bl]{\(3\)}}
\put(15,26){\makebox(0,0)[bl]{\(2\)}}
\put(21,26){\makebox(0,0)[bl]{\(0\)}}
\put(27,26){\makebox(0,0)[bl]{\(2\)}}
\put(33,26){\makebox(0,0)[bl]{\(3\)}}

\put(09,22){\makebox(0,0)[bl]{\(3\)}}
\put(15,22){\makebox(0,0)[bl]{\(2\)}}
\put(21,22){\makebox(0,0)[bl]{\(2\)}}
\put(27,22){\makebox(0,0)[bl]{\(3\)}}
\put(33,22){\makebox(0,0)[bl]{\(3\)}}

\put(09,18){\makebox(0,0)[bl]{\(3\)}}
\put(15,18){\makebox(0,0)[bl]{\(2\)}}
\put(21,18){\makebox(0,0)[bl]{\(0\)}}
\put(27,18){\makebox(0,0)[bl]{\(3\)}}
\put(33,18){\makebox(0,0)[bl]{\(3\)}}

\put(09,14){\makebox(0,0)[bl]{\(3\)}}
\put(15,14){\makebox(0,0)[bl]{\(2\)}}
\put(21,14){\makebox(0,0)[bl]{\(2\)}}
\put(27,14){\makebox(0,0)[bl]{\(3\)}}
\put(33,14){\makebox(0,0)[bl]{\(3\)}}

\put(09,10){\makebox(0,0)[bl]{\(3\)}}
\put(15,10){\makebox(0,0)[bl]{\(2\)}}
\put(21,10){\makebox(0,0)[bl]{\(2\)}}
\put(27,10){\makebox(0,0)[bl]{\(3\)}}
\put(33,10){\makebox(0,0)[bl]{\(3\)}}

%\put(09,06){\makebox(0,0)[bl]{\(5\)}}
%\put(15,06){\makebox(0,0)[bl]{\(5\)}}
\put(21,06){\makebox(0,0)[bl]{\(0\)}}
\put(27,06){\makebox(0,0)[bl]{\(3\)}}
\put(33,06){\makebox(0,0)[bl]{\(3\)}}

%\put(09,02){\makebox(0,0)[bl]{\(5\)}}
%\put(15,02){\makebox(0,0)[bl]{\(5\)}}
\put(21,02){\makebox(0,0)[bl]{\(3\)}}
\put(27,02){\makebox(0,0)[bl]{\(2\)}}
\put(33,02){\makebox(0,0)[bl]{\(2\)}}

\end{picture}
%\end{center}
%
%\begin{center}
%\begin{picture}(25,27)
\begin{picture}(35,27)

\put(00,0){\line(1,0){25}} \put(00,14){\line(1,0){25}}
\put(00,20){\line(1,0){25}}

\put(00,0){\line(0,1){20}} \put(07,0){\line(0,1){20}}
\put(25,0){\line(0,1){20}}

\put(01,10){\makebox(0,0)[bl]{\(X_{1}\)}}
\put(01,06){\makebox(0,0)[bl]{\(X_{2}\)}}
\put(01,02){\makebox(0,0)[bl]{\(X_{3}\)}}

\put(13,14){\line(0,1){6}} \put(19,14){\line(0,1){6}}

\put(07.5,16){\makebox(0,0)[bl]{\(F_{1}\)}}
\put(13.5,16){\makebox(0,0)[bl]{\(F_{2}\)}}
\put(19.5,16){\makebox(0,0)[bl]{\(F_{3}\)}}

\put(09,10){\makebox(0,0)[bl]{\(3\)}}
\put(15,10){\makebox(0,0)[bl]{\(3\)}}
\put(21,10){\makebox(0,0)[bl]{\(3\)}}

\put(09,06){\makebox(0,0)[bl]{\(3\)}}
\put(15,06){\makebox(0,0)[bl]{\(3\)}}
\put(21,06){\makebox(0,0)[bl]{\(3\)}}

\put(09,02){\makebox(0,0)[bl]{\(0\)}}
\put(15,02){\makebox(0,0)[bl]{\(3\)}}
\put(21,02){\makebox(0,0)[bl]{\(2\)}}

%======================================

\put(00,32){\line(1,0){25}} \put(00,42){\line(1,0){25}}
\put(00,48){\line(1,0){25}}

\put(00,32){\line(0,1){16}} \put(07,32){\line(0,1){16}}
\put(25,32){\line(0,1){16}}

\put(01,38){\makebox(0,0)[bl]{\(O_{0}\)}}
\put(01,34){\makebox(0,0)[bl]{\(O_{1}\)}}

\put(13,42){\line(0,1){6}} \put(19,42){\line(0,1){6}}

\put(07.5,44){\makebox(0,0)[bl]{\(L_{0}\)}}
\put(13.5,44){\makebox(0,0)[bl]{\(L_{1}\)}}
\put(19.5,44){\makebox(0,0)[bl]{\(L_{2}\)}}

\put(09,38){\makebox(0,0)[bl]{\(3\)}}
\put(15,38){\makebox(0,0)[bl]{\(2\)}}
\put(21,38){\makebox(0,0)[bl]{\(2\)}}

\put(09,34){\makebox(0,0)[bl]{\(1\)}}
\put(15,34){\makebox(0,0)[bl]{\(3\)}}
\put(21,34){\makebox(0,0)[bl]{\(3\)}}

%======================================

\end{picture}
%\end{center}
%
%\begin{center}
\begin{picture}(24,25)

\put(0.5,21){\makebox(0,0)[bl]{Table 6. Compatibility
\cite{lev12shop}}}

\put(00,0){\line(1,0){19}} \put(00,14){\line(1,0){19}}
\put(00,20){\line(1,0){19}}

\put(00,0){\line(0,1){20}} \put(07,0){\line(0,1){20}}
\put(19,0){\line(0,1){20}}

\put(01,10){\makebox(0,0)[bl]{\(K_{0}\)}}
\put(01,06){\makebox(0,0)[bl]{\(K_{1}\)}}
\put(01,02){\makebox(0,0)[bl]{\(K_{2}\)}}

\put(13,14){\line(0,1){6}}

\put(07.5,16){\makebox(0,0)[bl]{\(G_{0}\)}}
\put(13.5,16){\makebox(0,0)[bl]{\(G_{1}\)}}

\put(09,10){\makebox(0,0)[bl]{\(3\)}}
\put(15,10){\makebox(0,0)[bl]{\(0\)}}

\put(09,06){\makebox(0,0)[bl]{\(2\)}}
\put(15,06){\makebox(0,0)[bl]{\(3\)}}

\put(09,02){\makebox(0,0)[bl]{\(1\)}}
\put(15,02){\makebox(0,0)[bl]{\(2\)}}

\end{picture}
%\end{center}
%
%\begin{center}
\begin{picture}(19,25)

\put(00,0){\line(1,0){19}} \put(00,14){\line(1,0){19}}
\put(00,20){\line(1,0){19}}

\put(00,0){\line(0,1){20}} \put(07,0){\line(0,1){20}}
\put(19,0){\line(0,1){20}}

\put(01,10){\makebox(0,0)[bl]{\(U_{0}\)}}
\put(01,06){\makebox(0,0)[bl]{\(U_{1}\)}}
\put(01,02){\makebox(0,0)[bl]{\(U_{2}\)}}

\put(13,14){\line(0,1){6}}

\put(07.5,16){\makebox(0,0)[bl]{\(V_{0}\)}}
\put(13.5,16){\makebox(0,0)[bl]{\(V_{1}\)}}

\put(09,10){\makebox(0,0)[bl]{\(3\)}}
\put(15,10){\makebox(0,0)[bl]{\(0\)}}

\put(09,06){\makebox(0,0)[bl]{\(0\)}}
\put(15,06){\makebox(0,0)[bl]{\(2\)}}

\put(09,02){\makebox(0,0)[bl]{\(0\)}}
\put(15,02){\makebox(0,0)[bl]{\(3\)}}

\end{picture}
\end{center}

 The evident design scheme for morphological hierarchy is the following:

~~

 {\it Stage 1.} Design of system hierarchy (one of the methods
 above).

 {\it Stage 2.} Generation of DAs for each leaf node of the system
 model (expert judgement and/or usage of information bases).

  {\it Stage 3.} Generation of scales for assessment of DAs
  and assessment of DAs (usually multicriteria description of DAs
  is used at a preliminary phase).

  {\it Stage 4.} generation  of an ordinal scale for compatibility
  among DAs and assessment of the compatibility.

\section{SOME APPROACHES TO MODIFICATION}

\subsection{Modification of
 Tree via Condensing of Weighted  Edges}

 This section
  \footnote{
 From (with amendments):\\
  (i) M.Sh. Levin, An extremal problem of organization of
 data. {\it Eng. Cybern.}, 19(5), 1981, 87-95.\\
 (ii) M.Sh. Levin, {\it Combinatorial Engineering of Decomposable
 Systems}, Springer, 1998, Chapter 2.
  }
%%%%%%%%%%%%%%%%%
  describes transformation of a tree
 (with weights of vertices and weights of edge/arcs)
 via integration (condensing)
 of some neighbor vertices
 while taking into account
 a constraint for a total weight of the maximum tree tail
 (i.e., length from root to a leaf vertex).
 The problem was firstly formulated
  for designing an overlay structure of
  a  modular software system in \cite{lev81}.
 The integration of software modules requires additional memory, but
 allows to decrease a time (i.e., frequency) of loading some
 corresponding modules.
 Other applications of the problem can be examined as well, e.g.,
  hierarchical structure of data,
  call problem,
  hierarchical information structure of Web-sites.
  This problem  is illustrated in
 Fig. 28 and Fig. 29 by an example for designing the over-lay structure
 on the basis of module integration, when different software or data modules
 can apply the same parts of RAM.
 A  new kind of FPTAS for the above-mentioned combinatorial optimization
 problem (a generalization of multiple choice problem over a
 tree-like structure and special constraints)
 was suggested in \cite{lev81}.

\begin{center}
\begin{picture}(95,41)

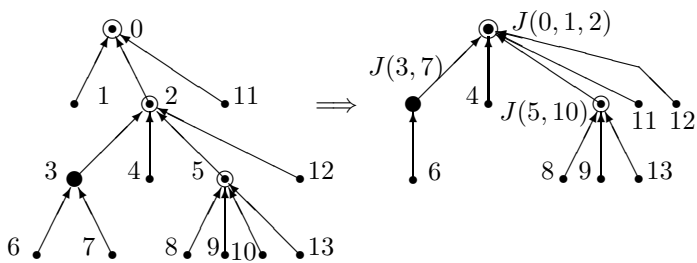
\put(00,00){\makebox(0,0)[bl] {Fig. 28. Integration of software
 modules (over-lay structure)}}

\put(5,07){\circle*{1}} \put(1,07){\makebox(0,0)[bl]{\(6\)}}

\put(15,07){\circle*{1}} \put(11,07){\makebox(0,0)[bl]{\(7\)}}

\put(05,07){\vector(1,2){4.5}} \put(15,07){\vector(-1,2){4.5}}

\put(25,07){\circle*{1}} \put(22,07){\makebox(0,0)[bl]{\(8\)}}

\put(35,07){\circle*{1}}
\put(30.7,06.3){\makebox(0,0)[bl]{\(10\)}}

\put(40,07){\circle*{1}} \put(41,07){\makebox(0,0)[bl]{\(13\)}}
\put(40,07){\vector(-1,1){9.5}}

\put(25,07){\vector(1,2){4.5}} \put(35,07){\vector(-1,2){4.5}}

\put(30,07){\circle*{1}} \put(27.6,07){\makebox(0,0)[bl]{\(9\)}}
\put(30,07){\vector(0,1){09}}

\put(10,17){\circle*{2}} \put(6,17){\makebox(0,0)[bl]{\(3\)}}
\put(20,17){\circle*{1}} \put(17,17){\makebox(0,0)[bl]{\(4\)}}
\put(30,17){\circle*{1}} \put(30,17){\circle{2}}
\put(25.5,17){\makebox(0,0)[bl]{\(5\)}}

\put(40,17){\circle*{1}} \put(41,17){\makebox(0,0)[bl]{\(12\)}}
\put(40,17){\vector(-2,1){19}}

\put(10,17){\vector(1,1){09}} \put(30,17){\vector(-1,1){09}}
\put(20,17){\vector(0,1){09}}

\put(20,27){\circle*{1}} \put(20,27){\circle{2.1}}
\put(22,27){\makebox(0,0)[bl]{\(2\)}}

\put(30,27){\circle*{1}} \put(31,27){\makebox(0,0)[bl]{\(11\)}}
\put(30,27){\vector(-3,2){14}}

\put(10,27){\circle*{1}} \put(13,27){\makebox(0,0)[bl]{\(1\)}}
\put(10,27){\vector(1,2){4.5}} \put(20,27){\vector(-1,2){4.5}}

\put(15,37){\circle*{1.2}} \put(15,37){\circle{2.3}}
\put(17.4,36){\makebox(0,0)[bl]{\(0\)}}

%========================================== Right part

\put(65,37){\circle*{1.3}} \put(65,37){\circle{2.5}}
\put(68,36){\makebox(0,0)[bl]{\(J(0,1,2)\)}}

\put(55,27){\circle*{2}}
\put(49,30){\makebox(0,0)[bl]{\(J(3,7)\)}}

\put(65,27){\circle*{1}} \put(62,27){\makebox(0,0)[bl]{\(4\)}}

\put(80,27){\circle*{1}} \put(80,27){\circle{2}}
\put(80,27){\vector(-3,2){14}}
\put(66.5,24){\makebox(0,0)[bl]{\(J(5,10)\)}}

\put(85,27){\circle*{1}} \put(84,23.7){\makebox(0,0)[bl]{\(11\)}}
\put(85,27){\vector(-2,1){19}}

\put(90,27){\circle*{1}} \put(89,23.7){\makebox(0,0)[bl]{\(12\)}}
\put(90,27){\line(-1,1){5}} \put(85,32){\vector(-4,1){19}}

\put(55,27){\vector(1,1){09}} \put(65,27){\vector(0,1){09}}

\put(55,17){\circle*{1}} \put(57,17){\makebox(0,0)[bl]{\(6\)}}
\put(55,17){\vector(0,1){09}}

\put(75,17){\circle*{1}} \put(72,17){\makebox(0,0)[bl]{\(8\)}}
\put(75,17){\vector(1,2){4.5}}

\put(80,17){\circle*{1}} \put(77,17){\makebox(0,0)[bl]{\(9\)}}
\put(80,17){\vector(0,1){9}}

\put(85,17){\circle*{1}} \put(86,17){\makebox(0,0)[bl]{\(13\)}}
\put(85,17){\vector(-1,2){4.5}}

%%%%%%%%%%%%%%%%%%%%%%%%%%%%%%%%%%%%

\put(42,26){\makebox(0,0)[bl]{\(\Longrightarrow\)}}

\end{picture}
\end{center}

\begin{center}
\begin{picture}(97,63)
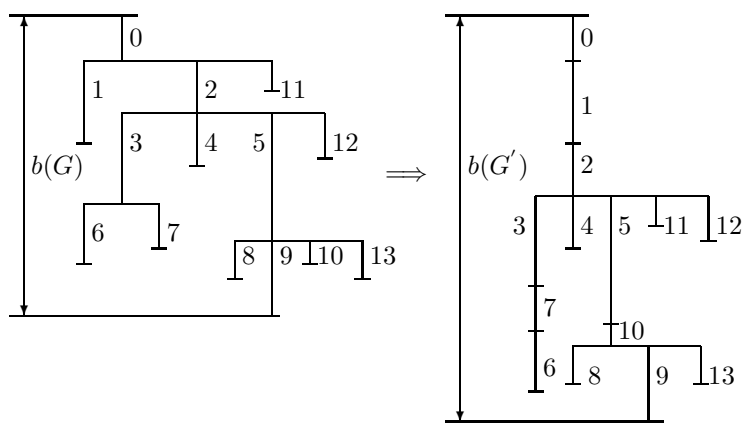
\put(17,00){\makebox(0,0)[bl] {Fig. 29. Usage of memory (RAM)}}

%%%%%%%%%%%%%%%% Right part

\put(58,07){\line(1,0){29}}

\put(085,07){\line(0,1){10}}
\put(086,12){\makebox(0,0)[bl]{\(9\)}}

\put(92,17){\line(0,-1){5}} \put(93,12){\makebox(0,0)[bl]{\(13\)}}
\put(085,17){\line(1,0){07}} \put(91,12){\line(1,0){02}}

\put(075,12){\line(0,1){5}} \put(077,12){\makebox(0,0)[bl]{\(8\)}}
\put(074,12){\line(1,0){2}}

\put(075,17){\line(1,0){10}}

\put(080,17){\line(0,1){3}}
\put(081,18){\makebox(0,0)[bl]{\(10\)}}
\put(079,20){\line(1,0){2}}

\put(080,20){\line(0,1){17}}
\put(081,32){\makebox(0,0)[bl]{\(5\)}}

\put(86,37){\line(0,-1){4}} \put(87,32){\makebox(0,0)[bl]{\(11\)}}
\put(80,37){\line(1,0){06}} \put(85,33){\line(1,0){02}}

\put(93,37){\line(0,-1){6}} \put(94,32){\makebox(0,0)[bl]{\(12\)}}
\put(86,37){\line(1,0){07}} \put(92,31){\line(1,0){02}}

\put(75,30){\line(0,1){7}} \put(76,32){\makebox(0,0)[bl]{\(4\)}}
\put(74,30){\line(1,0){2}}

\put(70,37){\line(1,0){10}}

\put(070,25){\line(0,1){12}}
\put(067,32){\makebox(0,0)[bl]{\(3\)}} \put(069,25){\line(1,0){2}}

\put(070,19){\line(0,1){6}} \put(071,21){\makebox(0,0)[bl]{\(7\)}}
\put(069,19){\line(1,0){2}}

\put(070,11){\line(0,1){8}} \put(071,13){\makebox(0,0)[bl]{\(6\)}}
\put(069,11){\line(1,0){2}}

\put(075,37){\line(0,1){7}} \put(076,40){\makebox(0,0)[bl]{\(2\)}}
\put(074,37){\line(1,0){2}}

\put(075,44){\line(0,1){11}}
\put(076,48){\makebox(0,0)[bl]{\(1\)}} \put(074,44){\line(1,0){2}}

\put(075,55){\line(0,1){6}} \put(076,57){\makebox(0,0)[bl]{\(0\)}}
\put(074,55){\line(1,0){2}}

\put(58,61){\line(1,0){19}}

\put(60,12){\vector(0,-1){5}} \put(60,12){\vector(0,1){49}}
\put(61,39){\makebox(0,0)[bl]{\(b(G^{'})\)}}

%%%%%%%%%%%%%%%%%%%%%

\put(50,39){\makebox(0,0)[bl]{\(\Longrightarrow\)}}

%%%%%%%%%%%%%%%%%%%%%%%%%%%%Left p[art

\put(0,61){\line(1,0){17}} \put(0,21){\line(1,0){36}}

\put(2,32){\vector(0,-1){11}} \put(2,32){\vector(0,1){29}}

\put(3,39){\makebox(0,0)[bl]{\(b(G)\)}}

\put(015,55){\line(0,1){6}} \put(16,57){\makebox(0,0)[bl]{\(0\)}}
\put(010,55){\line(1,0){15}}

\put(010,44){\line(0,1){11}}
\put(011,50){\makebox(0,0)[bl]{\(1\)}} \put(09,44){\line(1,0){2}}

\put(025,55){\line(0,-1){7}}
\put(026,50){\makebox(0,0)[bl]{\(2\)}}
\put(015,48){\line(1,0){20}}

\put(035,55){\line(0,-1){4}}
\put(036,50){\makebox(0,0)[bl]{\(11\)}}
\put(025,55){\line(1,0){10}} \put(034,51){\line(1,0){02}}

\put(035,48){\line(0,-1){17}}
\put(032.4,43){\makebox(0,0)[bl]{\(5\)}}
\put(030,31){\line(1,0){10}}

\put(042,48){\line(0,-1){6}}
\put(043,43){\makebox(0,0)[bl]{\(12\)}}
\put(035,48){\line(1,0){07}} \put(041,42){\line(1,0){02}}

\put(025,48){\line(0,-1){7}}
\put(026,43){\makebox(0,0)[bl]{\(4\)}} \put(024,41){\line(1,0){2}}

\put(015,48){\line(0,-1){12}}
\put(016,43){\makebox(0,0)[bl]{\(3\)}}
\put(010,36){\line(1,0){10}}

\put(020,36){\line(0,-1){6}}
\put(021,31){\makebox(0,0)[bl]{\(7\)}} \put(019,30){\line(1,0){2}}

\put(010,36){\line(0,-1){8}}
\put(011,31){\makebox(0,0)[bl]{\(6\)}} \put(09,28){\line(1,0){2}}

\put(035,31){\line(0,-1){10}}
\put(036,28){\makebox(0,0)[bl]{\(9\)}} \put(034,21){\line(1,0){2}}

\put(030,31){\line(0,-1){5}}
\put(031,28){\makebox(0,0)[bl]{\(8\)}} \put(029,26){\line(1,0){2}}

\put(040,31){\line(0,-1){3}}
\put(041,28){\makebox(0,0)[bl]{\(10\)}}
\put(039,28){\line(1,0){2}}

\put(047,31){\line(0,-1){5}}
\put(048,28){\makebox(0,0)[bl]{\(13\)}}
\put(040,31){\line(1,0){07}} \put(046,26){\line(1,0){02}}

\end{picture}
\end{center}

 Let
 \(G=(A,\Upsilon)\) be an oriented tree, where
 \(A\) is a set of vertices (software or data modules).
 \(\Upsilon\) is a multi-valued mapping of \(A\) into \(A\).
 Arcs of \(G\)  are oriented from the root \(a_{o} \in A\)
 to leaf vertices. Each vertex \(a\in  A\) has a positive weight
 (required volume of RAM) \(\beta(a)>0\).
 Each arc \((a^{'},a^{''})\)
 (\(a^{'},a^{''}\in A\) ¨ \(a^{''}\in \Upsilon a^{'}\))
 has a weight (i.e., an initial frequency of loading into RAM)
 \(w(a^{'},a^{''})>0\).
 This arc weight corresponds to the frequency of calling (and loading)
 from module \(a^{'}\) to module \(a^{''}\).

 Let \( \pi  (a^{1},a^{l})=  <a^{1},...,a^{i},...,a^{l}>  \)
 be a path
 (\(a^{j+1}\in\Upsilon a^{j}, j=1,...,l-1\)).
 We propose for each path a weight
 \(\lambda (\pi(a^{1},a^{l})= \sum_{i=1}^{l}\lambda(a^{i}) \).
 Denote by a weight of graph \(G\) the value
 \[\lambda(G)=max_{a^{''}\in A^{o}}\{\lambda(\pi(a_{o},a^{''}))\},\]
 where
 \( A^{o}=\{a\in A\mid \Upsilon a=\O \}\) is a set of leaf vertices.
 Let
 \(G_{a}=(A_{a},\Upsilon)\) is a subtree with  root \(a\in A\),
 and \(A_{a}\) contains vertex \(a\) and all other vertices,
 which can be reached from \(a\).
 Graph \((A_{a}\setminus  a,\Upsilon  )\) is called
 {\em tail} of vertex \(a\), and  value
  \( \lambda^{-}(a)= \lambda(G_{a})- \lambda(a) \) is called a
 {\em tail} weight of vertex \(a\).
 Clearly that
 \[ \lambda(a)=max_{a^{'}\in \Upsilon a} \{\lambda(G_{a^{'}})\}. \]
 We examine weight \(w(a)\) and binary variable \(x(a)\)
 \(\forall a\in A\setminus a_{o}\)
 (\(1\) corresponds to a situation when the arc, directed to
 \(a\), is condensed).
 Now let us define a transformation of graph
 \(G\) on the basis of integrating the vertices
 \(a^{'}\) and \(a^{''}\) as follows:

(a) vertex \(a^{'}\) is changed into \(J(a^{'},a^{''})\)
 with the following properties:

 \(\lambda (J(a^{'},a^{''}))=\lambda(a^{'})+\lambda(a^{''})\) ¨
 \(\Upsilon J(a^{'},a^{''})=(\Upsilon a^{'}
 \bigcup \Upsilon a^{''})\setminus a^{''}\);

 (b) vertex \(a^{''}\) and arcs, which are oriented from the vertex,
 are deleted.

 For graph \(G\), we propose a binary vector
 \(\kappa(a)\) that involves all
 \(x(a) ~\forall a\in A\setminus a_{o}\).
 Thus, we examine the weights of vertex
  \(a\) and its {\em tail}
  as functions of vector \(\kappa\):
 \(\lambda (a,\kappa ), \lambda^{-}(a,\kappa )\).
 Now let us consider a problem (kind 1):
 \[\max ~W(\kappa )=\sum_{a\in A\setminus a_{o}} x(a)w(a)
 ~~~ s.t. ~\lambda (a_{o},\kappa)+\lambda^{-}(a_{o},\kappa) \leq b , \]
 where \(b\) is a positive constant
 (i.e., a volume of accessible RAM).
 In general,
 this problem formulation corresponds to the example in Fig. 28 and Fig. 29.

 In addition, we examine analogical problem (kind 2)
 with other constraints as follows:
 \[ \lambda (a_{o},\kappa) \leq b^{-}, ~ \lambda^{-}(a_{o},\kappa) \leq
 b^{+}, ~
 b^{-}+b^{+}=b.\]
 Note, illustrations
 of the class of considered combinatorial problems are presented
 in Fig. 30 (basic knapsack problem and multiple choice problem)
 and in Fig. 31, correspondence of problem to illustration is pointed out in Table 7.

 Now consider some simple cases of the problems (kind 1 and kind 2).
 Let \(\Upsilon a_{o}=\{a_{1},...,a_{i},...,a_{m}\}\)
 (and \(u_{i}\)) corresponds to an arc \((a_{o},a_{i})\)~
 (\(w(u_{i})=w_{i}\)).
%
%%%%%%%%%%%%%%%%%%%
 Then, corresponding problem
 (problem 1, an equivalent to knapsack problem, Fig. 30i) is:
 \[\max ~ \sum_{i=1}^{m} x_{i} w_{i}
 ~~~~ s.t. ~~ \lambda(a_{o})+\sum_{i=1}^{m}x_{i}\lambda(a_{i})\leq
 b, ~~
 x_{i} \in \{0,1\}.\]
 The objective function in other simple cases
  (1.1 - Fig. 31a, 1.2 - Fig. 31b, 1.3 - Fig. 31c, 1.4 - Fig. 31d),
 which are based on knapsack problem (problem 1, Fig. 30i)
 is analogical,
 and only constraints will be presented for them.

%%%%%%%%%%%%%%%%%
 Problem 1.1 (Fig. 31a) has the following constraint:
 \[\lambda(a_{o})+ \sum_{i=1}^{m}x_{i} \lambda(a_{i}) +
 \max_{1 \leq i \leq m}((1-x_{i}) \lambda(a_{i})) \leq b.\]
 This problem corresponds to a  ``kernel'' load in many software packages.

%
%%%%%%%%%%%%%%
%
 Problem 1.2 (Fig. 31b) of kind 2 is the following:
 \[\lambda(a_{o})+ \sum_{i=1}^{m}x_{i} \lambda(a_{i}) \leq b^{-},
 ~~ \max_{1 \leq i \leq m}((1-x_{i} \lambda(a_{i})) \leq b^{+}.\]
 Problem 1.3 (Fig. 31c) is:
 \[\lambda(a_{o}+ \sum_{i=1}^{m} x_{i}\lambda^{-}(a_{i})+
 \max_{1\leq i\leq m}((1-x_{i}\lambda^{-}(a_{i})) +\lambda^{+}(a_{i}))\leq b.\]
 It is reasonable to point the following properties of this problem:

 (a) \(a_{i}\)  (\(\forall a_{i}\in \Upsilon a_{o}\))
 has the weight \(\lambda^{-}(a_{i})\);

 (b) \(a_{i}\)  (\(\forall a_{i} \in \Upsilon a_{o}\))
 has the only one son with a weight
 \(\lambda^{+}(a_{i})\), and the value is the
 {\em tail} weight; and

 (c) only condensing the following arcs \((a_{o},a_{i}) ~ (i=\overline{1,m})\)
 is admissible.

 As a result,
 a sequence of simple problems  based on
 knapsack problem can be examined:
  1 (basic knapsack problem, Fig. 30i),
  1.1 (analogue of knapsack problem, Fig. 31a),
   1.2 (Fig. 31b), 1.3 (Fig. 31c), 1.4 (Fig. 31d).

 In the same way,
  a sequence of auxiliary problems based on multiple choice
   problem can be considered:
    2 (basic multiple choice  problem, Fig. 30ii),
   2.1 (analogue of multiple choice problem, Fig. 31e),
   2.2 (Fig. 31f),
   2.3 (Fig. 31g),
   2.4 (Fig. 31h).

%%%%%%%%%%%
 In our case, multiple choice problem or problem 2 (Fig. 28ii)
 is the following:
 \[\max ~W(\{x_{ij}\})=\sum_{i=1}^{m} \sum_{j=1}^{q_{i}} w(a_{ij}) x_{ij}
 ~~~ s.t. ~\lambda(a_{o}) \sum_{i=1}^{m} \sum_{j=1}^{q_{i}}
 x_{ij}\lambda(a_{ij}) \leq b,
 ~ \sum_{j=1}^{q_{i}} x_{ij}=1, ~i= \overline{1,m}; ~x_{ij} \in \{0,1\}.\]

 Here, the following set of Boolean vectors in auxiliary problems is used:
 \[X=\{\kappa=(x^{1}_{ij};x^{2}_{ij})|x^{1}_{ij},x^{2}_{ij}\in \{0,1\};
 ~j=\overline{1,q_{i}};~ i=\overline{1,m}\}\]
 In addition, the following constraint has to be taken into account
 in all auxiliary problems:
 \[\sum_{j=1}^{q_{i}}x^{1}_{ij}=1, ~\forall i; ~x^{2}_{ij}\leq x^{1}_{ij},
 ~ \forall i,j.\]
 Also, the following modified objective function is used:
 \[W(X) = \sum_{i=1}^{m} \sum_{j=1}^{q_{i}} (x^{1}_{ij}w^{-}(a_{ij}) +
 x^{2}_{ij}w(a_{ij})).\]
%
%%%%%%%%%%%%%%%
%
 Now, for example, auxiliary problem 2.4 is considered
 that corresponds to kind 2 above (Fig. 31h):
 \[\lambda(a_{o})+\sum_{i=1}^{m} \sum_{j=1}^{q_{i}} x^{2}_{ij}
 \lambda^{-}(a_{ij})\leq b^{-},
  ~~ \max_{i,j}((1-x^{2}_{ij}) \lambda^{-}(a_{ij})+\lambda^{+}(a_{ij}))\leq
 b^{+}.\]
 For the sequence of simple problems above, we can apply
 approximation algorithms, which are based on an
 \(\epsilon\)-approximate algorithm
 (\(\epsilon\in[0,1]\)) for knapsack problem
 (e.g.,  \cite{keller04}, \cite{mar90},
  \cite{sahni75}, \cite{sahni78}).
 In the case of these algorithms, an estimate of an operation number
 is similar for knapsack problem
 (e.g., \cite{keller04}, \cite{mar90}), and equals
 \(O(\frac{m^{2}}{\epsilon})\) \cite{lev81}.
 The algorithms apply ordering of elements from set
 \(\Upsilon a\) by non-decreasing of
 \(\lambda(a_{i})\) or \((\lambda^{-}(a_{i})+\lambda^{+}(a_{i}))\).

 The solving process of auxiliary problems is based on similar
 approximation approach
 to multiple choice problem with the following estimates of
 number of operations and required memory accordingly
 (e.g., \cite{keller04}, \cite{mar90}):~
 \(O(\frac {m}{\epsilon} \sum_{i=1}^{m}q_{i}),
 ~ O(\frac {m^{2}}{\epsilon} \max_{1 \leq i \leq m} \{q_{i}\})\).
 Unfortunately,
 we could not construct an algorithm with similar
 estimates for the auxiliary problem 2.3 (Fig. 31g) \cite{lev81}.
 As a result,
  \((\epsilon,\delta)\)-approximate
 algorithms with the following estimates
 (number of operations, and required memory) were
% have been
 suggested \cite{lev81}: ~
 \(O(\frac {m}{\epsilon \delta} \sum_{i=1}^{m}q_{i}),
 ~~ O(\frac {m^{2}}{\epsilon} \max_{1 \leq i \leq m}
 \{q_{i}\})\),
 where \(\delta\) is a relative error for constraints.

\begin{center}
\begin{picture}(114,34)
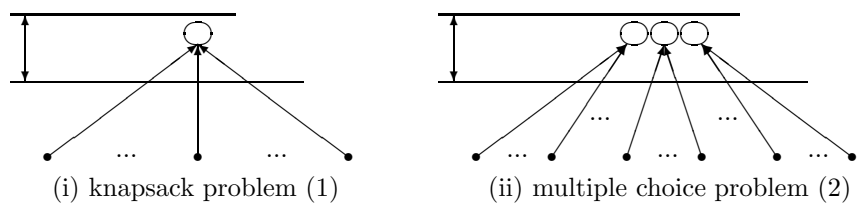
\put(02.4,00){\makebox(0,0)[bl]{Fig. 30. Illustration for knapsack
and multiple choice problems \cite{lev98}}}

%----------------------------- I
\put(05.6,06){\makebox(0,0)[bl]{(i) knapsack problem (1)}}

\put(05,12){\circle*{1}} \put(25,12){\circle*{1}}
\put(45,12){\circle*{1}}

\put(14,12){\makebox(0,0)[bl]{...}}
\put(34,12){\makebox(0,0)[bl]{...}}

\put(25,28.5){\oval(3.6,3)}

\put(25,12){\vector(0,1){15}} \put(05,12){\vector(4,3){20}}
\put(45,12){\vector(-4,3){20}} \put(00,31){\line(1,0){30}}
\put(00,22){\line(1,0){39}}

\put(02,25){\vector(0,-1){3}} \put(02,25){\vector(0,1){06}}

%---------------------- II
\put(63.6,06){\makebox(0,0)[bl]{(ii) multiple choice problem (2)}}

\put(62,12){\circle*{1}} \put(72,12){\circle*{1}}
\put(82,12){\circle*{1}} \put(92,12){\circle*{1}}
\put(102,12){\circle*{1}} \put(112,12){\circle*{1}}

\put(66,12){\makebox(0,0)[bl]{...}}
\put(86,12){\makebox(0,0)[bl]{...}}
\put(106,12){\makebox(0,0)[bl]{...}}
\put(77,17){\makebox(0,0)[bl]{...}}
\put(94,17){\makebox(0,0)[bl]{...}}

\put(83,28.5){\oval(3.6,3)} \put(87,28.5){\oval(3.6,3)}
\put(91,28.5){\oval(3.6,3)}

\put(82,12){\vector(1,3){5}} \put(92,12){\vector(-1,3){5}}
\put(62,12){\vector(4,3){20}} \put(72,12){\vector(2,3){10}}
\put(112,12){\vector(-4,3){20}} \put(102,12){\vector(-2,3){10}}
\put(57,31){\line(1,0){40}} \put(57,22){\line(1,0){49}}

\put(59,25){\vector(0,-1){3}} \put(59,25){\vector(0,1){06}}

\end{picture}
\end{center}

\begin{center}
\begin{picture}(114,140)
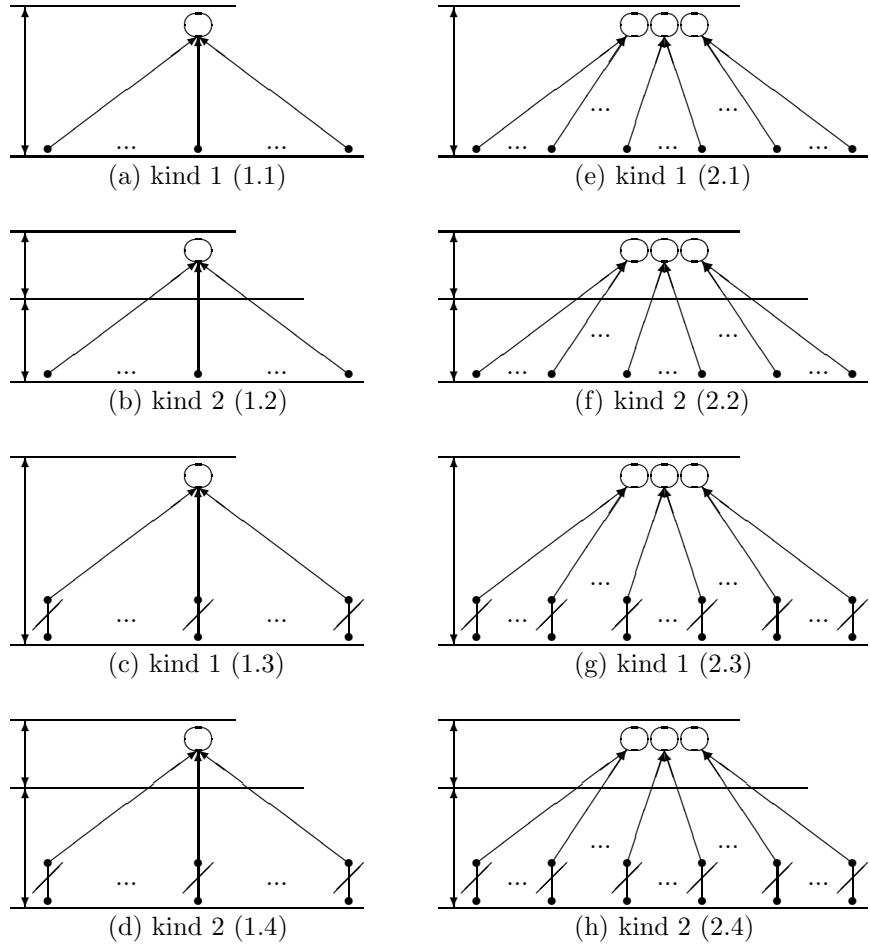
\put(09,04){\makebox(0,0)[bl]{Fig. 31. Illustration for simplest
and auxiliary problems \cite{lev98}}}

%------------------ d
\put(13,09.5){\makebox(0,0)[bl]{(d) kind 2 (1.4)}}

\put(05,15){\circle*{1}} \put(25,15){\circle*{1}}
\put(45,15){\circle*{1}}

\put(05,15){\line(0,1){5}} \put(25,15){\line(0,1){5}}
\put(45,15){\line(0,1){5}} \put(03,16){\line(1,1){4}}
\put(23,16){\line(1,1){4}} \put(43,16){\line(1,1){4}}
\put(05,20){\circle*{1}} \put(25,20){\circle*{1}}
\put(45,20){\circle*{1}}

\put(14,17){\makebox(0,0)[bl]{...}}
\put(34,17){\makebox(0,0)[bl]{...}}

%\put(25,36){\circle{2}}

\put(25,36.5){\oval(3.6,3)}

\put(25,20){\vector(0,1){15}} \put(05,20){\vector(4,3){20}}
\put(45,20){\vector(-4,3){20}} \put(00,39){\line(1,0){30}}
\put(00,30){\line(1,0){39}} \put(00,14){\line(1,0){47}}
\put(02,20){\vector(0,-1){6}} \put(02,20){\vector(0,1){10}}
\put(02,35){\vector(0,-1){5}} \put(02,35){\vector(0,1){4}}

%---------------------------c
\put(13,44.5){\makebox(0,0)[bl]{(c) kind 1 (1.3)}}

\put(05,50){\circle*{1}} \put(25,50){\circle*{1}}
\put(45,50){\circle*{1}} \put(05,50){\line(0,1){5}}
\put(25,50){\line(0,1){5}} \put(45,50){\line(0,1){5}}
\put(03,51){\line(1,1){4}} \put(23,51){\line(1,1){4}}
\put(43,51){\line(1,1){4}} \put(05,55){\circle*{1}}
\put(25,55){\circle*{1}} \put(45,55){\circle*{1}}

\put(14,52){\makebox(0,0)[bl]{...}}
\put(34,52){\makebox(0,0)[bl]{...}}

%\put(25,71){\circle{2}}

\put(25,71.5){\oval(3.6,3)}

\put(25,55){\vector(0,1){15}} \put(05,55){\vector(4,3){20}}
\put(45,55){\vector(-4,3){20}} \put(00,74){\line(1,0){30}}
\put(00,49){\line(1,0){47}} \put(02,55){\vector(0,-1){6}}
\put(02,55){\vector(0,1){19}}

%----------------------------- b
\put(13,79.5){\makebox(0,0)[bl]{(b) kind 2 (1.2)}}

\put(05,85){\circle*{1}} \put(25,85){\circle*{1}}
\put(45,85){\circle*{1}}

\put(14,85){\makebox(0,0)[bl]{...}}
\put(34,85){\makebox(0,0)[bl]{...}}

%\put(25,101){\circle{2}}

\put(25,101.5){\oval(3.6,3)}

\put(25,85){\vector(0,1){15}} \put(05,85){\vector(4,3){20}}
\put(45,85){\vector(-4,3){20}} \put(00,104){\line(1,0){30}}
\put(00,095){\line(1,0){39}} \put(00,84){\line(1,0){47}}
\put(02,90){\vector(0,-1){06}} \put(02,90){\vector(0,1){05}}
\put(02,098){\vector(0,-1){3}} \put(02,098){\vector(0,1){06}}

%------------------------------ a
\put(13,109.5){\makebox(0,0)[bl]{(a) kind 1 (1.1)}}

\put(05,115){\circle*{1}} \put(25,115){\circle*{1}}
\put(45,115){\circle*{1}} \put(14,115){\makebox(0,0)[bl]{...}}
\put(34,115){\makebox(0,0)[bl]{...}}

%\put(25,131){\circle{2}}

\put(25,131.5){\oval(3.6,3)}

\put(05,115){\vector(4,3){20}} \put(25,115){\vector(0,1){15}}
\put(45,115){\vector(-4,3){20}} \put(00,134){\line(1,0){30}}
\put(00,114){\line(1,0){47}} \put(02,128){\vector(0,-1){14}}
\put(02,128){\vector(0,1){06}}

%-------------------- h
\put(75,09.5){\makebox(0,0)[bl]{(h) kind 2 (2.4)}}

\put(62,15){\circle*{1}} \put(72,15){\circle*{1}}
\put(82,15){\circle*{1}} \put(92,15){\circle*{1}}
\put(102,15){\circle*{1}} \put(112,15){\circle*{1}}
\put(62,15){\line(0,1){5}} \put(72,15){\line(0,1){5}}
\put(82,15){\line(0,1){5}} \put(92,15){\line(0,1){5}}
\put(102,15){\line(0,1){5}} \put(112,15){\line(0,1){5}}
\put(60,16){\line(1,1){4}} \put(70,16){\line(1,1){4}}
\put(80,16){\line(1,1){4}} \put(90,16){\line(1,1){4}}
\put(100,16){\line(1,1){4}} \put(110,16){\line(1,1){4}}
\put(62,20){\circle*{1}} \put(72,20){\circle*{1}}
\put(82,20){\circle*{1}} \put(92,20){\circle*{1}}
\put(102,20){\circle*{1}} \put(112,20){\circle*{1}}

\put(66,17){\makebox(0,0)[bl]{...}}
\put(86,17){\makebox(0,0)[bl]{...}}
\put(106,17){\makebox(0,0)[bl]{...}}
\put(77,22){\makebox(0,0)[bl]{...}}
\put(94,22){\makebox(0,0)[bl]{...}}

%\put(83,36){\circle{2}} \put(87,36){\circle{2}}
%\put(91,36){\circle{2}}

\put(83,36.5){\oval(3.6,3)} \put(87,36.5){\oval(3.6,3)}
\put(91,36.5){\oval(3.6,3)}

\put(82,20){\vector(1,3){5}} \put(92,20){\vector(-1,3){5}}
\put(62,20){\vector(4,3){20}} \put(72,20){\vector(2,3){10}}
\put(112,20){\vector(-4,3){20}} \put(102,20){\vector(-2,3){10}}
\put(57,39){\line(1,0){40}} \put(57,30){\line(1,0){49}}
\put(57,14){\line(1,0){57}} \put(59,20){\vector(0,-1){6}}
\put(59,20){\vector(0,1){10}} \put(59,35){\vector(0,-1){5}}
\put(59,35){\vector(0,1){4}}

%----------------------- g
\put(75,44.5){\makebox(0,0)[bl]{(g) kind 1 (2.3)}}

\put(62,50){\circle*{1}} \put(72,50){\circle*{1}}
\put(82,50){\circle*{1}} \put(92,50){\circle*{1}}
\put(102,50){\circle*{1}} \put(112,50){\circle*{1}}
\put(62,50){\line(0,1){5}} \put(72,50){\line(0,1){5}}
\put(82,50){\line(0,1){5}} \put(92,50){\line(0,1){5}}
\put(102,50){\line(0,1){5}} \put(112,50){\line(0,1){5}}
\put(60,51){\line(1,1){4}} \put(70,51){\line(1,1){4}}
\put(80,51){\line(1,1){4}} \put(90,51){\line(1,1){4}}
\put(100,51){\line(1,1){4}} \put(110,51){\line(1,1){4}}
\put(62,55){\circle*{1}} \put(72,55){\circle*{1}}
\put(82,55){\circle*{1}} \put(92,55){\circle*{1}}
\put(102,55){\circle*{1}} \put(112,55){\circle*{1}}

\put(66,52){\makebox(0,0)[bl]{...}}
\put(86,52){\makebox(0,0)[bl]{...}}
\put(106,52){\makebox(0,0)[bl]{...}}
\put(77,57){\makebox(0,0)[bl]{...}}
\put(94,57){\makebox(0,0)[bl]{...}}

%\put(83,71){\circle{2}} \put(87,71){\circle{2}}
%\put(91,71){\circle{2}}

\put(83,71.5){\oval(3.6,3)} \put(87,71.5){\oval(3.6,3)}
\put(91,71.5){\oval(3.6,3)}

\put(82,55){\vector(1,3){5}} \put(92,55){\vector(-1,3){5}}
\put(62,55){\vector(4,3){20}} \put(72,55){\vector(2,3){10}}
\put(112,55){\vector(-4,3){20}} \put(102,55){\vector(-2,3){10}}
\put(57,74){\line(1,0){40}} \put(57,49){\line(1,0){57}}
\put(59,55){\vector(0,-1){6}} \put(59,55){\vector(0,1){19}}

%---------------------- f
\put(75,79.5){\makebox(0,0)[bl]{(f) kind 2 (2.2)}}

\put(62,85){\circle*{1}} \put(72,85){\circle*{1}}
\put(82,85){\circle*{1}} \put(92,85){\circle*{1}}
\put(102,85){\circle*{1}} \put(112,85){\circle*{1}}

\put(66,85){\makebox(0,0)[bl]{...}}
\put(86,85){\makebox(0,0)[bl]{...}}
\put(106,85){\makebox(0,0)[bl]{...}}
\put(77,90){\makebox(0,0)[bl]{...}}
\put(94,90){\makebox(0,0)[bl]{...}}

%\put(83,101){\circle{2}} \put(87,101){\circle{2}}
%\put(91,101){\circle{2}}

\put(83,101.5){\oval(3.6,3)} \put(87,101.5){\oval(3.6,3)}
\put(91,101.5){\oval(3.6,3)}

\put(82,85){\vector(1,3){5}} \put(92,85){\vector(-1,3){5}}
\put(62,85){\vector(4,3){20}} \put(72,85){\vector(2,3){10}}
\put(112,85){\vector(-4,3){20}} \put(102,85){\vector(-2,3){10}}
\put(57,104){\line(1,0){40}} \put(57,095){\line(1,0){49}}
\put(57,84){\line(1,0){57}} \put(59,90){\vector(0,-1){06}}
\put(59,90){\vector(0,1){05}} \put(59,098){\vector(0,-1){3}}
\put(59,098){\vector(0,1){06}}

%----------------- E
\put(75,109.5){\makebox(0,0)[bl]{(e) kind 1 (2.1)}}

\put(62,115){\circle*{1}} \put(72,115){\circle*{1}}
\put(82,115){\circle*{1}} \put(92,115){\circle*{1}}
\put(102,115){\circle*{1}} \put(112,115){\circle*{1}}

\put(66,115){\makebox(0,0)[bl]{...}}
\put(86,115){\makebox(0,0)[bl]{...}}
\put(106,115){\makebox(0,0)[bl]{...}}
\put(77,120){\makebox(0,0)[bl]{...}}
\put(94,120){\makebox(0,0)[bl]{...}}

%\put(83,131){\circle{2}} \put(87,131){\circle{2}}
%\put(91,131){\circle{2}}

\put(83,131.5){\oval(3.6,3)} \put(87,131.5){\oval(3.6,3)}
\put(91,131.5){\oval(3.6,3)}

\put(82,115){\vector(1,3){5}} \put(92,115){\vector(-1,3){5}}
\put(62,115){\vector(4,3){20}} \put(72,115){\vector(2,3){10}}
\put(112,115){\vector(-4,3){20}} \put(102,115){\vector(-2,3){10}}
\put(57,134){\line(1,0){40}} \put(57,114){\line(1,0){57}}
\put(59,128){\vector(0,-1){14}} \put(59,128){\vector(0,1){06}}

\end{picture}
\end{center}

 \begin{center}
\begin{picture}(80,56)

\put(03.5,52){\makebox(0,0)[bl]{Table 7.
 Knapsack-like problems - illustrations}}

\put(00,0){\line(1,0){80}} \put(00,43){\line(1,0){80}}
\put(00,50){\line(1,0){80}}

\put(00,0){\line(0,1){50}} \put(61,0){\line(0,1){50}}
\put(80,0){\line(0,1){50}}

\put(01,45){\makebox(0,0)[bl]{Problem}}
\put(62,45){\makebox(0,0)[bl]{Illustration }}

%------------------------------------
%------------------------------------

\put(01,37.4){\makebox(0,0)[bl]{4.Basic knapsack problem
 ({\it problem 1})}}
\put(62,38){\makebox(0,0)[bl]{Fig. 30i}}

%-------------
\put(01,34){\makebox(0,0)[bl]{5.{\it Problem 1.1} }}
%\put(56,02){\makebox(0,0)[bl]{}}
\put(62,34){\makebox(0,0)[bl]{Fig. 31.a}}

%-------------
\put(01,30){\makebox(0,0)[bl]{5.{\it Problem 1.2} }}
%\put(56,02){\makebox(0,0)[bl]{}}
\put(62,30){\makebox(0,0)[bl]{Fig. 31.b}}

%-------------
\put(01,26){\makebox(0,0)[bl]{5.{\it Problem 1.3} }}
%\put(56,02){\makebox(0,0)[bl]{}}
\put(62,26){\makebox(0,0)[bl]{Fig. 31.c}}

%-------------
\put(01,21.6){\makebox(0,0)[bl]{5.{\it Problem 1.4} }}
%\put(56,02){\makebox(0,0)[bl]{}}
\put(62,22){\makebox(0,0)[bl]{Fig. 31.d}}

%------------------------------------

\put(01,17.4){\makebox(0,0)[bl]{4.Multiple choice problem
 ({\it problem 2})}}
%\put(56,6){\makebox(0,0)[bl]{}}
\put(62,18){\makebox(0,0)[bl]{Fig. 30ii}}

%-------------
\put(01,14){\makebox(0,0)[bl]{5.{\it Problem 2.1} }}
%\put(56,02){\makebox(0,0)[bl]{}}
\put(62,14){\makebox(0,0)[bl]{Fig. 31.e}}

%-------------
\put(01,10){\makebox(0,0)[bl]{5.{\it Problem 2.2} }}
%\put(56,02){\makebox(0,0)[bl]{}}
\put(62,10){\makebox(0,0)[bl]{Fig. 31.f}}

%-------------
\put(01,06){\makebox(0,0)[bl]{5.{\it Problem 2.3} }}
%\put(56,02){\makebox(0,0)[bl]{}}
\put(62,06){\makebox(0,0)[bl]{Fig. 31.g}}

%-------------
\put(01,01.6){\makebox(0,0)[bl]{5.{\it Problem 2.4} }}
%\put(56,02){\makebox(0,0)[bl]{}}
\put(62,02){\makebox(0,0)[bl]{Fig. 31.h}}

\end{picture}
\end{center}

% In more general case,
 When \(G\) is a \(k\)-level tree,
 the algorithm is based on
 cascade-like 'Bottom-Up' process  (Fig. 32) \cite{lev81}:

%~~

  {\em Step 1.} Problem 1.2.

  ~~~~~~ {\bf .~ .~ .}

  {\em Step} \(j\)~ \((j= \overline{2,k-2})\). Problem 2.4

  ~~~~~~ {\bf .~ .~ .}

  {\em Step} \((k-1)\). Problem 2.3.

%~~

 Estimates of the algorithms are as follows
 (i.e., operations, and memory):~
 \(O(\frac {n^{2}\eta^{5}(a_{o})}{\epsilon \delta^{4}}),
 ~ O(\frac {m^{2}\eta^{4}}{\epsilon \delta^{4}})\),
 where \(m(a)=|\Upsilon (a)|\), \(m=max_{a\in A} m(a)\),
 \(\eta (a)=\mid  A_{a} \setminus \{a^{'}\in A_{a} \mid
 \Upsilon a^{'}=\O \}\mid \).

 In the case of \(3\)-level tree,
 the estimate of the operation number is:~~
 \(O(\frac {n^{2}\eta^{4}(a_{o})}{\epsilon \delta^{3}})\).

 Prospective generalizations of the problem may involve the following:
 (a) multicriteria descriptions of the elements,
 (b) more complicated structure (e.g., parallel-series graph),
 (c) uncertainty, and
 (d) dynamics.

\begin{center}
\begin{picture}(103,75)
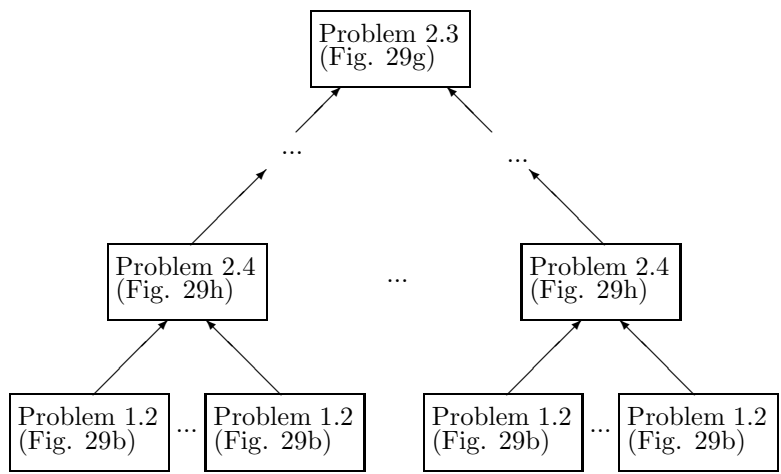
\put(05,00){\makebox(0,0)[bl]{Fig. 32. 'Bottom-Up' solving scheme
 for tree-like case \cite{lev98}}}

\put(00,10){\framebox(21,10)[bl]{}}
\put(01,16){\makebox(0,0)[bl]{Problem 1.2}}
\put(01,12){\makebox(0,0)[bl]{(Fig. 29b)}}
\put(22,15){\makebox(0,0)[bl]{...}} \put(11,20){\vector(1,1){10}}

\put(26,10){\framebox(21,10)[bl]{}}
\put(27,16){\makebox(0,0)[bl]{Problem 1.2}}
\put(27,12){\makebox(0,0)[bl]{(Fig. 29b)}}
\put(36,20){\vector(-1,1){10}}

\put(13,30){\framebox(21,10)[bl]{}}
\put(14,36){\makebox(0,0)[bl]{Problem 2.4}}
\put(14,32){\makebox(0,0)[bl]{(Fig. 29h)}}
\put(24,40){\vector(1,1){10}}

\put(50,35){\makebox(0,0)[bl]{...}}

\put(55,10){\framebox(21,10)[bl]{}}
\put(56,16){\makebox(0,0)[bl]{Problem 1.2}}
\put(56,12){\makebox(0,0)[bl]{(Fig. 29b)}}
\put(77,15){\makebox(0,0)[bl]{...}} \put(66,20){\vector(1,1){10}}

\put(81,10){\framebox(21,10)[bl]{}}
\put(82,16){\makebox(0,0)[bl]{Problem 1.2}}
\put(82,12){\makebox(0,0)[bl]{(Fig. 29b)}}
\put(91,20){\vector(-1,1){10}}

\put(68,30){\framebox(21,10)[bl]{}}
\put(69,36){\makebox(0,0)[bl]{Problem 2.4}}
\put(69,32){\makebox(0,0)[bl]{(Fig. 29h)}}
\put(79,40){\vector(-1,1){10}}

\put(40,61){\framebox(21,10)[bl]{}}
\put(41,67){\makebox(0,0)[bl]{Problem 2.3}}
\put(41,63){\makebox(0,0)[bl]{(Fig. 29g)}}

\put(38,55){\vector( 1,1){6 }} \put(64,55){\vector(-1,1){6 }}

\put(36,52){\makebox(0,0)[bl]{...}}
\put(66,51){\makebox(0,0)[bl]{...}}

\end{picture}
\end{center}

\subsection{Hotlink Assignment Problem}

 In recent decade,
 the concept of ``hotlinks'' has been introduced
 for decreasing the complexity of  access
 in web directories (or similar information structures)
 via usage (inserting) of a limited set of additional
 hyperlinks  (i.e., ``hotlinks'')  to data
 (e.g., \cite{bose01},\cite{dou10},\cite{perk99}).
 In  general, ``hotlink assignment problem''
 is a network upgrade problem  (e.g., \cite{fuhr01}):

%~~

 {\it Find additional new arc(s) to the initial graph
 in order to insert shortcuts and decrease the expected path
 length.}

%~~

 Mainly, the problem is examined for trees.
 Let \(T = (A,E)\) be a directed tree with maximum degree \(d\),
 rooted at a node \(r_{0}\in A\)
 (elements of \(A\) correspond to Web sites, elements of \(E\) correspond to hyperlinks).
 A node weight equals its access (search) frequency (probability).
 It is assumed that required information is contained at
 leaf nodes (for simplicity).
 The length of the search for node \(v \in A\) equals
 the number of links in the path from \(r_{0}\) to \(v\).

 Let \(T_{u} = (A_{u},E_{u})\) be a subtree of \(T\)
 (\(A_{u} \subseteq A,  E_{u} \subseteq E \)), rooted at node \(u\in A\)
 (here \(u\) is not the son of \(r_{0}\)).
 Thus, additional direct link (``hotlink'') will be as follows:
 (\(r_{0},u\)).
 In this case, a path to all leaf nodes in \(T_{u}\) will be smaller.

 Fig. 33, Fig. 34, Fig. 35, and Fig. 36 illustrate
 the simplest versions of ``hotlink assignment'' problem.
 Note, Fig. 36 depicts the usage of internal nodes as additional roots
 (i.e., \(r_{1}\)).

\begin{center}
%\begin{picture}(50,40)
\begin{picture}(52,36)

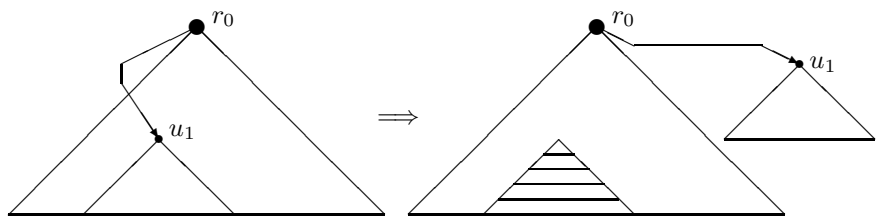
\put(20,00){\makebox(0,0)[bl]{Fig 33. Hotlink assignment problem
(one hotlink)}}

%-----

\put(25,31){\circle*{2}}

\put(27,31){\makebox(0,0)[bl]{\(r_{0}\)}}

\put(00,06){\line(1,1){25}} \put(50,06){\line(-1,1){25}}

\put(00,06){\line(1,0){50}}

%-----

\put(25,31){\line(-2,-1){10}} \put(15,26){\line(0,-1){2.5}}
\put(15,23.5){\vector(2,-3){5}}
%\put(25,31){\line(-2,1){10}}

%-----

\put(20,16){\circle*{1.2}}

\put(21.3,16){\makebox(0,0)[bl]{\(u_{1}\)}}

\put(10,06){\line(1,1){10}} \put(30,06){\line(-1,1){10}}

\put(10,06){\line(1,0){20}}

%%%%%%%%%%%%%%%%%%%%%

\put(49,18){\makebox(0,0)[bl]{\(\Longrightarrow\)}}

\end{picture}
%\end{center}
%
%\begin{center}
%\begin{picture}(50,40)
\begin{picture}(62,36)

%-----

\put(25,31){\circle*{2}}

\put(27,31){\makebox(0,0)[bl]{\(r_{0}\)}}

\put(00,06){\line(1,1){25}} \put(50,06){\line(-1,1){25}}

\put(00,06){\line(1,0){50}}

%-----

%\put(25,31){\line(-2,-1){10}} \put(15,26){\line(0,-1){2.5}}
%\put(15,23.5){\vector(2,-3){5}}

%-----

%\put(20,16){\circle*{1.2}}
%\put(31,44){\makebox(0,0)[bl]{\(r_{1}\)}}

\put(10,06){\line(1,1){10}} \put(30,06){\line(-1,1){10}}

%--

\put(10,06){\line(1,0){20}}

%----

\put(18,14){\line(1,0){04}} \put(16,12){\line(1,0){08}}
\put(14,10){\line(1,0){12}} \put(12,08){\line(1,0){16}}

%-----

\put(25,31){\line(2,-1){5}} \put(30,28.5){\line(1,0){17}}
\put(47,28.5){\vector(2,-1){5}}

%---

\put(52,26){\circle*{1.2}}

\put(53.3,25){\makebox(0,0)[bl]{\(u_{1}\)}}

\put(42,16){\line(1,1){10}} \put(62,16){\line(-1,1){10}}

\put(42,16){\line(1,0){20}}

%-------------

\end{picture}
\end{center}

\begin{center}
%\begin{picture}(50,40)
\begin{picture}(58,36)

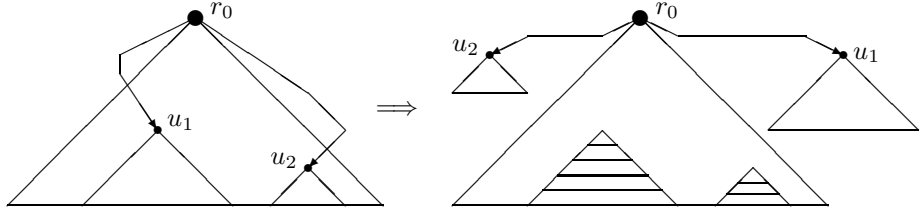
\put(20,00){\makebox(0,0)[bl]{Fig 34. Hotlink assignment problem
(two hotlinks)}}

%-----

\put(25,31){\circle*{2}}

\put(27,31){\makebox(0,0)[bl]{\(r_{0}\)}}

\put(00,06){\line(1,1){25}} \put(50,06){\line(-1,1){25}}

\put(00,06){\line(1,0){50}}

%-----

\put(25,31){\line(-2,-1){10}} \put(15,26){\line(0,-1){2.5}}
\put(15,23.5){\vector(2,-3){5}}
%\put(25,31){\line(-2,1){10}}

%-----

\put(20,16){\circle*{1.2}}

\put(21.3,16){\makebox(0,0)[bl]{\(u_{1}\)}}

\put(10,06){\line(1,1){10}} \put(30,06){\line(-1,1){10}}

\put(10,06){\line(1,0){20}}

%-----

\put(25,31){\line(3,-2){15}} \put(40,21){\line(1,-1){5}}
\put(45,16){\vector(-1,-1){5}}

%--

\put(40,11){\circle*{1.2}}

\put(35,11){\makebox(0,0)[bl]{\(u_{2}\)}}

\put(35,06){\line(1,1){05}} \put(45,06){\line(-1,1){05}}

\put(35,06){\line(1,0){10}}

%%%%%%%%%%%%%%%%%%%%%

\put(49,18){\makebox(0,0)[bl]{\(\Longrightarrow\)}}

\end{picture}
%\end{center}
%
%\begin{center}
%\begin{picture}(50,40)
\begin{picture}(62,36)

%-----

\put(25,31){\circle*{2}}

\put(27,31){\makebox(0,0)[bl]{\(r_{0}\)}}

\put(00,06){\line(1,1){25}} \put(50,06){\line(-1,1){25}}

\put(00,06){\line(1,0){50}}

%-----

%\put(25,31){\line(-2,-1){10}} \put(15,26){\line(0,-1){2.5}}
%\put(15,23.5){\vector(2,-3){5}}

%-----

%\put(20,16){\circle*{1.2}}
%\put(31,44){\makebox(0,0)[bl]{\(r_{1}\)}}

\put(10,06){\line(1,1){10}} \put(30,06){\line(-1,1){10}}

%--

\put(10,06){\line(1,0){20}}

%----

\put(18,14){\line(1,0){04}} \put(16,12){\line(1,0){08}}
\put(14,10){\line(1,0){12}} \put(12,08){\line(1,0){16}}

%-----

\put(25,31){\line(2,-1){5}} \put(30,28.5){\line(1,0){17}}
\put(47,28.5){\vector(2,-1){5}}

%---

\put(52,26){\circle*{1.2}}

\put(53.3,25){\makebox(0,0)[bl]{\(u_{1}\)}}

\put(42,16){\line(1,1){10}} \put(62,16){\line(-1,1){10}}
\put(42,16){\line(1,0){20}}

%-------------

%\put(40,11){\circle*{1.2}}
%\put(35,11){\makebox(0,0)[bl]{\(u_{2}\)}}

\put(35,06){\line(1,1){05}} \put(45,06){\line(-1,1){05}}
\put(35,06){\line(1,0){10}}

%--

\put(36.5,07.5){\line(1,0){7}} \put(38,09){\line(1,0){4}}

\put(25,31){\line(-2,-1){5}} \put(20,28.5){\line(-1,0){10}}
\put(10,28.5){\vector(-2,-1){5}}

%--

\put(05,26){\circle*{1.2}}
\put(00,26){\makebox(0,0)[bl]{\(u_{2}\)}}

\put(00,21){\line(1,1){05}} \put(10,21){\line(-1,1){05}}

\put(00,21){\line(1,0){10}}

\end{picture}
\end{center}

\begin{center}
%\begin{picture}(50,40)
\begin{picture}(58,36)

\put(20,00){\makebox(0,0)[bl]{Fig 35. Hotlink assignment problem
(three hotlinks)}}

%-----

\put(25,31){\circle*{2}}

\put(27,31){\makebox(0,0)[bl]{\(r_{0}\)}}

\put(00,06){\line(1,1){25}} \put(50,06){\line(-1,1){25}}

\put(00,06){\line(1,0){50}}

%-----TO u1

\put(25,31){\line(-2,-1){10}} \put(15,26){\line(0,-1){2.5}}
\put(15,23.5){\vector(2,-3){5}}

%-----TO u3 (20,11)

\put(25,31){\line(-3,-1){15}} \put(10,26){\line(0,-1){10}}
\put(10,16){\vector(2,-1){10}}

%-----

\put(20,16){\circle*{1.2}}

\put(21.3,16){\makebox(0,0)[bl]{\(u_{1}\)}}

\put(10,06){\line(1,1){10}} \put(30,06){\line(-1,1){10}}

\put(10,06){\line(1,0){20}}

%-----

\put(25,31){\line(3,-2){15}} \put(40,21){\line(1,-1){5}}
\put(45,16){\vector(-1,-1){5}}

%--

\put(40,11){\circle*{1.2}}

\put(35,11){\makebox(0,0)[bl]{\(u_{2}\)}}

\put(35,06){\line(1,1){05}} \put(45,06){\line(-1,1){05}}

\put(35,06){\line(1,0){10}}

%--

\put(20,11){\circle*{1.2}}

\put(18.6,12){\makebox(0,0)[bl]{\(u_{3}\)}}

\put(15,06){\line(1,1){05}} \put(25,06){\line(-1,1){05}}

\put(15,06){\line(1,0){10}}

%%%%%%%%%%%%%%%%%%%%%

\put(49,18){\makebox(0,0)[bl]{\(\Longrightarrow\)}}

\end{picture}
%\end{center}
%
%\begin{center}
%\begin{picture}(50,40)
\begin{picture}(62,36)

%-----

\put(25,31){\circle*{2}}

\put(27,31){\makebox(0,0)[bl]{\(r_{0}\)}}

\put(00,06){\line(1,1){25}} \put(50,06){\line(-1,1){25}}

\put(00,06){\line(1,0){50}}

%-----

%\put(25,31){\line(-2,-1){10}} \put(15,26){\line(0,-1){2.5}}
%\put(15,23.5){\vector(2,-3){5}}

%-----

%\put(20,16){\circle*{1.2}}
%\put(31,44){\makebox(0,0)[bl]{\(r_{1}\)}}

\put(10,06){\line(1,1){10}} \put(30,06){\line(-1,1){10}}

%--

\put(10,06){\line(1,0){20}}

%----

\put(18,14){\line(1,0){04}} \put(16,12){\line(1,0){08}}
\put(14,10){\line(1,0){12}} \put(12,08){\line(1,0){16}}

%-----

\put(25,31){\line(2,-1){5}} \put(30,28.5){\line(1,0){17}}
\put(47,28.5){\vector(2,-1){5}}

%---

\put(52,26){\circle*{1.2}}

\put(53.3,25){\makebox(0,0)[bl]{\(u_{1}\)}}

\put(42,16){\line(1,1){10}} \put(62,16){\line(-1,1){10}}
\put(42,16){\line(1,0){20}}

%--u1 minus u3

%\put(40,11){\circle*{1.2}}
%\put(35,11){\makebox(0,0)[bl]{\(u_{2}\)}}

\put(47,16){\line(1,1){05}} \put(57,16){\line(-1,1){05}}
\put(47,16){\line(1,0){10}}

%--

\put(48.5,17.5){\line(1,0){7}} \put(50,19){\line(1,0){4}}

%-------------

%--

%\put(40,11){\circle*{1.2}}
%\put(35,11){\makebox(0,0)[bl]{\(u_{2}\)}}

\put(35,06){\line(1,1){05}} \put(45,06){\line(-1,1){05}}
\put(35,06){\line(1,0){10}}

%--

\put(36.5,07.5){\line(1,0){7}} \put(38,09){\line(1,0){4}}

%-------------TO u2

\put(25,31){\line(-2,-1){5}} \put(20,28.5){\vector(-2,0){10}}
%--

\put(10,28.5){\circle*{1.2}}
\put(005,28.5){\makebox(0,0)[bl]{\(u_{2}\)}}

\put(05,23.5){\line(1,1){05}} \put(15,23.5){\line(-1,1){05}}
\put(05,23.5){\line(1,0){10}}

%-------------TO u3

\put(25,31){\line(-2,1){5}}

\put(20,33.5){\line(-1,0){15}} \put(05,33.5){\line(-2,-1){5}}

\put(00,31){\line(0,-1){4}}

\put(00,27){\vector(1,-1){5}}
%--

\put(05,22){\circle*{1.2}}
\put(06.4,20.5){\makebox(0,0)[bl]{\(u_{3}\)}}

\put(00,17){\line(1,1){05}} \put(10,17){\line(-1,1){05}}
\put(00,17){\line(1,0){10}}

\end{picture}
\end{center}

\begin{center}
%\begin{picture}(50,40)
\begin{picture}(58,36)

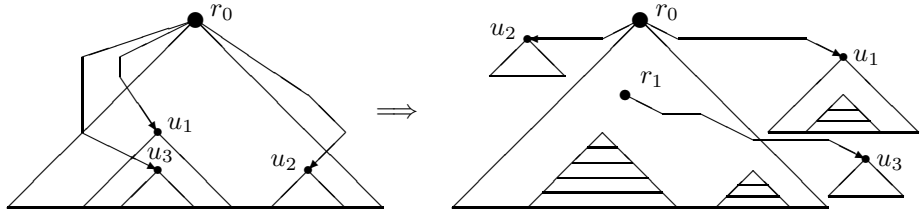
\put(06,00){\makebox(0,0)[bl]{Fig 36. Hotlink assignment problem
(three hotlinks, one internal root)}}

%-----

\put(25,31){\circle*{2}}

\put(27,31){\makebox(0,0)[bl]{\(r_{0}\)}}

\put(00,06){\line(1,1){25}} \put(50,06){\line(-1,1){25}}

\put(00,06){\line(1,0){50}}

%-----TO u1

\put(25,31){\line(-2,-1){10}} \put(15,26){\line(0,-1){2.5}}
\put(15,23.5){\vector(2,-3){5}}

%-----TO u3 (20,11)

\put(25,31){\line(-3,-1){15}} \put(10,26){\line(0,-1){10}}
\put(10,16){\vector(2,-1){10}}

%-----

\put(20,16){\circle*{1.2}}

\put(21.3,16){\makebox(0,0)[bl]{\(u_{1}\)}}

\put(10,06){\line(1,1){10}} \put(30,06){\line(-1,1){10}}

\put(10,06){\line(1,0){20}}

%-----

\put(25,31){\line(3,-2){15}} \put(40,21){\line(1,-1){5}}
\put(45,16){\vector(-1,-1){5}}

%--

\put(40,11){\circle*{1.2}}

\put(35,11){\makebox(0,0)[bl]{\(u_{2}\)}}

\put(35,06){\line(1,1){05}} \put(45,06){\line(-1,1){05}}

\put(35,06){\line(1,0){10}}

%--

\put(20,11){\circle*{1.2}}

\put(18.6,12){\makebox(0,0)[bl]{\(u_{3}\)}}

\put(15,06){\line(1,1){05}} \put(25,06){\line(-1,1){05}}

\put(15,06){\line(1,0){10}}

%%%%%%%%%%%%%%%%%%%%%

\put(49,18){\makebox(0,0)[bl]{\(\Longrightarrow\)}}

\end{picture}
%\end{center}
%
%\begin{center}
%\begin{picture}(50,40)
\begin{picture}(62,36)

%-----

\put(25,31){\circle*{2}}

\put(27,31){\makebox(0,0)[bl]{\(r_{0}\)}}

\put(00,06){\line(1,1){25}} \put(50,06){\line(-1,1){25}}

\put(00,06){\line(1,0){50}}

%-----

\put(10,06){\line(1,1){10}} \put(30,06){\line(-1,1){10}}

%--

\put(10,06){\line(1,0){20}}

%----

\put(18,14){\line(1,0){04}} \put(16,12){\line(1,0){08}}
\put(14,10){\line(1,0){12}} \put(12,08){\line(1,0){16}}

%-----

\put(25,31){\line(2,-1){5}} \put(30,28.5){\line(1,0){17}}
\put(47,28.5){\vector(2,-1){5}}

%---

\put(52,26){\circle*{1.2}}

\put(53.3,25){\makebox(0,0)[bl]{\(u_{1}\)}}

\put(42,16){\line(1,1){10}} \put(62,16){\line(-1,1){10}}
\put(42,16){\line(1,0){20}}

%--u1 minus u3

\put(47,16){\line(1,1){05}} \put(57,16){\line(-1,1){05}}
\put(47,16){\line(1,0){10}}

%--

\put(48.5,17.5){\line(1,0){7}} \put(50,19){\line(1,0){4}}

%------------

\put(35,06){\line(1,1){05}} \put(45,06){\line(-1,1){05}}
\put(35,06){\line(1,0){10}}

%--

\put(36.5,07.5){\line(1,0){7}} \put(38,09){\line(1,0){4}}

%-------------TO u2

\put(25,31){\line(-2,-1){5}} \put(20,28.5){\vector(-2,0){10}}
%--

\put(10,28.5){\circle*{1.2}}
\put(005,28.5){\makebox(0,0)[bl]{\(u_{2}\)}}

\put(05,23.5){\line(1,1){05}} \put(15,23.5){\line(-1,1){05}}
\put(05,23.5){\line(1,0){10}}

%-------------TO u3

\put(23,21){\circle*{1.4}}

\put(25,22){\makebox(0,0)[bl]{\(r_{1}\)}}

\put(23,21){\line(2,-1){5}} \put(28,18.5){\line(1,0){5}}

\put(33,18.5){\line(2,-1){7}}

\put(40,15){\line(1,0){10}}
%\put(45,19.5){\line(0,-1){4.5}}

\put(50,15){\vector(2,-1){5}}
%--

\put(55,12.5){\circle*{1.2}}
\put(56.4,11){\makebox(0,0)[bl]{\(u_{3}\)}}

\put(50,07.5){\line(1,1){05}} \put(60,07.5){\line(-1,1){05}}
\put(50,07.5){\line(1,0){10}}

\end{picture}
\end{center}

 The basic ``hotlink assignment'' problem
 consists in assignment of \(k\) additional ``hotlinks'' (from the
 root) to minimize the total number of steps to visit
 the required information nodes.
 On the other hand, it is necessary to find the set of \(k\)
 nodes (\(U=\{u\}\)) at the tree \(T\).
 In recent years, various ``hotlink assignment'' problems
 have been intensively studied (e.g., algorithm design,
 issues of complexity, approximation)
 (e.g., \cite{bose01},\cite{czy03},\cite{dou08},\cite{dou10},
 \cite{fuhr01},\cite{jacobs12},\cite{kra04},
 \cite{laber06},\cite{mat04},\cite{perk99}).
 Some versions of the  problem are presented in Table 8.
 Mainly, ``hotlink assignment'' problems belong to class of
 NP-hard problems
 (e.g.,
 \cite{jacobs12}).
 Many approximation algorithms have been suggested
 for the problems (including FPTAS)
(e.g.,
 \cite{kra01}, \cite{kra04}, \cite{mat04}).

 In the case of multiple attribute
 description of ``hotlinks'' or/and
 selected nodes/subtrees (i.e., nodes as \(u\)),
 multicriteria knapsack like problems or
 multicriteria generalized assignment problems
 may be used
 (e.g., \cite{levpet10a},
 \cite{levsaf10a}, \cite{lust10}, \cite{parra02}).
 It may be reasonable to examine some generalizations of the
 ``hotlink assignment'' problem, for example:
 (i) taking into account uncertainty,
 (ii) tree-inclusion problem (as ``hot-tree assignment''
 problem).
  Application of  ``hotlink assignment'' problems is very useful for many domains,
  e.g.,
  adaptive Web sites systems,
  knowledge bases,
  file systems,
  menus systems, asymmetric communication protocols
  (e.g.,
  \cite{bose01}, \cite{czy03},
  \cite{jacobs12}, \cite{laber06}, \cite{perk99}).

 \begin{center}
\begin{picture}(71,36)

\put(05.5,32){\makebox(0,0)[bl]{Table 8. Hotlink assignment
problems}}

\put(00,0){\line(1,0){71}} \put(00,23){\line(1,0){71}}
\put(00,30){\line(1,0){71}}

\put(00,0){\line(0,1){30}} \put(55,0){\line(0,1){30}}
\put(71,0){\line(0,1){30}}

\put(01,25){\makebox(0,0)[bl]{Problem}}
\put(56,25){\makebox(0,0)[bl]{Source}}

%------------------------------------

\put(01,18){\makebox(0,0)[bl]{1.Basic hotlink assignment}}
%\put(56,18){\makebox(0,0)[bl]{}}
\put(56,18){\makebox(0,0)[bl]{\cite{bose01},\cite{perk99}}}

%------------------------------------

\put(01,14){\makebox(0,0)[bl]{2.Single  hotlink assignment}}
%\put(56,34){\makebox(0,0)[bl]{}}
\put(56,14){\makebox(0,0)[bl]{\cite{dou10},\cite{kra01}}}

%------------------------------------

\put(01,10){\makebox(0,0)[bl]{3.Hotlinks only for leafs}}
%\put(56,10){\makebox(0,0)[bl]{}}
\put(56,10){\makebox(0,0)[bl]{\cite{jacobs12}}}

%------------------------------------

\put(01,6){\makebox(0,0)[bl]{4.Multiple hotlink assignment}}
%\put(56,6){\makebox(0,0)[bl]{}}
\put(56,6){\makebox(0,0)[bl]{\cite{dou10},\cite{fuhr01}}}

%-------------
\put(01,02){\makebox(0,0)[bl]{5.Dynamic hotlink assignment}}
%\put(56,02){\makebox(0,0)[bl]{}}
\put(56,02){\makebox(0,0)[bl]{\cite{dou08}}}

\end{picture}
\end{center}

 Generally, ``hotlink assignment'' problem
  is  a special case of
  ``graph augmentation problem''
  (e.g., \cite{esw76},\cite{khu97}).
 Here the goal is
 to modify an initial graph (e.g., by edges)
 such that the augmented graph will by satisfied
 some requirements (e.g., as increasing the connectivity).

\subsection{Scheme for Transformation of Tree to Steiner
Tree}

 Here a transformation of a tree
 \(T = (A,E)\)
  into Steiner tree
 \(S = (A',E')\)
  is considered
 as addition of Steiner points into an initial tree
 (or a preliminary built spanning tree)
 while taking into account the following:
 ``cost'' (required resource) of each Steiner point,
 ``profit'' of each Steiner point,
 total resource constraint (i.e., total ``cost'' of the selected Steiner points)
 \cite{levzam11}.
 A simplest case is considered when Steiner points for triangles are only
 examined.
 Evidently, vector-like
``cost'' and ``profit'' can be used as well.
 The solving scheme is the following:

~~

 {\it Stage 1.} Identification
 (e.g., expert judgment, clustering)
  of \(m\) regions
 (clusters, groups of neighbor nodes)
 in the initial tree \(T\)
 for possible addition of Steiner points.

 {\it Stage 2.} Generation of possible Steiner points (candidates)
 and their attributes (i.e., cost of addition, ``profit'').

 {\it Stage 3.} Formulation of multiple choice problem for
 selection of the best additional Steiner points
 while taking into account resource constraint(s):
 \[ \max~ \sum_{i=1}^{m}  \sum_{i=1}^{q_{i}}    ~ c_{ij} x_{ij}
  ~~~ s.t. ~~ \sum_{i=1}^{m} \sum_{j=1}^{q_{i}}  a_{ij} x_{ij} \leq b,
 ~~ \sum_{j=1}^{q_{i}} x_{ij} =  1,
 ~~ x_{ij} \in \{0, 1\};\]
 where
  \(i\) is the index of region (\(i = \overline{1,m}\)),
 \(q_{i}\) is the number of versions for addition of Steiner
 points in region \(i=\overline{1,m}\),
  \(j\) is the index of version for addition of Steiner points in region
 (\( j = \overline{1,q_{i}} \)  ),
 \(x_{ij}\) is binary variable that equals \(1\) if
 version \(j\) in region \(i\) is selected,
 \(b\) is a total constraint for the required resources
 (i.e., a total ``cost'').

 {\it Stage 4.} Solving the  multiple choice problem
 to obtain the resultant Steiner tree \(S\).

~

 A numerical illustrative example illustrates the scheme.
 Initial tree is (Fig. 37):
 \(T = (A,E)\), \(A = \{ 1,2,3,4,5,6,7,8,9,10,11 \}\).
 Four regions are defined (Fig. 37):
 region \(1\):  \(\{ 1,2,3,4 \}\);
 region \(2\):  \(\{ 4,6,7 \}\);
 region \(3\):  \(\{ 4,5,6,9,11 \}\);
 and
 region \(4\):  \(\{ 7,8,10 \}\).
 The considered Steiner points are the following (Fig. 38):
 region \(1\):  \(s_{11},s_{12}\);
 region \(2\):  \(s_{21} \);
 region \(3\):  \(s_{31},s_{32}\);
 and
 region \(4\):  \(s_{41}\).
 Table 9 contains  for multiple choice problem:
 binary variables and corresponding attributes
 (required resource as ``cost'', possible ``profit'').
  Thus, the problem is:
 \[ \max~ \sum_{i=1}^{4}  \sum_{i=1}^{q_{i}}    ~ c_{ij} x_{ij}
 ~~~~ s.t. ~~ \sum_{i=1}^{4} \sum_{j=1}^{q_{i}}
 a_{ij} x_{ij} \leq b,
 ~~ \sum_{j=1}^{q_{i}} x_{ij} =  1,
 ~~ x_{ij} \in \{0, 1\}.\]

\begin{center}
\begin{picture}(45,53)
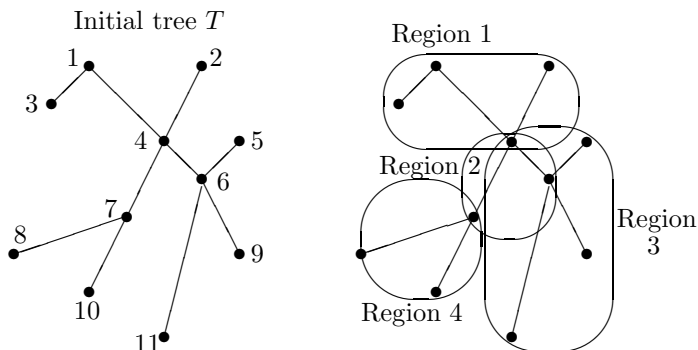
\put(012,00){\makebox(0,0)[bl]{Fig. 37. Initial tree and
 regions (clusters)}}

\put(08,48){\makebox(0,0)[bl]{Initial tree \(T\)}}

%--3(10)
\put(01.5,37){\makebox(0,0)[bl]{\(3\)}} \put(05,38){\circle*{1.3}}

%=======================10-9
\put(05,38){\line(1,1){5}}

%=======================9-8
%\put(10,43){\line(1,0){15}}

%=======================9-4
%\put(10,43){\line(1,-4){5}}

%=======================9-7
\put(10,43){\line(1,-1){10}}

%--1(9)
\put(07,43){\makebox(0,0)[bl]{\(1\)}} \put(10,43){\circle*{1.3}}

%--2(8)
\put(26,43){\makebox(0,0)[bl]{\(2\)}} \put(25,43){\circle*{1.3}}

%=======================8-7
\put(25,43){\line(-1,-2){5}}

%=======================8-5
%\put(25,43){\line(0,-1){15}}

%--4(7)
\put(16,32){\makebox(0,0)[bl]{\(4\)}} \put(20,33){\circle*{1.3}}

%=======================7-4
\put(20,33){\line(-1,-2){5}}

%=======================7-5
\put(20,33){\line(1,-1){5}}

%--9(6)
\put(31.5,17){\makebox(0,0)[bl]{\(9\)}} \put(30,18){\circle*{1.3}}

%--5(6')
\put(31.5,32){\makebox(0,0)[bl]{\(5\)}} \put(30,33){\circle*{1.3}}

%=======================6'-X
\put(30,33){\line(-1,-1){5}}

%=======================5-6
\put(25,28){\line(1,-2){5}}

%--6(5)
\put(27,26.5){\makebox(0,0)[bl]{\(6\)}} \put(25,28){\circle*{1.3}}

%%%10:(5,38)
%=======================4-10
%\put(15,23){\line(-2,3){10}}

%=======================4-5
%\put(15,23){\line(2,1){10}}

%--(7)4
\put(012,23){\makebox(0,0)[bl]{\(7\)}} \put(15,23){\circle*{1.3}}

%=======================4-3
\put(15,23){\line(-3,-1){15}}

%=======================4-2
\put(15,23){\line(-1,-2){5}}

%=======================4-1
%\put(15,23){\line(1,-3){5}}

%9:10,25
%=======================3-9
%\put(00,18){\line(1,2){5}} \put(05,28){\line(1,3){5}}

%--8(3)
\put(00,19.5){\makebox(0,0)[bl]{\(8\)}} \put(00,18){\circle*{1.3}}

%=======================2-6
%\put(10,13){\line(4,1){20}}

%=======================2-3
%\put(10,13){\line(-2,1){10}}

%--10(2)
\put(08,9.4){\makebox(0,0)[bl]{\(10\)}} \put(10,13){\circle*{1.3}}

%=======================1-5
\put(20,7){\line(1,4){5}}

%--11(1)
\put(016,5){\makebox(0,0)[bl]{\(11\)}} \put(20,7){\circle*{1.3}}

\end{picture}
%\end{center}
%%%%%%%%%%%%%%%%%%%%%%%%%%%%%%%%%%%%%%%%%%%%%%%%% Regions
%\begin{center}
\begin{picture}(45,53)

\put(04,45.6){\makebox(0,0)[bl]{Region \(1\)}}
%---------------------------------------------Region 1
\put(16,38.3){\oval(26,12.6)}

%--3(10) \put(01.5,37){\makebox(0,0)[bl]{\(3\)}}
\put(05,38){\circle*{1.3}}

%=======================10-9
\put(05,38){\line(1,1){5}}

%=======================9-8
%\put(10,43){\line(1,0){15}}

%=======================9-4
%\put(10,43){\line(1,-4){5}}

%=======================9-7
\put(10,43){\line(1,-1){10}}

%--1(9)  \put(07,43){\makebox(0,0)[bl]{\(1\)}}
\put(10,43){\circle*{1.3}}

%--2(8)  \put(26,43){\makebox(0,0)[bl]{\(2\)}}
\put(25,43){\circle*{1.3}}

%=======================8-7
\put(25,43){\line(-1,-2){5}}

%=======================8-5
%\put(25,43){\line(0,-1){15}}

\put(02.4,28){\makebox(0,0)[bl]{Region \(2\)}}
%------------------------------------------Region 2
\put(19.7,27){\oval(12.5,14)}

%--4(7)  \put(16,32){\makebox(0,0)[bl]{\(4\)}}
\put(20,33){\circle*{1.3}}

%=======================7-4
\put(20,33){\line(-1,-2){5}}

%=======================7-5
\put(20,33){\line(1,-1){5}}

\put(034,21){\makebox(0,0)[bl]{Region }}
\put(038,18){\makebox(0,0)[bl]{\(3\)}}
%----------------------------------------------Region 3
\put(25,20){\oval(17,30)}

%--9(6) \put(31.5,17){\makebox(0,0)[bl]{\(9\)}}
\put(30,18){\circle*{1.3}}

%--5(6')  \put(31.5,32){\makebox(0,0)[bl]{\(5\)}}
\put(30,33){\circle*{1.3}}

%=======================6'-X
\put(30,33){\line(-1,-1){5}}

%=======================5-6
\put(25,28){\line(1,-2){5}}

%--6(5) \put(27,26.5){\makebox(0,0)[bl]{\(6\)}}
\put(25,28){\circle*{1.3}}

%%%10:(5,38)
%=======================4-10
%\put(15,23){\line(-2,3){10}}

%=======================4-5
%\put(15,23){\line(2,1){10}}

%--7(4) \put(012,23){\makebox(0,0)[bl]{\(7\)}}
\put(15,23){\circle*{1.3}}

%=======================4-3
\put(15,23){\line(-3,-1){15}}

%=======================4-2
\put(15,23){\line(-1,-2){5}}

%=======================4-1
%\put(15,23){\line(1,-3){5}}

%9:10,25
%=======================3-9
%\put(00,18){\line(1,2){5}} \put(05,28){\line(1,3){5}}

%--8(3) \put(00,19.5){\makebox(0,0)[bl]{\(8\)}}
\put(00,18){\circle*{1.3}}

%=======================2-6
%\put(10,13){\line(4,1){20}}

%=======================2-3
%\put(10,13){\line(-2,1){10}}

\put(00,8.5){\makebox(0,0)[bl]{Region \(4\)}}
%----------------------------------------Region 4
\put(08,20){\oval(16,16)}

%--10(2) \put(08,9.4){\makebox(0,0)[bl]{\(10\)}}
\put(10,13){\circle*{1.3}}

%=======================1-5
\put(20,7){\line(1,4){5}}

%--11(1)\put(016,5){\makebox(0,0)[bl]{\(11\)}}
 \put(20,7){\circle*{1.3}}

\end{picture}
\end{center}

\begin{center}
\begin{picture}(45,53)
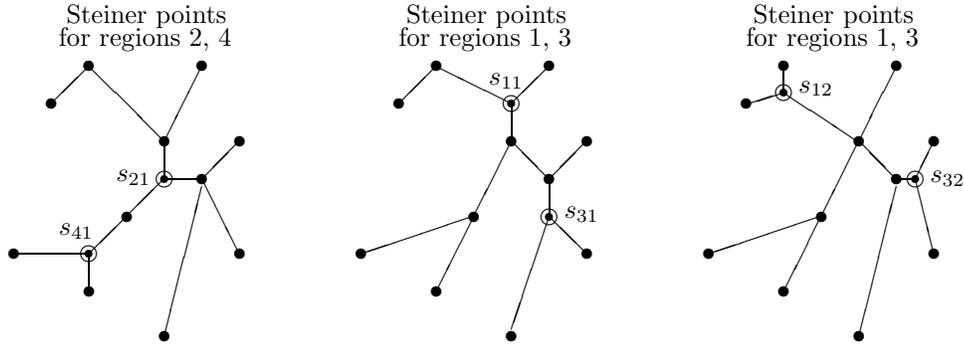
\put(026,00){\makebox(0,0)[bl]{Fig. 38. Additional Steiner points
 for regions}}

\put(07,48){\makebox(0,0)[bl]{Steiner points}}
\put(06,45){\makebox(0,0)[bl]{for regions 2, 4}}

%--10
\put(05,38){\circle*{1.3}}

%=======================10-9
\put(05,38){\line(1,1){5}}

%=======================9-8
%\put(10,43){\line(1,0){15}}

%=======================9-4
%\put(10,43){\line(1,-4){5}}

%=======================9-7
\put(10,43){\line(1,-1){10}}

%--9(6)
\put(10,43){\circle*{1.3}}

%--8
\put(25,43){\circle*{1.3}}

%=======================8-7
\put(25,43){\line(-1,-2){5}}

%=======================8-5
%\put(25,43){\line(0,-1){15}}

%--7
\put(20,33){\circle*{1.3}}

%=======================7-4
%\put(20,33){\line(-1,-2){5}}

%=======================7-5
%\put(20,33){\line(1,-1){5}}

%--6
\put(30,18){\circle*{1.3}}

%--6'
\put(30,33){\circle*{1.3}}

%=======================6'-X
\put(30,33){\line(-1,-1){5}}

%=======================5-6
\put(25,28){\line(1,-2){5}}

%--5
\put(25,28){\circle*{1.3}}

%%%%%%%%%%%%%%%%%%%%%%%%%%%%%%%%%%%%%%%%%%% Steiner point 1-4-5
\put(20,28){\circle*{1.0}} \put(20,28){\circle{2.0}}

\put(13.6,27.2){\makebox(0,0)[bl]{\(s_{21}\)}}

%=======================1-4-5->
\put(20,28){\line(0,1){5}} \put(20,28){\line(1,0){5}}
\put(20,28){\line(-1,-1){5}}
%%%%%%%%%%%%%%%%%%%%%%%%%%%%%%%%%%%%%%%%%%%%%%%%

%--4
\put(15,23){\circle*{1.3}}

%=======================4-3
%\put(15,23){\line(-3,-1){15}}

%=======================4-2
%\put(15,23){\line(-1,-2){5}}

%=======================4-1
%\put(15,23){\line(1,-3){5}}

%%%%%%%%%%%%%%%%%%%%%%%%%%%%%%%%%%%%%%%%%%% Steiner point 2-3-4
\put(10,18){\circle*{1.0}} \put(10,18){\circle{2.0}}

\put(06,20){\makebox(0,0)[bl]{\(s_{41}\)}}

%=======================2-3-4-> 9
\put(10,18){\line(-1,0){10}} \put(10,18){\line(0,-1){5}}
\put(10,18){\line(1,1){5}}
%%%%%%%%%%%%%%%%%%%%%%%%%%%%%%%%%%%%%%%%%%%%%%%%

%--3
\put(00,18){\circle*{1.3}}

%--2
\put(10,13){\circle*{1.3}}

%=======================1-5
\put(20,7){\line(1,4){5}}

%--1
 \put(20,7){\circle*{1.3}}

\end{picture}
%\end{center}
%
%\begin{center}
%
\begin{picture}(45,53)
\put(06,48){\makebox(0,0)[bl]{Steiner points}}
\put(05,45){\makebox(0,0)[bl]{for regions 1, 3}}

%--10
\put(05,38){\circle*{1.3}}

%=======================10-9
\put(05,38){\line(1,1){5}}

%--9
\put(10,43){\circle*{1.3}}

%--8
\put(25,43){\circle*{1.3}}

%%%%%%%%%%%%%%%%%%%%%%%%%%%%%%%%%%%%%%%%%%%Steiner point 7-8-9
\put(20,38){\circle*{1.0}} \put(20,38){\circle{2.0}}

\put(17,40){\makebox(0,0)[bl]{\(s_{11}\)}}

%=======================7-8-9-> 9
\put(20,38){\line(-2,1){10}} \put(20,38){\line(1,1){5}}
\put(20,38){\line(0,-1){5}}
%%%%%%%%%%%%%%%%%%%%%%%%%%%%%%%%%%%%%%%%%%%%%%%%

%--7
\put(20,33){\circle*{1.3}}

%=======================7-4
\put(20,33){\line(-1,-2){5}}

%=======================7-5
\put(20,33){\line(1,-1){5}}

%--6
\put(30,18){\circle*{1.3}}

%--6'
\put(30,33){\circle*{1.3}}

%=======================6'-X
\put(30,33){\line(-1,-1){5}}

%%%%%%%%%%%%%%%%%%%%%%%%%%%%%%%%%%%%%%%%%%% Steiner point 1-5-6
\put(25,23){\circle*{1.0}} \put(25,23){\circle{2.0}}

\put(27,22.3){\makebox(0,0)[bl]{\(s_{31}\)}}

%=======================1-5-6->
\put(25,23){\line(0,1){5}} \put(25,23){\line(1,-1){5}}
\put(25,23){\line(-1,-3){5}}
%%%%%%%%%%%%%%%%%%%%%%%%%%%%%%%%%%%%%%%%%%%%%%%%

%--5
\put(25,28){\circle*{1.3}}

%--4
\put(15,23){\circle*{1.3}}

%=======================4-3
\put(15,23){\line(-3,-1){15}}

%=======================4-2
\put(15,23){\line(-1,-2){5}}

%--3
\put(00,18){\circle*{1.3}}

%--2
\put(10,13){\circle*{1.3}}

%--1
 \put(20,7){\circle*{1.3}}

\end{picture}
%\end{center}
%
\begin{picture}(30,53)
\put(06,48){\makebox(0,0)[bl]{Steiner points}}
\put(05,45){\makebox(0,0)[bl]{for regions 1, 3}}

%--10
\put(05,38){\circle*{1.3}}

%--9
\put(10,43){\circle*{1.3}}

%--8
\put(25,43){\circle*{1.3}}

%%%%%%%%%%%%%%%%%%%%%%%%%%%%%%%%%%%%%%%%%%%Steiner point 1-3-4 (9-10-7)
\put(10,39.5){\circle*{1.0}} \put(10,39.5){\circle{2.0}}

\put(12,39){\makebox(0,0)[bl]{\(s_{12}\)}}

%=======================7-8-9-> 9
%\put(20,38){\line(-2,1){10}} \put(20,38){\line(1,1){5}}
\put(10,39.5){\line(-3,-1){5}} \put(10,39.5){\line(0,1){3.5}}
\put(10,39.5){\line(3,-2){10}}

%=======================8-7
\put(25,43){\line(-1,-2){5}}

%--7
\put(20,33){\circle*{1.3}}

%=======================7-4
\put(20,33){\line(-1,-2){5}}

%=======================7-5
\put(20,33){\line(1,-1){5}}

%--6
\put(30,18){\circle*{1.3}}

%--6'
\put(30,33){\circle*{1.3}}

%%%%%%%%%%%%%%%%%%%%%%%%%%%%%%%%%%%%%%%%%%% Steiner point 5-6-9 (5-6-6')
\put(27.5,28){\circle*{1.0}} \put(27.5,28){\circle{2.0}}

\put(29.5,27){\makebox(0,0)[bl]{\(s_{32}\)}}

%=======================1-5-6->
\put(27.5,28){\line(-1,0){2.5}} \put(27.5,28){\line(1,2){2.5}}
\put(27.5,28){\line(1,-4){2.5}}
%%%%%%%%%%%%%%%%%%%%%%%%%%%%%%%%%%%%%%%%%%%%%%%%

%--5
\put(25,28){\circle*{1.3}}

%--4
\put(15,23){\circle*{1.3}}

%=======================4-3
\put(15,23){\line(-3,-1){15}}

%=======================4-2
\put(15,23){\line(-1,-2){5}}

%--3
\put(00,18){\circle*{1.3}}

%=======================2-6
%\put(10,13){\line(4,1){20}}

%=======================2-3
%\put(10,13){\line(-2,1){10}}

%--2
\put(10,13){\circle*{1.3}}

%=======================1-5
\put(20,7){\line(1,4){5}}

%--1
 \put(20,7){\circle*{1.3}}

\end{picture}
\end{center}

\begin{center}
\begin{picture}(68,67)
\put(01.5,62){\makebox(0,0)[bl]{Table 9. Data for multiple
 choice problem}}

\put(00,0){\line(1,0){68}} \put(00,49){\line(1,0){68}}
\put(00,60){\line(1,0){68}}

%--

\put(00,00){\line(0,1){60}} \put(16,00){\line(0,1){60}}
\put(30,00){\line(0,1){60}} \put(44,00){\line(0,1){60}}
\put(55,00){\line(0,1){60}} \put(68,00){\line(0,1){60}}

%--

\put(01,54.6){\makebox(0,0)[bl]{Region}}

\put(17,54.6){\makebox(0,0)[bl]{Binary}}
\put(17,51.4){\makebox(0,0)[bl]{variable}}

\put(32,55){\makebox(0,0)[bl]{Steiner}}
\put(32,51){\makebox(0,0)[bl]{point}}

\put(44.5,55){\makebox(0,0)[bl]{``Cost''}}
\put(47,51){\makebox(0,0)[bl]{\(c_{ij}\)}}

\put(55.5,55){\makebox(0,0)[bl]{``Profit''}}
\put(57,51){\makebox(0,0)[bl]{\(a_{ij}\)}}

%%%%%%%%%%%%%%%%%%%%%%%

\put(01,43.5){\makebox(0,0)[bl]{Region \(1\)}}

\put(21,44){\makebox(0,0)[bl]{\(x_{11}\)}}
\put(33,44){\makebox(0,0)[bl]{None}}
\put(47,44){\makebox(0,0)[bl]{\(0.0\)}}
\put(59,44){\makebox(0,0)[bl]{\(0.0\)}}

\put(21,40){\makebox(0,0)[bl]{\(x_{12}\)}}
\put(35,40){\makebox(0,0)[bl]{\(s_{11}\)}}
\put(47,40){\makebox(0,0)[bl]{\(3.1\)}}
\put(59,40){\makebox(0,0)[bl]{\(1.5\)}}

\put(21,36){\makebox(0,0)[bl]{\(x_{13}\)}}
\put(35,36){\makebox(0,0)[bl]{\(s_{12}\)}}
\put(47,36){\makebox(0,0)[bl]{\(1.2\)}}
\put(59,36){\makebox(0,0)[bl]{\(1.4\)}}

%--

\put(00,34){\line(1,0){68}}

\put(01,29.5){\makebox(0,0)[bl]{Region \(2\)}}

\put(21,30){\makebox(0,0)[bl]{\(x_{21}\)}}
\put(33,30){\makebox(0,0)[bl]{None}}
\put(47,30){\makebox(0,0)[bl]{\(0.0\)}}
\put(59,30){\makebox(0,0)[bl]{\(0.0\)}}

\put(21,26){\makebox(0,0)[bl]{\(x_{22}\)}}
\put(35,26){\makebox(0,0)[bl]{\(s_{21}\)}}
\put(47,26){\makebox(0,0)[bl]{\(2.0\)}}
\put(59,26){\makebox(0,0)[bl]{\(1.3\)}}

%--

\put(00,24){\line(1,0){68}}

\put(01,19.5){\makebox(0,0)[bl]{Region \(3\)}}

\put(21,20){\makebox(0,0)[bl]{\(x_{31}\)}}
\put(33,20){\makebox(0,0)[bl]{None}}
\put(47,20){\makebox(0,0)[bl]{\(0.0\)}}
\put(59,20){\makebox(0,0)[bl]{\(0.0\)}}

\put(21,16){\makebox(0,0)[bl]{\(x_{32}\)}}
\put(35,16){\makebox(0,0)[bl]{\(s_{31}\)}}
\put(47,16){\makebox(0,0)[bl]{\(2.4\)}}
\put(59,16){\makebox(0,0)[bl]{\(1.4\)}}

\put(21,12){\makebox(0,0)[bl]{\(x_{33}\)}}
\put(35,12){\makebox(0,0)[bl]{\(s_{32}\)}}
\put(47,12){\makebox(0,0)[bl]{\(1.8\)}}
\put(59,12){\makebox(0,0)[bl]{\(1.3\)}}

%\put(21,12){\makebox(0,0)[bl]{\(x_{34}\)}}
%\put(31,12){\makebox(0,0)[bl]{\(s_{31} \& s_{32}\)}}
%\put(47,12){\makebox(0,0)[bl]{\(4.0\)}}
%\put(57,12){\makebox(0,0)[bl]{\(2.7\)}}

%--

\put(00,10){\line(1,0){68}}

\put(01,05.5){\makebox(0,0)[bl]{Region \(4\)}}

\put(21,06){\makebox(0,0)[bl]{\(x_{41}\)}}
\put(33,6){\makebox(0,0)[bl]{None}}
\put(47,6){\makebox(0,0)[bl]{\(0.0\)}}
\put(59,6){\makebox(0,0)[bl]{\(0.0\)}}

\put(21,02){\makebox(0,0)[bl]{\(x_{42}\)}}
\put(35,02){\makebox(0,0)[bl]{\(s_{41}\)}}
\put(47,02){\makebox(0,0)[bl]{\(1.5\)}}
\put(59,02){\makebox(0,0)[bl]{\(1.2\)}}

\end{picture}
\end{center}

 Some obtained solutions (i.e., as set of additional Steiner
 points) are the following (a simple greedy heuristic was used)
 (Fig. 39):

 (1) \(b_{1} = 2.9\):~
 \( \overline{x}_{b_{1}}\):~ \( x_{12}=1\), \(x_{21}=1\), \(x_{32}=1\), \(x_{41}=1\),
 Steiner points \(Z_{b_{1}} =  \{ s_{11},s_{31} \} \),

 total (additive) ``profit'' \( \overline{c} = 5.5 \);

 (2) \(b_{2} = 4.2\):~
 \( \overline{x}_{b_{2}} \):  \(x_{12}=1\), \(x_{22}=1\), \(x_{32}=1\),
 \(x_{42}=0\);
 Steiner points \(Z_{b_{2}} = \{ s_{11},s_{21},s_{31} \} \),

 total (additive) ``profit'' \( \overline{c} = 7.5 \);

 (3) \(b_{3} = 5.4\):~
 \( \overline{x}_{b_{3}}\):   \(x_{12}=1\), \(x_{22}=1\), \(x_{32}=1\), \(x_{42}=1\),
 Steiner points \(Z_{b_{3}} =  \{ s_{11},s_{21},s_{31},s_{41} \} \),

 total (additive) ``profit'' \( \overline{c} = 9.0\).

 Note, more complicated problem can be examined, for example,  while
 taking into account the following:
 (i) several additional Steiner points in the same
 region,
 (ii) compatibility of points additions in neighbor regions.

\begin{center}
\begin{picture}(45,53)
\put(47,00){\makebox(0,0)[bl]{Fig. 39. Solutions}}

\put(07.5,48){\makebox(0,0)[bl]{Solution \(S_{b_{1}}\):}}
\put(03,44){\makebox(0,0)[bl]{(2 Steiner points)}}

%--10
\put(05,38){\circle*{1.3}}

%=======================10-9
\put(05,38){\line(1,1){5}}

%--9
\put(10,43){\circle*{1.3}}

%--8
\put(25,43){\circle*{1.3}}

%%%%%%%%%%%%%%%%%%%%%%%%%%%%%%%%%%%%%%%%%%%Steiner point 7-8-9
\put(20,38){\circle*{1.0}} \put(20,38){\circle{2.0}}

\put(17,40){\makebox(0,0)[bl]{\(s_{11}\)}}

%=======================7-8-9-> 9
\put(20,38){\line(-2,1){10}} \put(20,38){\line(1,1){5}}
\put(20,38){\line(0,-1){5}}
%%%%%%%%%%%%%%%%%%%%%%%%%%%%%%%%%%%%%%%%%%%%%%%%

%--7
\put(20,33){\circle*{1.3}}

%=======================7-4
\put(20,33){\line(-1,-2){5}}

%=======================7-5
\put(20,33){\line(1,-1){5}}

%--6
\put(30,18){\circle*{1.3}}

%--6'
\put(30,33){\circle*{1.3}}

%=======================6'-X
\put(30,33){\line(-1,-1){5}}

%%%%%%%%%%%%%%%%%%%%%%%%%%%%%%%%%%%%%%%%%%% Steiner point 1-5-6
\put(25,23){\circle*{1.0}} \put(25,23){\circle{2.0}}

\put(27,22.3){\makebox(0,0)[bl]{\(s_{31}\)}}

%=======================1-5-6->
\put(25,23){\line(0,1){5}} \put(25,23){\line(1,-1){5}}
\put(25,23){\line(-1,-3){5}}
%%%%%%%%%%%%%%%%%%%%%%%%%%%%%%%%%%%%%%%%%%%%%%%%

%--5
\put(25,28){\circle*{1.3}}

%--4
\put(15,23){\circle*{1.3}}

%=======================4-3
\put(15,23){\line(-3,-1){15}}

%=======================4-2
\put(15,23){\line(-1,-2){5}}

%--3
\put(00,18){\circle*{1.3}}

%=======================2-6
%\put(10,13){\line(4,1){20}}

%=======================2-3
%\put(10,13){\line(-2,1){10}}

%--2
\put(10,13){\circle*{1.3}}

%=======================1-5
%\put(20,7){\line(1,4){5}}

%--1
 \put(20,7){\circle*{1.3}}

\end{picture}
%\end{center}
%
%\begin{center}
\begin{picture}(45,53)

\put(07.5,48){\makebox(0,0)[bl]{Solution \(S_{b_{2}}\):}}
\put(03,44){\makebox(0,0)[bl]{(3 Steiner points)}}

%--10
\put(05,38){\circle*{1.3}}

%=======================10-9
\put(05,38){\line(1,1){5}}

%=======================9-8
%\put(10,43){\line(1,0){15}}

%=======================9-4
%\put(10,43){\line(1,-4){5}}

%=======================9-7
%\put(10,43){\line(1,-1){10}}

%--9
\put(10,43){\circle*{1.3}}

%--8
\put(25,43){\circle*{1.3}}

%%%%%%%%%%%%%%%%%%%%%%%%%%%%%%%%%%%%%%%%%%%Steiner point 7-8-9
\put(20,38){\circle*{1.0}} \put(20,38){\circle{2.0}}

\put(17,40){\makebox(0,0)[bl]{\(s_{11}\)}}

%=======================7-8-9-> 9
\put(20,38){\line(-2,1){10}} \put(20,38){\line(1,1){5}}
\put(20,38){\line(0,-1){5}}
%%%%%%%%%%%%%%%%%%%%%%%%%%%%%%%%%%%%%%%%%%%%%%%%

%--7
\put(20,33){\circle*{1.3}}

%--6
\put(30,18){\circle*{1.3}}

%--6'
\put(30,33){\circle*{1.3}}

%=======================6'-X
\put(30,33){\line(-1,-1){5}}

%%%%%%%%%%%%%%%%%%%%%%%%%%%%%%%%%%%%%%%%%%% Steiner point 1-5-6
\put(25,23){\circle*{1.0}} \put(25,23){\circle{2.0}}

\put(27,22.3){\makebox(0,0)[bl]{\(s_{31}\)}}

%=======================1-5-6->
\put(25,23){\line(0,1){5}} \put(25,23){\line(1,-1){5}}
\put(25,23){\line(-1,-3){5}}
%%%%%%%%%%%%%%%%%%%%%%%%%%%%%%%%%%%%%%%%%%%%%%%%

%--5
\put(25,28){\circle*{1.3}}

%%%%%%%%%%%%%%%%%%%%%%%%%%%%%%%%%%%%%%%%%%% Steiner point 1-4-5
\put(20,28){\circle*{1.0}} \put(20,28){\circle{2.0}}

\put(13.6,27.2){\makebox(0,0)[bl]{\(s_{21}\)}}

%=======================1-4-5->
\put(20,28){\line(0,1){5}} \put(20,28){\line(1,0){5}}
\put(20,28){\line(-1,-1){5}}
%%%%%%%%%%%%%%%%%%%%%%%%%%%%%%%%%%%%%%%%%%%%%%%%

%--4
\put(15,23){\circle*{1.3}}

%=======================4-3
\put(15,23){\line(-3,-1){15}}

%=======================4-2
\put(15,23){\line(-1,-2){5}}

%--3
\put(00,18){\circle*{1.3}}

%--2
\put(10,13){\circle*{1.3}}

%=======================1-5
%\put(20,7){\line(1,4){5}}

%--1
 \put(20,7){\circle*{1.3}}

\end{picture}
%\end{center}
%
%\begin{center}
\begin{picture}(45,53)

\put(07.5,48){\makebox(0,0)[bl]{Solution \(S_{b_{3}}\):}}
\put(03,44){\makebox(0,0)[bl]{(4 Steiner points)}}

%+++++++++++++++++++++++++++++++++++++++++++++++++++++++++++++++

\put(13.6,27.2){\makebox(0,0)[bl]{\(s_{21}\)}}
\put(06,20){\makebox(0,0)[bl]{\(s_{41}\)}}
\put(27,22.3){\makebox(0,0)[bl]{\(s_{31}\)}}
\put(17,40){\makebox(0,0)[bl]{\(s_{11}\)}}

%--10
\put(05,38){\circle*{1.3}}

%=======================10-9
\put(05,38){\line(1,1){5}}

%--9
\put(10,43){\circle*{1.3}}

%--8
\put(25,43){\circle*{1.3}}

%%%%%%%%%%%%%%%%%%%%%%%%%%%%%%%%%%%%%%%%%%%Steiner point 7-8-9
\put(20,38){\circle*{1.0}} \put(20,38){\circle{2.0}}

%=======================7-8-9-> 9
\put(20,38){\line(-2,1){10}} \put(20,38){\line(1,1){5}}
\put(20,38){\line(0,-1){5}}
%%%%%%%%%%%%%%%%%%%%%%%%%%%%%%%%%%%%%%%%%%%%%%%%

%--7
\put(20,33){\circle*{1.3}}

%--6
\put(30,18){\circle*{1.3}}

%--6'
\put(30,33){\circle*{1.3}}

%=======================6'-X
\put(30,33){\line(-1,-1){5}}

%%%%%%%%%%%%%%%%%%%%%%%%%%%%%%%%%%%%%%%%%%% Steiner point 1-5-6
\put(25,23){\circle*{1.0}} \put(25,23){\circle{2.0}}

%=======================1-5-6->
\put(25,23){\line(0,1){5}} \put(25,23){\line(1,-1){5}}
\put(25,23){\line(-1,-3){5}}
%%%%%%%%%%%%%%%%%%%%%%%%%%%%%%%%%%%%%%%%%%%%%%%%

%--5
\put(25,28){\circle*{1.3}}

%%%%%%%%%%%%%%%%%%%%%%%%%%%%%%%%%%%%%%%%%%% Steiner point 1-4-5
\put(20,28){\circle*{1.0}} \put(20,28){\circle{2.0}}

%=======================1-4-5->
\put(20,28){\line(0,1){5}} \put(20,28){\line(1,0){5}}
\put(20,28){\line(-1,-1){5}}
%%%%%%%%%%%%%%%%%%%%%%%%%%%%%%%%%%%%%%%%%%%%%%%%

%%%%%%%%%%%%%%%%%%%%%%%%%%%%%%%%%%%%%%%%%%%Steiner point 4-7-5
\put(20,38){\circle*{1.0}} \put(20,38){\circle{2.0}}

%=======================4-7-5->
\put(20,38){\line(-2,1){10}} \put(20,38){\line(1,1){5}}
\put(20,38){\line(0,-1){5}}
%%%%%%%%%%%%%%%%%%%%%%%%%%%%%%%%%%%%%%%%%%%%%%%%

%--4
\put(15,23){\circle*{1.3}}

%%%%%%%%%%%%%%%%%%%%%%%%%%%%%%%%%%%%%%%%%%% Steiner point 2-3-4
\put(10,18){\circle*{1.0}} \put(10,18){\circle{2.0}}

%=======================2-3-4-> 9
\put(10,18){\line(-1,0){10}} \put(10,18){\line(0,-1){5}}
\put(10,18){\line(1,1){5}}
%%%%%%%%%%%%%%%%%%%%%%%%%%%%%%%%%%%%%%%%%%%%%%%%

%--3
\put(00,18){\circle*{1.3}}

%--2
\put(10,13){\circle*{1.3}}

%--1
 \put(20,7){\circle*{1.3}}

\end{picture}
\end{center}

\subsection{Towards Restructuring Problems}

 In recent several years,
 a special class of combinatorial optimization
 problems as  ``reoptimization''
 has been intensively
 studied for several well-known problems:
 (a) minimum spanning tree problem
 \cite{boria10},
 (b) traveling salesman problems
  (\cite{archetti03}, \cite{aus09}),
 (c) Steiner tree problems
 (\cite{bilo08},\cite{esc09}),
 (d)  covering problems
 \cite{bilo08a}.

 In general, the reoptimization problem is formulated
 as follows:

~~

 {\it Given}:
 (i) an instance of the combinatorial problem over a graph
 and corresponding optimal solution,
 (ii) some ``small'' perturbations
 (i.e., modifications) on this instance
 (e.g., node-insertion, node-deletion).

 {\it Question}:~~
 {\it Is it possible to compute a new good (optimal or near-optimal)
 solution
 subject to minor modifications?}

~~

 A survey of complexity issues for reoptimization problems
 is presented in \cite{bok08}.
 Mainly, the problems belong
 to class of NP-hard problems and
 various approximation algorithms have been suggested

 In \cite{lev11restr},
 the author has suggested another approach to modification in
 combinatorial optimization problems as ``restructuring''.
 The approach corresponds to many applied
 reengineering (redesign) problems in  existing modular systems.
 The restructuring process
 is illustrated in Fig. 40 \cite{lev11restr}.
 Here modifications are based on insertion/deletion of elements
 (i.e., elements, nodes, arcs)
 and changes of a structure as well.
 Two main features of the restructuring process are examined:
 (i) a cost of the initial problem solution restructuring
 (i.e., cost of the selected modifications),
 (ii) a closeness the obtained restructured solution to a goal solution.

\begin{center}
\begin{picture}(80,52)

\put(01,00){\makebox(0,0)[bl]{Fig. 40. Illustration for
restructuring process \cite{lev11restr}}}

\put(00,09){\vector(1,0){80}}

\put(00,7.5){\line(0,1){3}} \put(11,7.5){\line(0,1){3}}
\put(69,7.5){\line(0,1){3}}

\put(00,05){\makebox(0,0)[bl]{\(0\)}}
\put(11,05){\makebox(0,0)[bl]{\(\tau_{1}\)}}
\put(67,05){\makebox(0,0)[bl]{\(\tau_{2}\)}}

\put(79,05.3){\makebox(0,0)[bl]{\(t\)}}

%=================================== S1

\put(00,41){\line(1,0){22}} \put(00,51){\line(1,0){22}}
\put(00,41){\line(0,1){10}} \put(22,41){\line(0,1){10}}

\put(0.5,46){\makebox(0,0)[bl]{Requirements}}
\put(0.5,42){\makebox(0,0)[bl]{(for \(\tau_{1}\))}}

%--

\put(11,41){\vector(0,-1){4}}

\put(00,23){\line(1,0){22}} \put(00,37){\line(1,0){22}}
\put(00,23){\line(0,1){14}} \put(22,23){\line(0,1){14}}

\put(0.5,32){\makebox(0,0)[bl]{Optimization}}
\put(0.5,28){\makebox(0,0)[bl]{problem}}
\put(0.5,24){\makebox(0,0)[bl]{(for \(\tau_{1})\)}}

%--

\put(11,23){\vector(0,-1){4}}

\put(11,16){\oval(22,06)}

\put(02,15){\makebox(0,0)[bl]{Solution \(S^{1}\)}}

%%%%%%%%%%%%%%%%%%%%%%%%%%%%%%%%%%%%%%%%%%%%%%%%%%%%%%%%%%

\put(26,14){\line(1,0){28}} \put(26,49){\line(1,0){28}}
\put(26,14){\line(0,1){35}} \put(54,14){\line(0,1){35}}

\put(26.5,14.5){\line(1,0){27}} \put(26.5,48.5){\line(1,0){27}}
\put(26.5,14.5){\line(0,1){34}} \put(53.5,14.5){\line(0,1){34}}

\put(28,44){\makebox(0,0)[bl]{Restructuring:}}
\put(28,40){\makebox(0,0)[bl]{\(S^{1} \Rightarrow S^{*}  \)}}
\put(28,35.5){\makebox(0,0)[bl]{while taking}}
\put(28,32.5){\makebox(0,0)[bl]{into account:}}
\put(28,28){\makebox(0,0)[bl]{(i) \(S^{*}\) is close}}
\put(28,24.4){\makebox(0,0)[bl]{to \(S^{2}\),}}
\put(28,20){\makebox(0,0)[bl]{(ii) change of \(S^{1}\)}}
\put(28,16){\makebox(0,0)[bl]{into \(S^{*}\) is cheap.}}

%================================== S2

%--

\put(58,41){\line(1,0){22}} \put(58,51){\line(1,0){22}}
\put(58,41){\line(0,1){10}} \put(80,41){\line(0,1){10}}

\put(58.5,46){\makebox(0,0)[bl]{Requirements}}
\put(58.5,42){\makebox(0,0)[bl]{(for \(\tau_{2}\))}}

%--

\put(69,41){\vector(0,-1){4}}

\put(58,23){\line(1,0){22}} \put(58,37){\line(1,0){22}}
\put(58,23){\line(0,1){14}} \put(80,23){\line(0,1){14}}

\put(58.5,32){\makebox(0,0)[bl]{Optimization}}
\put(58.5,28){\makebox(0,0)[bl]{problem }}
\put(58.5,24){\makebox(0,0)[bl]{(for \(\tau_{2})\)}}

%--

\put(69,22){\vector(0,-1){4}}

\put(69,16){\oval(22,06)}

\put(59,15){\makebox(0,0)[bl]{Solution \(S^{2}\)}}

%===================================

\end{picture}
\end{center}

 This kind of problems corresponds to
 redesign/reconfiguration (improvement, upgrade)
 of modular systems and
 the situations can be faced
 in complex software, algorithm systems,
 communication networks, computer networks,
 information systems, manufacturing systems,
 constructions,
 etc.
 (e.g,
   \cite{levsib10}, \cite{lev11restr}).
  The optimization problem is solved for two time moments:
 \(\tau_{1}\) and  \(\tau_{2}\) to obtain corresponding solutions
 \(S^{1}\) and \(S^{2}\).
 The examined restructuring problem consists in
 a ``cheap'' transformation (change) of solution \(S^{1}\) to a solution \(S^{*}\) that
 is very close to \(S^{2}\).
 In \cite{lev11restr}, this restructuring approach is described and illustrated for the following combinatorial
 optimization problems:
 knapsack problem,
 multiple choice problem,
 assignment problem,
 spanning tree problems.
 Evidently,
  this restructuring problem may be used for many combinatorial
 optimization problems as changing a solution (e.g., subset,
 structure).
 Fig. 41 depicts  the restructuring problem \cite{lev11restr}.

 Let \(P\) be a combinatorial optimization problem with a solution as
 structure
 \(S\)
 (i.e., subset, graph),
 \(\Omega\) be initial data (elements, element parameters, etc.),
 \(f(P)\) be objective function(s).
 Thus, \(S(\Omega)\) be a solution for initial data \(\Omega\),
 \(f(S(\Omega))\) be the corresponding objective function.
 Let \(\Omega^{1}\) be initial data at an initial stage,
  \(f(S(\Omega^{1}))\) be the corresponding objective function.
 \(\Omega^{2}\) be initial data at next stage,
  \(f(S(\Omega^{2}))\) be the corresponding objective function.

 As a result,
 the following solutions can be considered:
 ~(a) \( S^{1}=S(\Omega^{1})\) with \(f(S(\Omega^{1}))\) and
 ~(b) \( S^{2}=S(\Omega^{2})\) with \(f(S(\Omega^{2}))\).
 In addition it is reasonable to examine a cost of changing
 a solution into another one:~
 \( H(S^{\alpha} \rightarrow  S^{\beta})\).
 Let \(\rho ( S^{\alpha}, S^{\beta} )\) be a proximity between solutions
  \( S^{\alpha}\) and \( S^{\beta}\),
  for example,
 \(\rho ( S^{\alpha}, S^{\beta} ) = | f(S^{\alpha}) -  f(S^{\beta}) |\).
 Note, function \(f(S)\) is often a vector function.
 Finally, the restructuring problem can be examine as follows (a basic version):

~~

 Find a solution \( S^{*}\) while taking into account the
 following:

  (i) \( H(S^{1} \rightarrow  S^{*}) \rightarrow \min \),
  ~(ii) \(\rho ( S^{*}, S^{2} )  \rightarrow \min  \) ~(or constraint).

~~

 Thus, the basic optimization model can be examined as the following:
   \[\min \rho ( S^{*}, S^{2} ) ~~~s.t.
 ~ H(S^{1} \rightarrow  S^{*})  \leq \widehat{h}, \]
 where \(\widehat{h}\) is a constraint for cost of the solution
 change.

 Proximity function  ~\(\rho (S^{*},S^{2}) \)~
  can be considered as a vector function
 (analogically for the solution change cost).
 The situation will lead to a multicriteria restructuring problem
 (i.e., searching for a Pareto-efficient solutions).

\begin{center}
\begin{picture}(73,55)
\put(00,00){\makebox(0,0)[bl]{Fig. 41. Illustration for
restructuring problem \cite{lev11restr}}}

\put(00,05){\vector(0,1){46}} \put(00,05){\vector(1,0){70}}

\put(71,04.7){\makebox(0,0)[bl]{\(t\)}}

\put(00,51){\makebox(0,0)[bl]{``Quality''}}

%----------------S1

\put(0.6,15){\makebox(0,0)[bl]{\(S^{1}\)}}
\put(5,15){\circle{1.7}}

\put(12,13){\makebox(0,0)[bl]{Initial}}
\put(12,10){\makebox(0,0)[bl]{solution}}
\put(12,06){\makebox(0,0)[bl]{(\(t=\tau^{1}\))}}

\put(11,10){\line(-4,3){5}}

\put(6,16){\vector(1,1){23}}

%----------------S2

\put(40,35){\circle*{2.7}}

\put(54,37){\makebox(0,0)[bl]{Goal}}
\put(54,34){\makebox(0,0)[bl]{solution}}
\put(54,30){\makebox(0,0)[bl]{(\(t=\tau^{2}\)): \(S^{2}\)}}

\put(53,35){\line(-1,0){10.5}}

\put(40,35){\oval(12,10)} \put(40,35){\oval(17,17)}
\put(40,35){\oval(24,22)}

%++++++++++++++++++++++++++++++Proximity

\put(40,35){\vector(-2,1){9}} \put(30,40){\vector(2,-1){9}}

%\put(40,35){\line(-2,1){10}}

%--

\put(26,17){\makebox(0,0)[bl]{Proximity}}
\put(26,13){\makebox(0,0)[bl]{~\(\rho (S^{*},S^{2})\)}}

\put(36,20){\line(0,1){16}}

\put(46,19){\makebox(0,0)[bl]{Neighborhoods }}
\put(49,16){\makebox(0,0)[bl]{of ~\(S^{2}\)}}

\put(56,22){\line(-1,1){10}}

\put(53,22){\line(-2,1){8}}

\put(30,40){\circle{2}} \put(30,40){\circle*{1}}

\put(11,46){\makebox(0,0)[bl]{Obtained}}
\put(11,43){\makebox(0,0)[bl]{solution \(S^{*}\)}}

\put(20,42.8){\line(4,-1){7}}

%--------------------------------------
\put(1,37.5){\makebox(0,0)[bl]{Solution }}
\put(1,33.5){\makebox(0,0)[bl]{change cost }}
\put(1,29.5){\makebox(0,0)[bl]{\(H(S^{1} \rightarrow S^{*})\)}}

\put(8,29){\line(1,-1){5}}

%--------------------------------------

\end{picture}
\end{center}

\section{CONCLUSION}

 The paper contains a generalized integrated glance to design of
 system hierarchies.
 The problems of this kind are often crucial ones in system analysis, in systems design.
 The following  basic hierarchical structures are considered:
 trees,
 a special morphological hierarchy,
 multi-layer structures.
 Evidently, the examination
 is based on combinatorial optimization models.
 The list of design methods involves the following:
 expert-based 'top-down' procedure
 ('divisive' strategy of hierarchical clustering),
 hierarchical clustering (agglomerative algorithm),
 spanning trees,
 ontological approach,
 'optimal' organization,
 multi-layer structures (including multilevel network
 design problem).

 A special attention is targeted to modification problems.
 Here, a simplified modification of trees is examined
 (e.g., augmentation problem for graphs is not considered
  \cite{esw76}).
 The presented modification problems are as follows:
 modification of a weighted tree via condensing neighbor nodes
 (i.e., design of over-lay structure of tree-like software),
 hotlink assignment problem,
 building of a Steiner tree based on an initial tree,
 new class of restructuring problems.

 Generally, our material is the first integrated attempt
 to the problem.
 As a result,
 many considered issues require
 additional studies and analysis of many applied examples,
 many other significant problems have to be described and studied.
 For example, issues of uncertainty have not been examined.
 Mainly, considered problems are described on the basis of a
 simplified frames (e.g.,
  engineering description,
  problem formulation,
  basic solving schemes,
  prospective  versions).

 The presented material is oriented to needs of applications.
 In the future,
 it may be interesting to use the described approaches in education
  (engineering, management, computer science and
 information technology, applied mathematics).
 This material can be considered and used as a tutorial.

\end{document}